\def\real{{\tt I\kern-.2em{R}}}%Typos corrected through and including chapter 6 (the end of bui), 
\def\nat{{\tt I\kern-.2em{N}}}%Seems correct as of 2/27/10. One correction for in proof on last page and of 1/11/11. The arxiv file is corrected as of the date of the file 8/3/11. As of 11/28/11, addtional statement as to how the bold font is handed and the other qualifiers such as ``greater'' is mentioned. 14 OCT 2012
%contents change to include appendix. On 12/13/12, the foundations are changed by including $\cal W$ as part of the set of atoms where the symbols used are all considered as distinct from other stoms used
\def\eps{\epsilon}
\def\realp#1{{\tt I\kern-.2em{R}}^#1}
\def\natp#1{{\tt I\kern-.2em{N}}^#1}
\def\hyper#1{\ ^*\kern-.2em{#1}}

\def\hyperreal{{^*{\real}}}
\def\hyperrealp#1{{\tt ^*{I\kern-.2em{R}}}^#1} 
\def\hypernat{{^*{\nat }}}
\def\hypernatp#1{{{^*{{\tt I\kern-.2em{N}}}}}^#1} 
\def\eskip{\hskip.25em\relax}

\def\Hyper#1{\hyper {\eskip #1}}
\def\leaderfill{\leaders\hbox to 1em{\hss.\hss}\hfill}
\def\srealp#1{{\rm I\kern-.2em{R}}^#1}
\def\sp{\vert\vert\vert} 

\def\power#1{{{\cal P}(#1)}}
\def\iff{\leftrightarrow}
\def\qed{{\vrule height6pt width3pt depth2pt}\par\medskip}
\def\pars{\par\smallskip}
\def\parm{\par\medskip}
\def\m#1{{\rm #1}}
\def\b#1{{\bf #1}}
\def\r#1{{\rm #1}}
\magnification=\magstep0
\tolerance 10000
\hoffset=.64in
\hsize 5.00 true in
\vsize 8.50 true in
\baselineskip=14pt
\nopagenumbers
\font\large=cmdunh10 at 18truept
\font\next=cmr10 at 14truept
\pageno=1
\headline={\ifnum\pageno=1\hfil\quad \hfil \else\ifodd\pageno\rightheadline 
\else\leftheadline\fi\fi}
\def\rightheadline{\hfil \bf The Theory of Ultralogics \hfil\tenbf\folio}
\def\leftheadline{{\tenbf\folio} \hfil\bf Robert A. Herrmann \hfil}
\voffset=1\baselineskip
\medskip
{\quad}
\vskip 0.85in
\centerline{\bf \next The Theory of}\bigskip
\centerline{\large ULTRALOGICS}\bigskip
\centerline{Standard Edition}
\vskip 0.25in
\centerline{\next Robert A. Herrmann}
\vskip 4.75in
\centerline{Mathematics Department}
\centerline{U. S. Naval Academy}
\centerline{572C Holloway Rd}
\centerline{Annapolis$,$ MD 21402-5002}
\centerline{1993}
\centerline{(Part I latest revision$,$ 3 JUL 2012.)}
\vfil
\eject

\def\id{\par\hangindent2\parindent\textindent}
\def\textindent#1{\indent\llap{#1}}
 {\quad}
\vskip 1.5in
{
\centerline{First formal announcement of many of these results appeared in}\par\medskip  
\id (1) Mathematical philosophy$,$ Vol. 2 (No. 6) (1981)$,$ \#81T-03-529$,$ p. 527.
\id (2) A useful *-real valued function$,$ Vol. 4 (No. 4) (1983)$,$ \#83T-26-280$,$ p.  318. 
\id (3) Nonstandard consequence operators I$,$ Vol. 5 (No. 1) (1984)$,$ \#84T-03-61$,$ p. 129.
\id (4) Nonstandard consequence operators II$,$ Vol. 5 (No. 2) (1984)$,$ \#84T-03-93$,$ p. 195.
\id (5) D-world alphabets I$,$ Vol. 5 (No. 4) (1984)$,$ \#84T-03-320$,$ p. 269.
\id (6) D-world alphabets II$,$ Vol. 5 (No. 5) (1984)$,$ \#84T-03-374$,$ p. 328. 
\id (7) A solution to the grand unification problem$,$ Vol. 7 (No. 2) (1986)$,$ \#86T-85-41$,$ p. 238. 
\par
\medskip
\centerline{Some of the refereed papers relative to MA-model} 
\centerline{concepts and its mathematical construction.}
\medskip
\id (1)  A special isomorphism between 
superstructures$,$ Kobe J. Math.$,$ 10(2)(1993)$,$ 125-129. 
\id (2)  Fractals and ultrasmooth microeffects$,$ J. Math. Physics$,$ 30(4)$,$ April 1989$,$ 805-808. 
\id (3)  Physics is legislated by a cosmogony$,$ Speculations in Science and Technology$,$ 11(1) 
                        (1988)$,$ 17-24.           
\id (4) Nonstandard consequence operators$,$ Kobe  
                        J. Math.$,$ 4(1)(1987)$,$ 1-14. 
\id (5) Mathematical philosophy and developmental processes$,$ Nature and System$,$ 5(1/2)(1983)$,$ 
                        17-36.\par}
\id{A}nd many more.
\vfil
\eject
\baselineskip=16pt 
\centerline{}
\centerline{\bf CONTENTS}
\bigskip
\indent Chapter 1\pars
{\bf Intuitive Concepts}\pars
\line{\indent 1.1 \ \ The Alphabet$,$ Words$,$ and Choice Sets\leaderfill 7}
\line{\indent 1.2 \ \ Readable Sentences\leaderfill 9}
\line{\indent 1.3 \ \ Human Deduction\leaderfill 11}
\bigskip
\indent Chapter 2\pars
{\bf The G-Structure}\pars
\line{\indent 2.1 \ \ A Basic Construction\leaderfill 17}
\line{\indent 2.2 \ \ A Remark About 2.1\leaderfill 21}
\line{\indent 2.3 \ \ The Nonstandard Structure\leaderfill 22}
\line{\indent 2.4 \ \ General Interpretations\leaderfill 25}
\line{\indent 2.5 \ \ Sets of Behavior Patterns\leaderfill 26}
\bigskip                                                      
\indent Chapter 3\pars
{\bf Deductive Processes}\pars
\line{\indent 3.1 \ \ Introduction\leaderfill 29} 
\line{\indent 3.2 \ \ The Identity Process\leaderfill 29}
\line{\indent 3.3 \ \ Adjective Reasoning\leaderfill 31}
\line{\indent 3.4 \ \ Propositional Reasoning\leaderfill 34}
\line{\indent 3.5 \ \ Modus Ponens Reasoning\leaderfill 36}
\line{\indent 3.6 \ \ Predicate Deduction\leaderfill 36}
\bigskip
\indent Chapter 4\pars
{\bf Special Deductive Processes}\pars
\line{\indent 4.1 \ \ Introduction\leaderfill 39}
\line{\indent 4.2 \ \ Reasoning From the Perfect Type W\leaderfill 39}
\line{\indent 4.3 \ \ Strong Reasoning From the Perfect\leaderfill 40}
\line{\indent 4.4 \ \ Order\leaderfill 45}
\bigskip
\indent Chapter 5\pars
{\bf Consequence Operators}\par
\line{\indent 5.1 \ \ Basic Definitions\leaderfill 47}
\line{\indent 5.2 \ \ Basic $\sigma$ properties\leaderfill 48}
\line{\indent 5.3 \ \ Major Results\leaderfill 49}
\line{\indent 5.4 \ \ Applications\leaderfill 51}
\vfil
\eject
\centerline{}
\indent Chapter 6\pars
{\bf Associated Material}\pars
\line{\indent 6.1 \ \ Perception\leaderfill 55}
\line{\indent 6.2 \ \ Existence\leaderfill 55}
\line{\indent 6.3 \ \ An Alternate Approach\leaderfill 56}
\line{\indent References for Chapter 1 ---6\leaderfill 61} 
\bigskip
\indent Chapter 7\pars
{\bf Developmental Paradigms}\pars
\line{\indent 7.1 \ \ Introduction\leaderfill 62}
\line{\indent 7.2 \ \ Developmental Paradigms\leaderfill 64}
\line{\indent 7.3 \ \ Ultrawords\leaderfill 67}
\line{\indent 7.4 \ \ Ultracontinuous Deduction\leaderfill 71}
\line{\indent 7.5 \ \ Hypercontinuous Gluing\leaderfill 72}
\line{\indent References for Chapter 7\leaderfill 76}
\bigskip
\indent Chapter 8\pars
{\bf A Special Application}\pars
\line{\indent 8.1 \ \ A Neutron Altering Process\leaderfill 79}
\line{\indent References for Chapter 8\leaderfill 81}
\bigskip    
\indent Chapter 9\pars
{\bf NSP-World Alphabets}\pars
\line{\indent 9.1 \ \ An Extension\leaderfill 83}
\line{\indent 9.2 \ \ NSP-World Alphabets\leaderfill 83}
\line{\indent 9.3 \ \ General Paradigms\leaderfill 84}
\line{\indent 9.4 \ \ Interpretations\leaderfill 85}
\line{\indent 9.5 \ \ A Barrier to Knowledge\leaderfill 88}
\line{\indent References for Chapter 9\leaderfill 89}
\bigskip
\indent Chapter 10 \pars
{\bf Laws$,$ Rules$,$ and Other Things}\pars
\line{\indent 10.1 \ \ More About Ultrawords\leaderfill 91}
\line{\indent 10.2 \ \ Laws and Rules\leaderfill 94}
\line{\indent References for Chapter 10\leaderfill 99}
\vfil
\eject
\centerline{} 
\indent Chapter 11 \pars
{\bf ``Things''}\pars
\line{\indent 11.1 \ \ Subparticles.\leaderfill 101}
\line{\indent 11.2 \ \ Ultraenergetic Subparticles\leaderfill 109}
\line{\indent 11.3 \ \ More on Subparticles\leaderfill 112}
\line{\indent 11.4 \ \ MA-Model\leaderfill 116}
\line{\indent References for Chapter 11\leaderfill 117}
\line{\indent Appendix\leaderfill 119}
\line{\indent Symbols\leaderfill 129}
\vskip1.00in
\centerline{NOTES}
\noindent Some of the following 
typographical errors (or confusions) may (do) appear 
in this edition.\pars
\noindent [1] There may be a very few times when it is  not apparent as to 
the necessary structure for 
the  results stated. Of course$,$ one can always stay with the EGS. I may have 
written $\cal M$ when I meant ${\cal M}_1$ or conversely.\pars 
\noindent [2] The actual reason for the inverse order discussed on page 10 is 
that it is useful when considering adjective reasoning discussed later in 
this edition.\pars                                                                                
\noindent [3] Much of this edition was written in the late 1980s with the 
exception of Chapter 2$,$ sections 2.1 --- 2.2 which was written in 1991. In 
the formal model$,$ I 
wanted to differentiate between the actual set used for the words and the sets 
used for other analytical purposes. This is done by  letting members of 
${\cal W}$ symbols that are differfent than those in sets such as $\nat$ from Chapter 2 on.  This is the usual approach.  However, as of the date of this version, by a special technique the set of symbols $\cal W$ or extended symbols ${\cal W}'$ are now included.\pars
\noindent [4] I give two different superstructure constructions and certain processes used to obtain nonstandard models. It would have been better to concentrated on the second construction which is the one actually used.\pars
\noindent [5] Page 80. In order to have $B \cap M_i = \emptyset,\ i 
\not=0,$ we can reserve a special symbol in our alphabet. Then only consider a 
$B$ that does not contain this special symbol within any of its members. Now 
use this special symbol$,$ with or without spacing$,$ to construct the $M_i,\ i 
\not= 0.$ Of course$,$ once again we interpret  $\land$ in the axiom system for 
$S$ as this special symbol with or without spacing. \pars 
\noindent [6] Relative to the alphabet $\cal A.$ It is trivially obvious that 
one can include within this alphabet the written symbols used by an 
intelligent life form that uses a written language and deductive rules similar 
to those used by humankind. If one goes to the extreme and requires 
infinitely many such intelligent life forms$,$ then generalized languages as 
discussed on page 87 can be utilized.\pars 
\noindent [7] Beginning in Chapter IV$,$ due to the complexity of the first-order statements, I adopt the process of replacing set-theoretically defined predicates (abbreviations) such as $x \in y,\ x = y,\ 1\leq x < z$ etc. by $(x \in y),\ (x = y),\ (1\leq x < z).$ Of course, the $=$ is interpreted as set-theoretic equality. \pars
\vfil\eject
\baselineskip14pt   
\centerline{\bf 1. INTUITIVE CONCEPTS}\par
\bigskip
\leftline{\bf 1.1 The Alphabet$,$ Words$,$ and Choice Sets.}\par 
\medskip
There exists at the instant of time you read this sentence a finite set of all 
the {symbols you have previously used throughout your life} for your various 
forms of {communication and human deduction}. You may  also include frames from 
sound {motion picture film}$,$ {TV tape}$,$ and the like$,$ if you wish$,$ in the event 
you require {visual or audio stimuli} for your deductive processes. Let 
{$A_h$} be 
the set of such symbols for a given human being at this instant of time and 
{$H_t$} the set of all human beings who exist at this instant of time. Now let 
${\cal A}= \bigcup\{A_h\bigm \vert h\in H_t\}$ be the 
{\it alphabet}  for 
humanity at this instant of time. The set $\cal A$ may be enlarged to include 
all of 
the symbols which humanity has ever used$,$ if we wished. However$,$ for our 
purposes the set $\cal A$ will suffice. Observe that the set $\cal A$ is 
finite and among the numerous references relative to alphabets$,$ I direct your 
attention to [3] [7] [12] [13] [21]. \par\smallskip
Assume that there is included in the set $\cal A$ a distinct symbol which 
represents a blank symbol$,$ say something like $\sp$. Now following the above 
references$,$ any finite string (with repetition) of elements from $\cal A$ 
which is nonempty is called a {\it word}.  Let {$\cal W$}  be the 
set of 
all words. In certain applications of these mathematical methods$,$ a word is 
also called an {\it {intuitive}}   or {\it {naive readable 
sentence}}.
Only a fragment of the set $\cal W$ is used in this investigation and$,$ in most 
cases$,$ this fragment will be the set of {\it {meaningful sentences}} 
or a portion of a {\it {formal language}}.  
The metalanguage may be assume  
to be  written in a different color than is $\cal W.$ Also$,$ the concept of the  
{empty word} will NOT be employed.\par\smallskip
Next$,$ apply the concept that Markov calls the  {\it abstraction of} 
{\it identity} [12] to two words $W_1$ and $W_2.$ The two words are 
{\it equal} if they  
are composed of the 
same symbols in the same order written left to right or right to left etc.
The finitary character of each word allows for such an identification. 
Another intuitive string concept required for this discussion is the notion 
of the {\it juxtaposition operation} or 
{\it join} of two 
strings such as $W_1,\ W_2$ and this is denoted by $W_1W_2$ or $W_2W_1.$ 
Apparently$,$ since about 1914 [22]$,$ these two intuitive  string concepts 
relative to
{word theory} have been accepted as a reasonable consequence 
of the
finitary character of such forms.\par\smallskip
Following Robinson's procedure in [15] applied to $\cal W$ rather than to a formal 
language$,$ assume that there exists an injection $i$ from $\cal W$
onto $\nat,$ the set 
of natural numbers. Intuitively$,$ such an injection
exists since $\cal W$ is 
countably infinite. Now the term ``intuitive'' is  
utilized to denote the words$,$ 
various grammatical rules$,$ the \noindent informal logical procedures$,$ and the like used 
in ordinary communication. This is to differentiate  our common modes 
of descriptive communication from the formal theory into which these intuitive
objects are mapped.
Select a {fixed injection} for this
\break
\vfil
\centerline{} 
\eject
\noindent and all other 
investigations. Notice that Robinson used any set $U$ of cardinality \noindent greater
than or equal to that of his set of well-formed formulas. The procedure of 
{intuitively mapping} objects such as the set $\cal W$ onto concrete mathematical 
objects is well established and has been a major procedure in geometry  since 
{Descartes}. The {Cantor Axiom} used by {Hilbert and Birkhoff} [2] for modern 
geometry assumes that such as map exists from the set of points in a straight 
line onto the real numbers. Evidently$,$ this injection $i\colon {\cal W} \to \nat$ 
falls into a category of similar content as these well-known geometric 
assumptions. Indeed$,$ {$i$} can be an into G\"odel coding.\pars
Prior to continuing this introductory section a brief discussion of the 
Set Theory being used and related matters appears useful. The general 
set theory  
being 
used is {{\bf ZFC}} =  {{\bf ZF}} + {{\bf AC}}. The {\bf ZF} {(i.e. 
Zermelo-Fraenkel axioms)} and the Axiom of Choice {\bf AC} 
may be found listed on 
pages 2 ---19 of [5] among hundreds of other references. Within this general 
set theory we are working within a model for the axiom system {{\bf ZFH}} = 
{\bf ZFA} + {\bf AC} + $A$ {\bf is  countably infinite}. The axiom system 
{{\bf ZFA}} is the Zermelo-Franenkel axiom system with {atoms} (i.e. 
{urelements or individuals}). The set of atoms is the $A$ in the above formulation. 
All of the axioms are expressible in a first-order language with the 
predicates $\in$ and $=$$,$ where $\emptyset$ and $A$ are constants. The almost 
completely written axioms or modifications to {\bf ZF} + {\bf AC} axioms that 
yield {\bf ZFA} + {\bf AC} may be found on pages 122 of [5]$,$ page 44 of [6] 
and$,$ due to the use of individuals$,$ the actual system studied throughout [20]. 
I point out that the language used is a first-order language with the logical 
axioms for equality. \pars
Assuming the consistency of {\bf ZF}$,$ G\"odel showed that there is a model in 
which the {\bf ZFC} axioms hold. Thus the consistency of {\bf ZF} implies the 
consistency of {\bf ZFC}. Using {\bf ZFC} our model construction for 
{\bf ZFH} is a slight modification of that which appears in problem 1 on page 
51 of [6]. In the modification$,$ let $C$ be a countably infinite set of 
infinite subsets of the $\omega$ in the {\bf ZFC} model.  Moreover$,$ we bijectively 
map the order relation for $\omega$ onto $C - {\{a_0\}}.$ Thus$,$ for a bijection 
$f\colon (C - {\{a_0\}})\to \omega$ define the well-ordering $<$ on  $C - {\{a_0\}}$ 
as follows: For each $ x,\ y \in  (C - {\{a_0\}}),\ x < y$ iff 
$f(x) \in f(y).$ By the construction of the model$,$ this well-ordering is also 
a member of this model. We also have need of an interpretation of $=$ for the 
set $(C - {\{a_0\}}).$ This is to be the logical equality and as such is 
interpreted to be the identity relation on  $(C - {\{a_0\}})$ which also exists 
within this model. The set $(C - {\{a_0\}})$ with its order and equality relation
is to be considered the set of natural numbers $\nat$ within this model. On 
the other hand$,$ we also have within this model its $\omega$ that will be used 
when certain constructions are considered. The interpretation of $=$ for 
objects not in $(C - {\{a_0\}})$ is set-theoretic equality and the $\in$ in this 
model is the same as the $\in$ in our model for {\bf ZFC.} All of this yields 
a model for the {\bf ZFH} axioms that are$,$ hence$,$ consistent relative to 
{\bf ZF}.  Set $\cal W$ will be intuitively considered an object in this 
model. It does not contain the empty set 
since the empty word is not used. The injection $i$ is to be considered as a 
intuitive map from $\cal W$ onto $\nat$ within this model for {\bf ZFH.} [See page 19 detailed refinements.]
\pars
All of the new results are obtained by using informal mathematical reasoning 
relative to {\bf ZFH}$,$ proving results by means of the  {``observer 
language''} 
and using acceptable mathematical procedures. We subscribe to the remarks 
quoted in chapter III of the work by Rosiowa and Sikorski [14] as well as the 
observations made by Stoll [18$,$ p. 228] relative to Rosser's philosophy of 
mathematics. That is that the procedures used$,$ even though informal in nature$,$ 
are capable of formalization$,$ and such things as the  ``formal'' proofs of the 
formal consistency of {\bf ZFH} relative to {\bf ZF} by means of model theory 
methods imply that the same informal but acceptable procedures which we use,
when convincingly presented$,$ will {not produce any contradictions}.\pars
Now$,$ intuitively$,$ consider that there is a fixed dictionary {$D$} that uses a 
subset of an individual's personal $Y \subset {\cal A}.$ Let {$W_Y$} be the 
set of all 
words generated by $Y.$ Then descriptions of natural processes$,$ 
of their behaviors and developments$,$ as well as psychological descriptions for 
human behavior$,$ philosophical descriptions for belief-systems$,$ life styles and 
any other descriptions that are of interest are elements of the set 
${\cal P}(W_Y),$ the power set of 
$W_Y.$ Such sets are also informally called  {\it {describing sets}.} The 
meanings of such sets are understandable by the individual and have great 
content. Now there exists a set $B_Y \subset {\cal P}(W_Y)$ composed of all and 
only these describing sets. Evidently$,$ the set $B_Y$ is finite. Consider the 
set ${\cal P}(B_Y).$ In some applications of the forthcoming mathematical 
results$,$ the set ${\cal P}(B_Y)$ is called a {{\it free (will) choice set}.}
For the {{\it universal free (will) choice set}}$,$ simply consider the set of 
all ${\cal P}(B_Y)$ as $Y$ varies over all humanity which exists at a given 
instant of time. This is still a finite set. Notice that all that has been 
done to obtain these free choice sets can easily be formulated with respect to 
the set $\nat.$ Simply consider $i[W_Y],$ and then the formal describing sets are the 
corresponding elements of ${\cal P}(i[W_Y])$ under the injection $i.$ The 
intuitive describing sets can be recaptured by considering the inverse image 
of $i$ on these elements of ${\cal P}(i[W_Y]).$ The (general) Axiom of Choice 
is not necessary in order to obtain these free choice sets. However$,$ in the 
sequel$,$ The Axiom of Choice is used in the construction of the NSP-structure. 
Consequently$,$ The Axiom of Choice is utilized when we interpret$,$ in some 
applications$,$ the free choice sets as elements of the NSP-structure.
\bigskip
\leftline{\bf 1.2 Readable Sentences.}\par 
\medskip
For the purposes of this research$,$ not only is the finitary concepts of the 
abstraction of identity and join accepted but evidently a third fundamental 
procedure needs to be introduced and investigated. Consider the symbol string 
`mathematics'. Now this can be obtained or ``read'' by numerous applications 
of the join operation with symbol strings of lesser length. For example$,$ let 
$W_1 =$ {\bf math}; $W_2 =$ {\bf e}; $W_3 =$ {\bf mat}; $W_4 =$ {\bf ics}. 
Then {\bf mathematics} $= W_1W_2W_3W_4.$ Observe that $W_1,\ W_2,\ W_3,\ W_4$
are the syllables for this word. Clearly$,$ for writing purposes$,$ we could 
consider {\bf mathematics} $= W_1W_2 \ldots W_{11},$ where $W_j,\ j = 1, 
\ldots,11,$ are the single letters in the word. The necessity to  
consider intuitively a symbol string as composed of words of various length joined 
together from left to right leads$,$ when $i[{\cal W}]$ is considered$,$ to 
the concept 
of the set of (special) partial sequences.\par\smallskip
Let nonempty $H \subset i[{\cal W}]$ and $n \in \nat,$ where $0\in \nat.$ For 
simplicity$,$ let $H^n= H^{[0,n]}$ be the set of all maps from 
the segment 
$[0,n]$ into $H.$ An element of $H^n$ is called a {{\it partial sequence}} 
even though this definition is a slight restriction of the one that usually 
appears in the literature. Now let $P_H = \bigcup\{H^n\bigm \vert n\in \nat \}.$ 
In general$,$ if $H = i[{\cal W}],$ then the $H$ notation will 
be deleted from such symbols as $P_H.$ The set $P_H$ is a set of partial sequences and is a subset of 
$P_{i[{\cal W}]} = P.$ \par\smallskip
Let $P_1,\ P_2$ denote the set-theoretic first and second coordinate 
projection maps. Then the following first-order sentences$,$ where the usual 
assortment of set-theoretic abbreviations for subset$,$ functions$,$ domains$,$ 
ranges$,$ etc. are used$,$ hold and represent two basic properties for the object 
$P_H.$\par\medskip
\line{\hfil$\forall x((x\in P_H) \to \exists y((y \in \nat) \land (P_1[x] = [0,y]) \land
(P_2[x] \subset H)));$\hfil}\smallskip
\line{\hfil$\forall x((x \subset \nat \times H) \land \exists y((y \in \nat) \land (P_1[x]= 
[0,y]\land$\hfil}\smallskip
\line{(1.2.1)\hfil$\forall w\forall z \forall z_1(((w\in \nat) \land
(z \in H) \land (z_1 \in H) \land ((w,z) \in x) \land$ \hfil }
\line{\hfil $((w,z_1)\in x))\to z = 
z_1)\to x \in P_H).$\hfil} \par\medskip
I point out that ever since {G\"odel used a natural number coding} for certain 
metamathematical concepts$,$ interpreting naive or intuitive processes involving 
symbol strings as concepts relative to $\nat$ has prevailed and become 
accepted by mathematicians. As Kleene writes: ``{\tt Metamathematics has 
become a branch of number theory.''} [7$,$ p. 205]  Consequently$,$ it is clearly 
justified to {use partial sequences} to discuss the ``ordering'' of symbol 
strings. This ordering will be associated with the finitary ordering obtained 
by {joining words by juxtapositioning} them in a specific intuitive order.\par\smallskip
Consider $f \in H^n,   n  \in \nat.$ Then $f(0),f(1),\ldots,f(n) \in H.$ The 
{{\it order induced by f}} is defined to be the simple inverse order 
generated by $f$ applied to the simple order of $[0,n].$ Formally$,$ for each 
$f(j),\ f(k) \in f[[0,n]],$ define the order induced by $f$ to be 
$f(k) \leq_f f(j)$ iff $j\leq k,$ where as usual the order $\leq $ is the 
simple order on $[0,n]$ induced by the simple order on $\nat.$ In general$,$ the 
notation {$\leq_f$} will NOT be specifically used to denote this $f$ induced 
order but$,$ rather$,$ the {order will be indicated} by writing the symbols 
$f(n),\ldots,f(0)$ from left to right in an ordered fashion. This corresponds 
to what the intuitive ordering would be under the inverse of $i$$,$ $i^{-1},$ 
when it is used to recapture the original symbol string.  Remember$,$ that 
everything done with the join operation$,$ partial sequences and the induced 
order is {finitary in character}. Apparently$,$ care is required in the selection 
of the $H$ if all the ways a given word may be partitioned for readability are 
to be investigated. It is required that $i^{-1}[H]$ contain enough symbols for 
this purpose. However$,$ an approach is now developed that eliminates this 
apparent difficulty.\par\smallskip
Relate the induced ordering for some $f \in P$ to the join in the 
following manner. The word $W_0$ to which $f\in (i[{\cal W}])^n,\ n \in \nat,$ 
corresponds is
$$W_0 = (i^{-1}((f(n)))(i^{-1}((f(n-1)))\cdots (i^{-1}((f(0))). \leqno 
(1.2.2)$$
Using the abstraction of identity and this ordered join concept$,$ a 
{basic equivalence  relation} is obvious. For $f,\ g\in P\ f\in H^n,\ g\in 
H^m,$ let \par\medskip
\line{\hfil$f \sim g\ {\rm iff}$\hfil}\smallskip
\line{\hfil$(i^{-1}((f(n)))\cdots(i^{-1}((f(0)))=(i^{-
1}((g(m)))\cdots(i^{-1}(((0))).$\hfil (1.2.3)} \par \medskip
Observe that this process is still finitary in nature and consequently 
effectively recognizable. Clearly$,$  {``$\sim$''} is an equivalence relation 
on $P$ since the identity  $=$ is such a relation on ${\cal W}.$ For each $f \in P$ let 
{$[f]$} denote the corresponding {{\it basic equivalence class.}} 
Now if the 
{{\it cardinality of an intuitive word}} $\vert W_0 \vert = m +1 =$ the total number of symbols in 
$W_0$ including repetitions and $W_0= (i^{-1}((f(m)))\cdots(i^{-1}((f(0))),$ 
then each equivalence class is a nonempty 
{\bf finite} set and theoretically each element in each equivalence class can 
be effectively recognized.\par\smallskip
Recall that many of the intuitive concepts associated with word theory and 
algorithms [7] [12] [13] [22] have much less rigorously defined concepts than 
this equivalence relation even though such word theory concepts have been 
extensively employed. Before proceeding$,$ however$,$ here are some of the simple 
facts about these equivalence classes. Consider $[f]$ and let the 
corresponding word for $f$ be $W_0,$ where $\vert W_0 \vert = m + 1.$ Then 
there exists two unique maps $f_m,\ f_0 \in P$ such that 
$f_m \sim f_0,\ f_m \in (i[{\cal W}])^m,\ f_0 \in (i[{\cal W}])^0$ and 
$$W_0 = (i^{-1}((f_m(m)))\cdots(i^{-1}((f_m(0)))= (i^{-1}(f_0(0))).\eqno 
(1.2.4)$$
Furthermore$,$ for each $k \in \nat$ such that $0\leq k \leq m,$ there exists 
$g \in [f]$ such that $g \in (i[{\cal W}])^k.$ And$,$ for each $k \in \nat$ such 
that $k>m$ there does not exist any $g \in (i[{\cal W}])^k$ such that 
$g \in [f].$ Finally$,$ how is a 
particular class $[f]$ to be intuitively interpreted? In order to interpret a class $[f]$ in ${\cal W},$ simply select 
any element in $[f],$ say $f_n,$ and effectively construct the word
$(i^{-1}((f_n(n)))\cdots (i^{-1}((f_n(0)))=W_0.$ The word $W_0$ is called the
{{\it intuitive or naive interpretation}} for the class $[f].$ \par\smallskip
Let {${\cal E}$} be {{\it set of all equivalence classes generated by
$\sim$ on $P.$}} The set ${\cal E}$ is called {{\it the set of (formal)
readable sentences.}} The term ``readable sentence'' is used in two contexts$,$ 
the intuitive readable sentence that is a member of ${\cal W}$ and the 
corresponding formal readable sentence in $\cal E.$ The terms  ``intuitive'' or 
``formal'' 
will not be used where no confusing would occur.\par
\bigskip
\leftline{\bf 1.3 Human Deduction.}\par 
\medskip
As has become the custom$,$ the concepts of human deduction (i.e. reasoning ) 
will first be discussed intuitively with respect to the set $\cal W.$ Certain
metatheorems relative to such processes will be established prior to 
associating these processes with the more formal set $\cal E.$ Recall
Tarski's [21] basic axioms for the undefined {(finite) consequence 
operators} {$C$}   
on a set of meaningful sentences $A.$  As usual$,$ let {$\vert A\vert=$}
the cardinality of $A,$ the symbol
{$\power A$ denote the powerset of $A,$} {$F(A)$ denote the set of all 
finite subsets of $A$}  and {$C\colon \power A \to \power A.$} Tarski 
first bounds the cardinality of $A$
$$({\rm 1})\ 0 < \vert A \vert \leq \aleph_0$$
Now the map $C$ is a {(finite) {\it consequence operator}} if $A \not= 
\emptyset$ and 
%===SIMPLE TABLE OF ONE BOX WITH MULTIPLE LINES 
%IT IS IN DISPLAYED FORM WHEN PRINTED
$$\vbox{\offinterlineskip
\hrule
\halign{&\vrule#&
\strut\quad#\hfil\quad\cr
height2pt&\omit&\cr
&${\rm (2)\  for \ each}\ B \subset A,\ {\rm then}\ B \subset C(B) \subset 
A,$&\cr
&\ &\cr
&${\rm (3)\  for \ each}\ B \subset A,\ {\rm then}\ C(C(B))=C(B),$&\cr
&\ &\cr
&${\rm (4)\  for \ each}\ B \subset A,\ {\rm then}\ C(B) = \bigcup \{C(F)\bigm 
\vert F \in F(B)\}$&\cr                
height2pt&\omit&\cr}
\hrule}$$     \par\smallskip
Axioms (2)$,$ (3)$,$ (4) appear to be considerably more significant than does 
axiom (1) in these applications and the cardinality of the set of 
sentences considered is not$,$ usually$,$ so bounded. Actually the consequence 
operator is generated by a slightly more general concept which is called the 
deductive process. Let $A$ be a nonempty set. Then any nonempty $k 
\subset F(A) \times A$ is a {{\it deductive process}.} The term {{\it 
deductive operator}} is also used. A deductive process $ k \subset F(A)\times 
A$ is 
{{\it total}} if for each nonempty $B \subset A$ there exists some $F\in 
F(B)$ and some $a \in A$ such that $(F,a) \in k.$ For a total $k \subset 
F(A)\times A,$ let \par\medskip
\line{\hfil(i) $C_k= \{(x,y)\bigm \vert x \in \power A\land 
y \in \power A\land y \not= \emptyset \land $\hfil}\smallskip
\line{(1.3.1)\hfil$\forall w(w \in A\land  w \in y \iff 
 \exists z(z \in F(x)\land (z,w) \in k))\}$\hfil }\smallskip
\line{\hfil And$,$  $(B,\emptyset)\in C_k \iff B = \emptyset \iff \emptyset 
\notin P_1[k].$\hfil}\par\medskip
\noindent Notice that the definition of $C_k$ implies 
that $C_k \colon \power A \to \power A.$ \par\smallskip
Assume that $A$ is a nonempty arbitrary set that corresponds to Tarski's set 
of
meaningful sentences with the exception  that axiom (1) need not hold for $A.$ 
\par \medskip
{\bf Theorem 1.3.1.} { \sl Let total $k \subset F(A) \times A.$ Then
$C_k$ satisfies Tarski's axiom (4).}\par\smallskip
Proof. Obviously$,$ if $B\subset A,$ then $C_k(B) \subset A.$ Let $C_k(B) 
=\emptyset.$ Then $B = \emptyset$ implies that $C_k(B) = \bigcup\{\emptyset 
\}= \bigcup\{C_k(F) \mid F \in F(\emptyset)\} = \bigcup\{C_k(\emptyset)\}.$  
Consequently$,$ assume that $C_k(B )\not= \emptyset$ and $x \in C_k(B).$ 
Then there 
exists some $F_0 \in F(B)$ such that $F_0 \subset B \subset A$ and $(F_0,x) 
\in k.$ From the definition of $C_k,$ this implies that $x \in C_k(F_0).$
Thus $C_k(B) \subset \bigcup\{C_k(F)\bigm\vert F \in F(B)\}.$ On the other 
hand$,$ assume that $x \in \bigcup\{C_k(F)\bigm\vert F \in F(B)\}.$ Then 
there exists some $F_1 \in F(B)$ such that $x  \in C_k(F_1).$ Hence there 
exists some $F_2 \in F(F_1) \subset F(B)$ and $(F_2,x) \in k.$ Consequently$,$ 
$x \in C_k(B)$ and this completes the proof. \qed
Many of the usual deduction processes$,$ such as propositional deduction$,$ 
satisfy the following additional properties. Let $k \subset F(A) \times A$ be 
{{\it ordinary}} whenever\par\smallskip
(i) if $(A_0,b) \in k$ and $A_0 \subset D\in F(A),$ then $(D,b) \in k$ 
and\par\smallskip
(ii) if $\{(A_0,b_1),\ldots,(A_0,b_n)\} \in k$ and 
$(\{b_1,\cdots,b_n\},c) \in k,$ then $(A_0,c) \in k.$\par\medskip
{\bf Theorem 1.3.2.} {\sl Let total and ordinary $k \subset F(A) \times A.$ If 
$C_k$ satisfies axiom (2) of Tarski$,$ then $C_k$ is a consequence 
operator.}\par\smallskip
Proof.
Since axiom (4) holds for $C_k,$ now proceed to establish axiom (3). Tarski 
[21] has shown that axiom (4) implies that if $B \subset D \subset A,$ then 
$C_k(A) \subset C_k(D).$ Thus assume $B \subset A.$ Then $B\subset C_k(B)$ 
implies that $C_k(B) \subset C_k(C_k(B)).$ Now if $C_k(C_k(B))=\emptyset,$ 
then this implies that $C_k(B) = \emptyset$ and $C_k(C_k(B))=C_k(B).$
Let $x \in C_k(C_k(B)).$ Then  there exists some finite $F_0 \subset C_k(B)$ 
such that $(F_0,x) \in k.$ If $F_0 = \emptyset,$ then $\emptyset \subset
B$ implies $x \in C_k(F_0) \subset C_k(B).$ Assume that $F_0 \not= \emptyset,$ 
say $F_0 = \{a_1,\ldots,a_n\}.$ Then there exists finitely many 
$F_i \subset B,\ i = 1,\ldots,n$ such that $(F_i,a_i) \in k.$ Let $F^\prime
= F_1 \cup \cdots \cup F_n \in F(B) \subset F(A).$ Then ordinary (i) implies 
that $(F^\prime,a_i) \in k,\ i = 1,\ldots,n.$ But $(F_0,x) \in k$ and ordinary
part (ii) implies that $(F^\prime,x) \in k.$ Consequently$,$ $x \in C_k(B)$ and 
this completes the proof. \qed
In passing note that if axiom (2) holds for $C_k$ and (i) of ordinary holds 
for $k,$ then for each $x \in A,\ (\{x\},x) \in k.$ Obviously$,$ 
if for each $x \in A,\ (\{x\},x) \in k,$ then axiom (2) holds for $C_k.$ In 
this case$,$ we say that the deduction process is {{\it singular}}. Combining 
the above results$,$ we have the next theorem.\par\medskip
{\bf Theorem 1.3.3.} {\sl If for nonempty $A$ the deductive process 
$k \subset F(A) \times A$ is a total$,$ ordinary and singular$,$ 
then $C_k$ is a consequence operator for $A.$}\par\medskip 
Is there an obvious deduction process generated by a given consequence 
operator? Let $C$ be a consequence operator on $\power A.$ Define
$k_c \subset F(A) \times A$ as follows: 
$$(F,a) \in k_c\ {\rm iff}\ F\in F(A) \ {\rm and} \  a\in C(F). \eqno 
(1.3.2)$$
By axiom (2)$,$ it follows that $k_c$ is total.\par \medskip
{\bf Theorem 1.3.4.} { \sl If $C$ is a consequence operator on $A,$  then
$C_{k_c} = C.$}\par\smallskip
Proof. First assume that $C(B) = \emptyset.$ Then axiom (2) implies that
$B = \emptyset.$ Hence $\emptyset \notin P_1(k_c).$ From this it follows that
$C_{k_c}(\emptyset) = \emptyset =C(B).$ Now let $B \subset A$ and suppose that 
$x \in C(B).$ Then there exists some $F\in F(B)$ such that $x \in C(F).$ Since 
$F\in F(A)$$,$ then $(F,x) \in k_c.$ From the definition of $C_{k_c},$ it follows 
that $x \in C_{k_c}(B);$ which yields that $C(B) \subset C_{k_c}.$ Now if
$C_{k_c}(B) = \emptyset,$ then  $C_{k_c}(B) \subset C(B).$ Hence suppose that
$y \in C_{k_c}(B).$ Then there exists $F \in F(B)$ such that $(F,y) \in k_c.$ 
Since $F \in F(A),$ then $y \in C(F)$ implies that  $C_{k_c}(B)\subset C(B).$ 
Therefore$,$ for each $B \in \power A ,\  C_{k_c} (B) = C(B)$ implies that
$ C_{k_c} = C$ and the proof is complete. \qed
From the above results$,$ a consequence operator can be thought of as being 
determined by a deductive process and$,$ in certain cases$,$ conversely. When 
deductive processes are described in a metalanguage$,$ then such properties  as
total$,$ ordinary or singular can often be easily established. Sometimes we say 
that a deductive process $ k\subset F(A) \times A$ satisfies the Tarski axioms 
(2)$,$ (3)$,$ or (4) if $C_k$ satisfies (2)$,$ (3)$,$ or (4)$,$ 
respectively.\par\smallskip
Our next task is to find an appropriate correspondence between any $k \subset
F({\cal W}) \times {\cal W}$ and some $K  \subset F({\cal E}) \times {\cal 
E}$ such that the axioms of Tarski and set-theoretic properties are preserved. 
It is clear that this relation should be defined relative to the quotient map 
determined by $\sim.$ Thus for each $w  \in {\cal W},$ let $f_w \in P$ be such 
that $f_w \in (i[{\cal W}])^0$ and $f_w(0) = i(w).$ Consider the intuitive 
bijection $\Theta \colon {\cal W} \to {\cal E},$ the {{\it basic 
bijection}$,$} defined by $\Theta (w) = 
[f_w]$ for each $w  \in {\cal W}.$ In the usual manner$,$ extend $\theta$ to its 
corresponding set functions and the like. Since it is not assumed that all 
readers of this book are aware of these definitions$,$ we present them in the 
context of consequence operator theory. However$,$ for the next few results$,$ the 
map will not be restricted to the specially defined map $\Theta$ but rather we 
establish them for any arbitrary injection $\beta \colon A \to X.$ Recall 
that\par\medskip 
(i) if $B \subset A,$ then $\beta [B] = \{\beta (x)\mid x \in B \},$ 
\par\smallskip
(ii) for
${\cal B} \subset \power A ,\ \beta [{\cal B}] = \{ \beta [x] \mid x \in {\cal 
B}\},$ \par\smallskip
(iii) for $k \subset F(A) \times A,$ let $\beta k = \{(\beta [x],\beta 
(y))\mid (x,y) \in k\} = \{(z,w)\mid z \in F(x) \land w \in X \land 
\exists x \exists y (x \in F(A) \land y \in A \land z = \beta [x] \land
w = \beta (y) \land (x,y) \in k) \},$\par\smallskip 
(iv) for $C\colon \power A \to \power A,$ let $\beta C = \{(\beta [x],\beta 
[y])\mid (x,y) \in C\}.$ \par\medskip
\noindent By considering the inverse$,$ $\beta^{-1},$ it is evident that if $k$ is total 
on $A,$ then $\beta k$ is total on $\beta [A];$ if $k$ is ordinary on $A,$ 
then $\beta k$ is ordinary on $\beta [A];$ and if $k$ is singular on $A,$ then $\beta 
k$ is singular on $\beta [A].$\par\smallskip
Numerous propositions are immediate consequences of the set-theoretic 
definitions associated with the injection $\beta.$ However$,$ even though the 
following propositions hold true from elementary algebraic 
results$,$ we prove them explicitly for they refer to intuitive structures and 
metatheoretric results. Furthermore$,$ they are of considerable importance to 
the foundations of this subject.\par \smallskip
As before$,$ let $\beta$ be an injection on $A$ into $X,$ let $k \subset F(A) 
\times A$ and $C\colon \power A \to \power A.$ \par\medskip
{\bf Theorem 1.3.5.} {\sl Let $B \subset \beta [A].$ Then\par\medskip 
\line{\hfil$\beta^{-1}[(\beta C)(B)] = C(\beta^{-1}[B]),\ {\rm and}$\hfil}\par
\medskip
\line{{\rm (1.3.3)}\hfil$(\beta C)(\beta[\beta^{-1}[B]]) = (\beta C)(B) =
 \beta[C(\beta^{-1}[B])].$\hfil }}\par
\medskip
Proof. Let $B \subset \beta [A].$ Then $\beta^{-1}[\beta [B]] \in \power 
A.$ Hence$,$ $(\beta^{-1}[B],C(\beta^{-1}[B])) \in C.$ Consequently,
$$(\beta[\beta^{-1}[B]],\beta[C(\beta^{-1}[B])]) \in \beta C\leqno (1.3.4)$$
implies since $\beta C$ is a map that
$$\beta[C(\beta^{-1}[B])] = (\beta C)(\beta[\beta^{-1}[B]]) = (\beta 
C)(B).\leqno (1.3.5)$$
Thus 
$$\beta^{-1} [\beta[C(\beta^{-1}[B])]] = \beta^{-1}[(\beta C)(B)] = C(\beta^{-
1}[B])\leqno (1.3.6)$$
and this completes the proof. \qed
The next two results are important consequences of Theorem 
1.3.5.\par\medskip{\bf Theorem 1.3.6.} {\sl If $C$ is a consequence operator on 
$\power A,$ (i.e. it satisfies axioms (2)$,$ (3)$,$ (4) of Tarski)$,$ then 
$\beta C$ is a consequence operator on $\power {\beta [A]}.$}\par
\smallskip
Proof. Observe that $\beta [\power A] = \power {\beta [A]}$ and that
$\beta C\colon \power {\beta [A]} \to \power {\beta [A]}.$\par \medskip
(i) Let $B\subset \beta [A].$ Then $\beta^{-1} [B] \subset A$ implies that 
$\beta^{-1} [B] \subset C(\beta^{-1} [B]) \subset A.$ Hence $B \subset 
\beta[C(\beta^{-1}[B])] = (\beta C)(B) \subset \beta [A].$ \par\medskip
(ii) Let $B \subset \beta [A].$ Then $\beta^{-1}[B] \subset A$ and 
$C(C(\beta^{-1}[B])) = C(\beta^{-1}[B]).$ Consequently\par\medskip
\line{\hfil$\beta[C(C(\beta^{-1}[B]))] = (\beta C)(\beta [C(\beta^{-
1}[B])])=$\hfil}\pars
\line{\qquad\qquad\hfil$(\beta C)(\beta C)(\beta [\beta^{-1}[B]]) = (\beta C)(\beta C) (B) 
=$\hfil (1.3.7)}\pars
\line{\hfil$\beta [C(\beta^{-1}[B])] = (\beta C)(\beta [\beta^{-1}[B]])=
(\beta C)(B).$\hfil}\par\medskip
(iii) First$,$ we show that $\beta[F(\beta^{-1}[B])] = f[B].$ Let $F \in 
F(\beta^{-1} [B]).$ Then $\beta[F] \subset B$ and $\beta[F] \in F(B).$ 
Therefore$,$ $\beta[F(\beta^{-1}[B])] \subset F(B).$ Conversely$,$ let $F \in 
f(B).$ Then $\beta^{-1} [F] \in F(\beta^{-1}[B]),$ since $\vert F \vert =
\vert \beta^{-1}[F] \vert.$ Thus $ \beta [F(\beta^{-1}[B])] = F(B).$\pars
For each $B \subset \beta [A],\ C(\beta^{-1}[B]) = \bigcup \{C(F) \mid 
F \in F(\beta^{-1} [B])\}.$ This implies that\parm
\line{\hfil $\beta [C(\beta^{-1}[B])] = \bigcup\{\beta [C(F)] \mid
F \in F(\beta ^{-1} [B])\} =$\hfil}\smallskip
\line{\quad\qquad\hfil $\bigcup \{(\beta C)(\beta [F])\mid F \in F(\beta^{-1}[B])\} 
=$\hfil (1.3.8)}\smallskip
\line{\hfil$\bigcup \{(\beta C)(F) \mid F\in F(B) \} = (\beta 
C)(B).$\hfil}\parm 
\noindent Results (i)$,$ (ii)$,$ (iii) imply that $(\beta C)$ is a consequence 
operator and the proof is complete. \qed
{\bf Theorem 1.3.7.} {\sl If $k \subset F(A) \times A$ is total$,$ then 
$\beta (C_k) = C_{\beta k}.$}\par \smallskip
Proof. Let $(x,y) \in \beta (C_k).$ Then $(\beta^{-1}[x],\beta^{-1}[y]) \in 
C_k.$ Notice that $\beta^{-1}[y] = \emptyset $ iff $y = \emptyset.$ Assume 
that $z \in \beta^{-1}[y].$ Then there exists some $ w_z \in F(\beta^{-
1}[x])$ such that $(w_z,z)\in k.$ Hence$,$ it follows that for each $\beta (z) 
\in y$ there exists some $\beta [w_z] \in F(x)$ such that $(\beta [w_z],\beta 
(z)) \in \beta k.$ This leads to $(x,y) \in C_{\beta k}.$ Now if $ y = 
\emptyset,$ then $\beta^{-1}[y] = \emptyset$ implies that $\beta^{-1}[x] = 
\emptyset$ and $ x = \emptyset.$ Hence$,$ $(\emptyset,\emptyset) \in \beta 
(C_k)$ implies that $(\emptyset,\emptyset) \in C_k.$ The definition gives 
$\emptyset \notin P_1(k)$ and $\beta [\emptyset] = \emptyset \notin P_1(\beta 
k).$ Thus $(\emptyset,\emptyset) \in C_{\beta k}.$\pars 
On the other hand$,$ for 
each $z \in y \not= \emptyset$ there exists some $F \in F(x)$ such that $(F,z) 
\in \beta k.$ Hence $(\beta^{-1}[F],\beta^{-1}(z)) \in k$ and again $ z \in y$ 
iff $\beta^{-1}(z) \in \beta^{-1}[y];\ \beta ^{-1}[F] \in F(\beta^{-1}[x])$ 
imply that $(\beta^{-1}[x],\beta^{-1}[y]) \in C_k.$ Also$,$ if $y = \emptyset,$ 
then $x = \emptyset$ (i.e. $(\emptyset,\emptyset) \in C_{\beta k}$) and 
$\emptyset \notin P_1(\beta k).$ Hence $\emptyset \notin P_1(k)$ and 
$(\emptyset,\emptyset) \in \beta (C_k)$ and the proof is complete. \qed
I will not continue with this piecemeal approach but rather use a more general 
result established within the next chapter$,$ where $\cal E$ will be defined 
on a set $\cal W$ 
${\cal W}\cap \nat = \emptyset.$ \par
\vfil
\eject
\centerline{NOTES}
\vfil\eject
 
\centerline{\bf 2. THE G-STRUCTURE}\par
\bigskip
\leftline{\bf 2.1 A Basic Construction.}\par 
\medskip
 A primary construction will be a {{\it superstructure}}. For 
$X,$ a superstructure is constructed as follows: Let {\it ground set} $X = X_0$. 
Then by 
induction$,$ let $X_{n+1} = X_n \cup \power {X_n}.$ Now let ${\cal X} =
\bigcup \{X_n \mid n \in \omega \}.$ The set $\cal X$ is a {{\it 
superstructure over $X_0$.}} Within our 
model for {\bf ZFA} 
a set $B$ is {$X_0$-{\it transitive}} if for each $ x \in B$ 
either $ x \in X_0$ or $x \subset B.$ \pars
{\bf Theorem 2.1.1} {\sl For each $ n \in \omega,$ the set $X_n$ is 
$X_0$-transitive.}\pars
Proof. The proof is by induction. Let $ n = 0.$ The $X_0$ is $X_0$-transitive
immediately from the definition. Assume that for $ n,$ the set 
$X_n$ is
$X_0$-transitive. Consider $X_{n+1}$ in the above construction. We need only 
check any $x \in X_{n+1} - X_0 = (X_n \cup \power {X_n}) - X_0.$ Hence$,$ either 
$x \in X_n - X_0$ or $x \in \power {X_n} -X_0.$ In the first case$,$ $ x \subset 
X_n$ by the induction hypothesis. In the second case$,$ $x \subset X_n$ by the 
definition of the power set operator. Since $X_n \subset X_{n+1},$ it follows 
that $x \subset X_{n+1}.$ Thus by induction the proof is complete. \qed
{\bf Theorem 2.1.2} {\sl For each $n \in \omega$$,$ if $y \in x \in X_{n+1} - X_0,$  
then $ y \in X_{n}.$}\pars
Proof. For $n =0,\ y\in x\subset X_0 \Rightarrow y\in X_0.$ Assume that it 
holds for  
$n- 1,\ n \geq 1.$ Let $y \in x \in X_{n+1} - X_0.$ Then 
$x \in (X_n \cup \power {X_n}) - X_0.$ If $ x \in X_n - X_0,$ then by the 
induction hypothesis $y \in X_{n-1}.$ But $X_{n-1} \subset X_n$ implies that 
$y \in X_n.$ If $x \in \power {X_n},$ then $x \subset X_n$ implies that $y \in 
X_n$. By induction the proof is complete. \qed
Obviously$,$ since we have only used facts about {\bf ZF} to establish Theorems 
2.1.1$,$ 2.1.2$,$ these theorems hold for superstructures within {\bf ZF}.  
Recall that if $A$ is a set of atoms$,$ then this means that 
if $x \in A,$ then $x \not= \emptyset$ and  $y \in x$ is not defined. 
A nonempty  ground set $X_0$ is {{\it  $n$-atomic}} 
if $x \in X_0$ implies that $x \not= \emptyset$ and if $ y \in x \in X_0,$ then 
$ y \notin X_n.$ Two important observations relative to $n$-atomic. If $X_0$ is 
a set of atoms$,$ then $X_0$ is $n$-atomic for each $ n \in  \omega.$ If 
$X_0$ is $n$-atomic$,$ then $X_0$ is $k$-atomic for each $k$ such that $0\leq 
k\leq n.$   For each $X_n,\ n \geq 0,$ let
$M_{X_n} = \{(x,y) \mid (x \in X_n) \land (y \in X_n) \land (x \in y)\}$ and 
$E_{X_n}=\{(x,y) \mid (x \in X_n) \land (y \in X_n) \land (x =y)\},$ where the = is 
set-theoretic equality on sets and the identity on atoms. In like manner$,$ 
for ground set $Z,$ defined $M_{Z_n}$ and $E_{Z_n}$ for the respective $Z_n.$ For $n\geq 1,$ an 
isomorphism $\beta_n$ from $\langle X_n,M_{X_n},E_{X_n}, \nat,\emptyset\rangle$ 
onto  $\langle Z_n,M_{Z_n},E_{Z_n},\nat,\emptyset\rangle$ is {{\it special}} if 
$\beta_n (X_{k}) = Z_{k},\ 0\leq k\leq n-1,$ where $\nat$ is a set of atoms. 
Observe that since $X_k \in X_{k+1} 
\subset X_n,$ it follows that $X_k \in X_n.$ \parm{\bf Theorem 2.1.3} {\sl  Let $A$ be a set of atoms. Suppose that 
for nonempty sets $X,\ Z,\ X\cap A = Z \cap A$  
and there exists a bijection $\beta\colon X \to Z= \beta [X],$ where $\beta$ 
is a set-theoretic bijection on sets and the identity on any atoms in $X 
\cap A.$ Consider the sets $X_0 = X \cup A,\ Z_0 = Z \cup A,$ and 
$A$ and $\emptyset$ as the constants that denote a set of atoms and the 
empty set in our {\bf ZFA} model.\parm
\centerline{}
\eject
(i)  If $X_0,\ Z_0$ are 0-atomic$,$ then the structures 
$\langle X_0,M_{X_0},E_{X_0}\rangle$ 
and $\langle Z_0,M_{Z_0},E_{Z_0}\rangle$ are isomorphic.\pars
(ii) For each $n\geq 1,$ if $X_0,\ Z_0 $ are 
$n$-atomic$,$ then the structures $\langle X_n,M_{X_n},E_{X_n}, A,
\emptyset\rangle$ 
and $\langle Z_n,M_{Z_n},E_{Z_n},A,\emptyset\rangle$ are isomorphic 
and the isomorphism is special.}\pars
Proof. By $\in$-recursion$,$ define the map $\zeta$ on {\cal X} as follows:\par
{\leftskip 1.5in \noindent For $ x \in X,\ \zeta(x) = \beta(x);$\par}\pars
{\leftskip 1.5in \noindent For $x \in X_0 - X,\ \zeta (x) = x;$\par}\pars
{\leftskip 1.5in \noindent For $x \in {\cal X} - X_0,\ \zeta(x) = 
\zeta[x].$\par}\pars
\noindent Let $ \beta_n = \zeta\vert X_n,$ where $n \in \omega.$ We need 
only show that for each $n \in \omega,$ if $X_0$ and $Z_0$ are $n$-atomic$,$ 
then
$${\rm (A)}\ \beta_n\ {\rm is}\ {\rm an}\ {\rm isomorphism}\ {\rm from}\ 
\langle X_n, M_{X_n}, E_{X_n}\rangle$$ $${\rm onto}\  \langle Z_n, 
M_{Z_n},E_{Z_n}\rangle;$$
$${\rm (B)}\ {\rm if}\ n \geq 1,\ {\rm then}\ \beta_n(A) = A,\ \beta_n(\emptyset) = 
\emptyset,$$ $${\rm and}\ \beta_n(X_k) = \beta_n(Z_k) (0 \leq k < n).$$\pars
Clearly $\beta_0$ is a bijection from $X_0$ onto $Z_0.$ Since $X_0$ and $Z_0$ 
are 0-atomic$,$ $M_{X_0} = M_{Z_0}=\emptyset.$ Therefore (A) and$,$ obviously$,$ (B) hold 
for $n = 0.$\pars
Assume that (A) and (B) hold for $n,$ where $X_0$ and $Z_0$ are $n$-atomic. We 
show that (A) and (B) hold  for $n+1,$ where now $X_0$ and $Z_0$ are 
$(n+1)$-atomic. \pars
Notice that
$${\rm for}\ {\rm any}\ {\rm set}\ x,\ x \in X_0\ {\rm implies}\ x \not\subset 
X_n,\leqno {\rm [\dagger]}$$
for if $x \in X_0,$ then $x \not= \emptyset$ and $x \cap X_n = \emptyset$ by 
the $n$-atomicity of  $X_0.$ Hence$,$ it cannot be that $x \subset X_0.$ 
Similarly$,$  
$${\rm for}\ {\rm any}\ {\rm set}\ z,\ z \in Z_0\ {\rm implies}\ z \not\subset 
Z_n.\leqno {\rm [\dagger\dagger]}$$ \pars
Clearly$,$ $\beta_{n+1}$ is a map from $X_{n+1}$ into $Z_{n+1}.$ Suppose that 
$x,\ y \in X_{n+1}$ and $\zeta(x) = \zeta(y).$ Then $\zeta(x) = \zeta (y) 
\notin Z_0$ or $\zeta(x) = \zeta (y) 
\in Z_0.$ For the first case$,$ $x,\ y \notin X_0.$ Hence$,$ by Theorem 2.1.2$,$ 
$x,\ y \subset X_n$ and $\zeta [x] = \zeta(x) = \zeta (y) = \zeta[y].$ 
Since $\beta_n$ is an injection$,$ it follows that $x = y.$ In the second case$,$ 
it follows from $[\dagger\dagger]$ that $\zeta(x) \not\subset Z_n.$ But as was 
shown in the course of the first case$,$ $x \notin X_0$ implies that 
$\zeta(x) = \zeta[x] \subset Z_n.$ Hence $x \in X_0.$ The same argument 
shows that $y \in X_0.$ Again the injectivity of $\beta_n$ implies that
$ x = y.$ Consequently$,$ $\beta_{n+1}$ is an injection from $X_{n+1}$ into 
$Z_{n+1}.$\pars
To show that $\beta_{n+1}$ is a surjection$,$ let $z \in Z_{n+1}.$ If $z \in 
Z_n$$,$ then the surjectivity of $\beta_n$ yields an $x \in X_n\subset X_{n+1}$ 
such that $z = \zeta(x)\in Z_{n+1}.$ If $z \notin Z_n,$ then $z \subset Z_n.$ Hence  
again by the surjectivity of $\beta_n,$ we have that $z = \zeta[x],$ where
$x = \beta^{-1}_n[z] \subset X_n.$ 
By $[\dagger]$$,$ $x \notin X_0.$ Hence $z = \zeta[x] = \zeta(x)\in 
Z_{n+1}.$\pars
If $x,\ y \in X_{n+1},$ and $x \in y,$ 
then $y \notin X_0$ by the 
$(n+1)$-atomicity of $X_0.$ Hence$,$ $\zeta(x) \in \zeta[y] = \zeta(y).$ 
Conversely$,$ since $\zeta| X_{n+1}=\beta_{n+1}$ is a bijection onto 
$Z_{n+1},$   
it suffices to assume that $x,\ y \in X_{n+1};\ \zeta(x),\ \zeta(y) \in 
Z_{n+1}$ and $\zeta(x) \in \zeta(y).$
Then from the $(n+1)$-atomicity of $Z_0$$,$  $\zeta(y)\notin Z_0$ implies that
$y \notin X_0.$ Hence$,$ $\zeta(x)\in \zeta[y],$ and$,$ thus$,$ $\zeta(x) = 
\zeta(x^\prime)$ for some $x^\prime \in y.$ By Theorem 2.1.2$,$ $x^\prime \in 
X_n.$ Since $\beta_{n+1}$ is an injection$,$ $x = x^\prime.$ Thus $x \in y.$ 
It follows immediately from the definition of $\zeta$ that $\beta_{n+1}(x) = 
\beta_{n+1}(y)$ if and only if $ x = y.$ Consequently$,$ (A) is established for 
$n+1.$\pars
In general$,$ since $A \subset X_0 \subset X_n,$ we have that $A \notin X_0$ by 
$[\dagger]$. Therefore$,$ $\beta_{n+1}(A) = \zeta[A] = A.$ The remainder of (B) 
is easily verified for $n+1$ and by induction the proof is complete. \qed 

A criterion as to when a set $X_0$ is 
$n$-atomic for all $n \in \omega$ is very useful. Obviously$,$ if $X_0$ is a set 
of atoms$,$ then $X_0$ is $n$-atomic for all $n \in \omega.$ For the definition of $TC$, see page 54.\parm
{\bf Theorem 2.1.4} {\sl Suppose that $\emptyset \notin TC(X_0).$ If there exists 
$n \in \omega$ such that $X_0$ is not $n$-atomic$,$ then there exists some $y 
\in X_0$ such that $TC(y) \cap X_0 \not= \emptyset.$}\pars
Proof. Observe that a straightforward inductive argument shows that for each $n 
\in \omega,$ if $x \in X_n$ and  $\emptyset \notin TC(\{x\}),$ then  
$TC(\{x\})\cap X_0  \not= \emptyset.$ Assume the hypotheses of the theorem. 
Since $X_0$ is not $n$-atomic$,$ there exists $x,\ y$ such that 
$x \in y \in X_0$ and $x \in X_n.$ Since $TC(\{x\}) \subset TC(X_0)$ and 
$\emptyset \notin TC(X_0),$ we have $\emptyset \notin TC(\{x\}).$ Hence
$TC(\{x\}) \cap X_0 \not= \emptyset.$ Since $TC(\{x\}) \subset TC(y),$ it 
follows that $TC(y) \cap X_0 \not= \emptyset.$ \qed 
Application of Theorems 2.1.3 and 2.1.4 can eliminate a great deal 
of tedious work. Intuitively$,$ words in a language behave$,$ in many respects$,$ as 
it they are themselves atoms. We discuss sets of them$,$ subsets of sets of 
them$,$ etc. Since the symbol strings carry a positioning$,$ unless we extend the 
intuitive set-theoretic structure to a much more complex one$,$ it would be 
difficult to discuss the internal construction of a word in the most 
simplistic of set-theoretic languages. 
After all$,$ as a set of elements
$\{{\rm BOOTS}\}= \{{\rm BOTS}\}$ have considerably different meanings. 
This is why the 
actual intuitive ordering is indicated by the partial sequences. On the other 
hand$,$ if words seem to behave like atoms within our basic logic$,$ then 
certain statements about the number of steps in a formal deduction or the 
``number'' of words used for some purpose needs to be represented by relations 
with respect to the natural numbers. \pars

 Let $\nat$ be a set of individuals in our model for {\bf ZFH} that is isomorphic to $\omega$. The set ${\cal W}$ is assumed to have symbols that represent aspects of the theory of natural numbers (or rational, real, etc.) In the usual manner, these are assumed to be different than those symbols from $\nat$ (or other formal sets) used to analyze the set $\cal W$. Since the specific type of entity being employed is always obvious, a symbolic distinction will not, generally, be made. Relative to the symbols in countably infinite ${\cal W},$ ${\cal W} \cap \nat = \emptyset$, $\cal W$ is a set of atoms and  
$\nat$ is a disjoint countably infinite set of atoms.  
The set $\nat$ is the 
natural numbers within the ``intuitive'' and the ``formal''  portion of this model. [See note [1] at end of this section.] Let $X_0 = {\cal W} \cup \nat$.  
It is a simply matter to show that separating the original set of atoms in 
this fashion is consistent relative to {\bf ZF}. 
 Since ${\cal W} \cup\nat$ 
are atoms$,$ $X_0$ is $n$-atomic for all $n \in 
\omega.$  \pars

We now show that the set 
${\cal E}\cup \nat = Z_0$ satisfies the hypotheses of the contrapositive of Theorem 
2.1.4. 
First consider $\cal E.$ Note that 
each member of
$\cal E$ is a nonempty set and is a finite set of partial functions. That is a 
finite set of nonempty sets of ordered pairs.  
Consider any 
$y \in \cal E.$ Let $x_0 \in y.$ Then 
$x_0$ is a nonempty finite set of ordered pairs. Let $x_1 \in x_0.$ Then $x_1$ 
is a 
nonempty finite set containing one singleton and one doubleton set. Then 
if $x_2 \in x_1,$ then $x_2$ 
is a nonempty finite set of atoms. Hence$,$ if $x_3 \in x_2,$ then $x_3 \in
i[{\cal W}].$ Now none of these sets is the empty set and for each $y \in {\cal E}, \
TC(y) \cap {\cal E} = \emptyset$. For each $y \in \nat,\ 
TC(y)= \emptyset.$ Since ${\cal W} \cup \nat$ are atoms, then $Z_0$ is $n$-atomic 
for each $n\in \omega.$ Thus$,$ for our superstructure construction$,$ Theorem 2.1.3 
now applies for each $n\in \omega.$\pars 

The above finite argument is considered an effective procedure as are 
inductive definitions.  What can be claimed to be the effective procedure? 
Even though 
some might accept the effective procedure as the inductive 
definition of members in $\cal E,$ in reality$,$ it is the concept of {finite 
recognizability}
and the fact that members of $\cal E$ can be constructed from a concrete 
physical symbol model. Finite recognizability is the same concept that allows 
for the acceptance that {G\"odel numbering} generates an effective injection 
into $\nat.$ If  we assign  $g($``(''$) = 3$$,$ $g($``,''$) = 5,$ and  
$g($``)''$) = 7$$,$ then unless it is accepted (i.e. recognized) that the 
string  ``(,)'' is different from the  string ``(),'' the relation determined 
by assigning to the strings $2^35^57^7$ and $2^33^75^5$ would not be a map.
Using a concrete symbol model$,$ then from the construction 
of $\cal E,$ no object that 
is either an atom$,$ a nonempty  set composed of one or two atoms$,$ an ordered 
pair composed from these previous sets$,$ or a nonempty set 
of such ordered pairs$,$ is equal to any 
nonempty  finite set of sets of such ordered pairs.
 Thus$,$ $Z_0 = {\cal E}\cup \nat$ 
is $n$-atomic for every $n \in \omega.$ Due to (1.2.4), there is a bijection   
$\theta \colon i[{\cal W}] \to {\cal E}$ that associates each member of $i[{\cal W}]$ with a unique member of $\cal E$. This is coupled with the identity map on ${\cal W}.$ This composition yields that bijection needed for Theorem 2.1.3. Consequently$,$ by Theorem 2.1.3$,$ 
for each $n \in \omega$ the 
structures $\langle X_n,\in,=, \nat, \emptyset \rangle$ and 
$\langle Z_n,\in,=, \nat,\emptyset \rangle$ are isomorphic. \pars

For each $n\in \omega$, let  $({\cal E}\cup \nat)_n$ be the 
$n$'th level in a superstructure based upon ground set ${\cal E}\cup \nat.$ Note that relative to a superstructure based upon ${\cal W} \cup\nat$ and $({\cal E} \cup \nat) = X_n,$ for each $n \in \omega,$ there is a $m \in \omega$ such that $({\cal E}\cup \nat)_n \subset X_m \subset X_{m+1}$. Thus, we also have that $({\cal E} \cup \nat)_n \in X_{m+1}.$

The intuitive properties 
for the deductive processes with which we are concerned can be described 
within a first-order language and all hold within some 
particular $({\cal W} \cup \nat)_n.$ Hence$,$ the same properties hold in the 
corresponding 
$({\cal E}\cup \nat)_n$ through application of the isomorphism which exists 
between these 
two structures. It is$,$ in reality$,$ by means of $i$ and $\theta,$  that the basic 
logical properties 
within our intuitive theory become properties within the formal 
mathematical theory based upon 
$\nat.$   (The term ``informal'' means a restriction to superstructure entities determined by ${\cal W}$. The term ``formal'' means the entire superstructure.) 
Assuming finite recognizability, the injection $i$ is created and 
used to pass informal information into formal 
information about members of $[f_{i_w}]$ since everything is finitary in 
character. The set $[f_{i_w}]$ is finite, each $g \in [f_{i_w}]$ is finite.
Each $x \in g$ is finite, etc. This 
intuitive finitary process that is employed when  formal statements are made
within the formal portion of our {\bf ZFH} about the structure of the ordering 
of the words. \pars

In all that follows$,$ rather than continually mentioning the existence of 
isomorphisms and applying them to obtain a corresponding property in some 
${\cal E}\cup \nat)_n$ a special approach is followed. When viewed as a models, every object in the 
superstructure based upon $X_0$ has a constant name. These objects are 
uniquely determined by their set-theoretic construction. Among these constants 
are the constants $\nat$ and $\emptyset$ that are used to represent the 
natural numbers and the empty set in this section. Since on any specific 
$X_n$ we have an isomorphism $\beta_n$ from $\langle X_n,\in,=,\nat,\emptyset 
\rangle$ onto $\langle Z_n\in,=,\nat,\emptyset 
\rangle,$ if $x \in X_n,$ then $\beta_n(x)$ is the corresponding element in 
$Z_n$ characterized by the same set-theoretic property. In the same way$,$ every 
member of $Z_n$ has a constant name within our language. \pars

The following 
convention is used.  The injection $i \colon {\cal W} \to \nat$ is extended, in the usual manner, to subsets of $\cal W.$ 
Certain constant symbols used to name objects with specific 
properties in the intuitive part of the superstructure, except for $\nat$, its elements and 
$\emptyset$ are mapped by extended $i$  into a formal superstructure such as $\cal X$, where the ground set is $X_0 = {\cal W}\cup \nat$. The map $\theta$ is also extended in the same manner as $i$. Where applicable, the composition of $i$ followed by $\theta$ is denoted by {\bf bold} type face.  Also, except for members of such sets as $\nat$ and variables, most of the informal notation for functions and the like are also represented in the standard model by {\bf bold} font. For this example, let $\r L = {\cal W}$  and the consequence operator $\r C \colon \power {\r L} \to  \power {\r L}$. Then $\b C \colon \power{\b L}\to \power{\b L}$ is also a consequence operator. This notational convention is followed throughout  the remainder of this book.\pars

From these results$,$ if $\r A \subset i[{\cal W}],$ then any intuitive deductive process 
$\r k \subset F(\r A) \times \r A$ or any consequence operator $\r C\colon \power \r A \to \power 
\r A$ becomes under the isomorphism  a deductive process ${\bf k} \subset F({\bf A)} \times {\bf A}$ 
or a consequence operator ${\bf C} \colon \power {\bf A} \to \power {\bf A}.$
Notice that we do not need to consider the isomorphism on the operators $F$ or 
$\cal P$ since $\r B \in F(\r A)$ if and only if a sentence$,$ with appropriate 
constants$,$  of the following type 
holds. $\r B = \emptyset\lor \forall x(x \in \ rB \iff
(x = \r a_1\land \r a_1 \in \r A) \lor (x = \r a_2\land \r a_2 \in \r A) \lor \cdots\lor (x = \r a_n
\land \r a_n \in \r A)).$ Hence ${\bf B} \in {\bf F(A)}$ if and only if
${\bf B} = \emptyset\lor \forall x(x \in {\bf B} \iff
(x = {\b a_1\land \b a_1 \in A}) \lor (x = {\b a_2\land \b a_2 \in A}) \lor \cdots\lor 
(x = {\b a_n \land \b a_n \in A})),$  where the isomorphism 
does map $\emptyset$ onto $\emptyset$ at level $n = 1.$ In like manner$,$ the 
power set operator. (In most cases since it reveals an order, only $\cal E$ is employed.) Let $\r A \subset {\cal W}$ and let $\r K_{\r A}$ denote the set of all
deductive processes defined for $\r A.$ Now let $\r C_A$ denote the set of all 
consequence operators defined on $\power {\r A}.$ The set $\r R_A = \r K_A \cup \r C_A$ is a set 
of all {{\it intuitive human reasoning processes}} while ${\bf R_A} ={\bf K_A \cup 
C_A}$ is a set of {{\it formal human reasoning processes.}} \pars

\vfil\eject
\leftline{\bf 2.2 A Remark About 2.1}
\medskip
The basic intuitive procedure in establishing a formal model is not relative 
to structures with a universe $({\cal E}\cup \nat).$ What most be done is to 
express in a structure such as $\langle ({\cal W}\cup \nat)_n,\in,=,\nat,\emptyset 
\rangle$ informal statements about our language ${\cal W},$ where $\cal W$ is termed as informal ground set of atoms disjoint from $\nat.$ 
For named objects within such a superstructure$,$ the same {\bf bold} face convention is used for the corresponding objects within any particular $({\cal E} \cup \nat)_n$  that involves only the members of ${\cal E}.$ \pars
One additional remark is in order. In 1978 when the following  concepts within 
the discipline termed nonstandard analysis were developed$,$ they were in the 
mainstream of complexity. Today$,$ many who work in this area would consider 
them to be very simplistic in nature. To the neophyte$,$ however$,$ they may seem 
to be somewhat difficult.\par
\bigskip
\leftline{ \bf 2.3 The Nonstandard Structure}
\medskip
Now that the general and basic concepts for the deductive processes and 
consequence operators have been developed$,$ its necessary to consider  
${\cal W} \cup \nat$  as embedded into an additional structure. The same concept
that every member of the following type of superstructure  
corresponds to a constant within our language is to be used. With respect to 
the previous convention$,$ many of these constants will be denoted in bold.\pars

Recall for a moment how $\nat$ is obtained. Let set $A$ be our countably 
infinite set of atoms$,$ disjoint from ${\cal W},$ and $f\colon A \to \omega$ the 
bijection which exists from $A$ onto the set $\omega$ of natural numbers in 
our model for {\bf ZFA}. Consider $f^{-1} \colon \omega \to A$ and use $f^{-
1}$ to pass the order relation (and other necessary operations) 
on $\omega$ from $\omega$ to $A.$ For example$,$ 
this yields for each $x,\ y \in A,\ x\leq y$ iff $ f(x)\subset f(y)$ and $A$ 
inherits all the order properties for $\omega.$ Notice that since $f$ is a 
bijection that $f$ preserves equality. Denote this set $A$ by $\nat.$ \pars

We obtain a nonstandard model for a slightly different superstructure with ground set ${\cal W} \cup \nat$ than considered in section 2.1.
This is one of the two basic constructions that appear in the literature. The superstructure levels are slightly different [10. p. 
40]$,$ [17$,$ p. 110]$,$ [19$,$ p. 23]. Let $X_0 = {\cal W} \cup \nat$ and by induction$,$ let $X_{n+1} 
= \power {\bigcup \{X_k \mid k = 0,\ldots,n \}}.$ Finally$,$ let ${\cal N} = 
\bigcup \{X_n \mid n \in \omega\}.$ Consider a {$\kappa$-adequate} 
{ultrafilter}
$\cal U,$ where $\kappa > \vert {\cal N} \vert.$ By Theorem 7.5.2 in [19] or 
Theorem 1.5.1 in [9] such an ultrafilter exists in our {\bf ZFH} and is 
determined by the indexing set $J = F(\power \kappa).$ \pars
Consider the structure ${\cal M} =\langle {\cal N},\in, = \rangle.$ [Note: Since every member of ${\cal N}$ is named by a constant$,$ including the customary ones for specific objects$,$ these constants are suppressed in the notation.] By Theorem 7.5.3 in 
[19] or Theorem 1.5.2 in [19]  the {ultrapower} construction yields by definition 
3.8.1 in [9] a structure ${\cal M}_1 = \langle {\cal N}^J,\in_{\cal U},=_{\cal 
U} \rangle$ which is a nonstandard model for all sentences$,$ $K_0,$ in a 
first-order language $L$ with equality and predicates  for $\in$ and $=$ which 
hold in 
$\langle {\cal N}, \in,= \rangle.$ Assume that the cardinality of the set of 
constants of $L \geq \vert {\cal N} \vert.$ Moreover$,$ by means of sequences 
from $J$ onto $\cal N,$ the structure $\langle {\cal N}, \in, = \rangle$ may 
be considered as ismorophically embedded into ${\cal M}_1$ so that ${\cal M}_1$
is also an elementary extension of the embedded  $\langle {\cal N}, \in,= \rangle.$          
The structure ${\cal M}_1$ is also an {{\it enlargement}} of $\langle {\cal N}, \in,= \rangle.$          
A proof of The Fundamental Result may be found on page 39 of [19] (Theorem 
3.8.3) among other places. Now in [10]$,$ Theorem 3.8 establishes this for 
{BOUNDED} sentences which hold in $\langle {\cal N}, \in,= \rangle.$ Notice 
that the interpretation map from the language onto $\langle {\cal N}, \in,= 
\rangle$ has been suppressed and each member of $\cal N$ is simply to be 
considered as named by the constants in $L.$ \pars
The next step is to realize either by analysis of the ultrapower 
construction directly or by interpreting the appropriate sentences [17$,$ p. 
119]$,$ that $=_{\cal U}$ is an equivalence relation with the
substitution 
property for $\in_{\cal U}.$ Thus passing to the equivalence class
$[x]$ for 
each $x \in {\cal N}^J,$ define $[x]\in^\prime [y]$ iff $x \in_{\cal U} 
y$ for each $x,\ y \in {\cal N}^J.$ w.\pars
Now let $(X_n)^\prime$ be the objects in $\{[x]\mid x \in {\cal N}^J \}$
that correspond to $X_n$ in ${\cal N}$ under the interpretation map $I$ 
followed by the quotient map for the equivalence relation as determined by 
$=_{\cal U}$ (i.e. the   ``prime'' mapping.) It follows that $(X_0)^\prime$
behaves like atoms (urelements) and each 
$(X_n)^\prime,\  n >0$ is well-founded with respect to $\in^\prime.$ This 
comes from interpreting the appropriate bounded sentences
such as the results of Lemma 2.1 (iv) [10$,$ p. 40] where $R = X_0$ or property 
(iii) on page 23 of [19] in order to obtain the $\in^\prime$ well-founded for 
each $(X_n)^\prime,\ n \geq 1.$ For example$,$ for each $n \geq 1,$ the 
following sentence
$$\forall x((x \in \nat) \to \neg \exists y((y\in X_n) \land (y \in x))) \eqno 
(2.3.1)$$
holds in the structure $\langle {\cal N},\in, = \rangle.$ Lastly$,$ each $(X_n)^\prime,\ n \geq 1,$ is 
well-founded with respect to $\in^\prime$ since ``If $x \in y\in X_n$$,$ 
then $x \in X_0 
\cup X_{n-1}$'' $(n \geq 1)$ holds in $\langle {\cal N}, \in,=\rangle.$ 
Consequently$,$ the Mostowski Collapsing Lemma [1$,$ p. 247] or [17$,$ p. 120] 
can be inductively applied to each $(X_n)^\prime$ and obtain a corresponding 
set $\hyper X_n.$ Specify the set $\hyper {X_0}$ to correspond to $(X_0)^\prime$ 
and we have a unique collapse. As a result of this$,$ the structure $\langle 
\bigcup \{\hyper X_n \mid 
n \in \omega \},\in,= \rangle = \hyper {{\cal M}} = \langle \hyper {{\cal 
N}},\in, = \rangle$  is a set-theoretic model for all bounded sentences that 
hold in $\langle {\cal N}, \in, = \rangle.$ Recall that a bounded sentence in 
a first-order language $L$ is one for which each quantified variable is 
restricted to an element of ${\cal N}.$ The composition  of the interpretation 
map $I^{\cal U},$ the quotient map $\prime$ and the collapse yield the * map 
from the structure $\langle {\cal N},\in,=\rangle$ into $\hyper {{\cal M}}$ 
and maps any element $ a \in {\cal N}$ to the element $\hyper a$ preserving 
all of the usual properties for a normal$,$ enlarging and comprehensive 
monomorphism. For each $B \in {\cal 
N},$ let $^\sigma B = \{ \hyper x\mid x \in B\}.$ (This definition does not correspond to that used by some other authors.)
The * notation is also 
not placed on elements of $X_0$ when they are considered as mapped into 
$\hyper {X_0}$ by the * map. Observe that for each $B \in {\cal N},\ 
^\sigma B \subset \hyper B.$ Technically, where used, $B \subset \hyper B$ also means 
$^\sigma B \subset \hyper B$).\pars
Now to complete the construction, begin with the set 
$Y_0 = \hyper {X_0}$ and construct a superstructure with $Y_0$ as the ground set 
as {\bf defined in this first example}. Let $Y_{n+1} = \power {\bigcup \{Y_n\mid
n \in \omega \}}$; and let $Y = \bigcup \{ Y_n \mid n \in \omega\}.$
\pars
   For the above$,$ some general principles such as the Mostowski Collapsing 
Lemma have been used in order to obtain $\hyper {\cal M};$ however$,$ an 
explicit construction appears on pages 44 and 45 of [19]. In actuality for the next constructed superstructure$,$ the one used in this book$,$ we 
intend to
use only a small portion of $\hyper {{\cal M}}.$ Indeed$,$ we apparently need to 
use a small hierarchy of the $\hyper {X_n}$  objects. You could$,$ if you 
wished$,$ restrict the G-structure to say only the $n < 100$ levels. However$,$ 
this will not be done for fear of not selecting a correct upper bound for 
$n.$\pars
 For the results in this book$,$ I advocate for our superstructure a construction as defined in section 2.1$,$ where $X_0 = {\cal W} \cup \nat,$ and the nonstandard model as constructed on pages 83 - 88 and Theorem 6.3 in Hurd$,$ A. E. and P. A. Loeb$,$ (1985)$,$ ``An Introduction to Nonstandard Real Analysis,'' Academic Press$,$ Orlando. [Note: This construction also appears on pages 42-49 in Loeb and Wolff$,$ (eds) (2000)$,$    ``Nonstandard Analysis for the Working Mathematician,'' Kluwer Academic Publishers$,$ London. Also $X_n({\rm 2.1}) = X_0 \cup X_n({\rm 2.3}),\ n\geq 1$.\ddag] This construction simply needs to be restricted to our language with $\in$ and $=$$,$ where $=$ is interpreted as set-theoetic equality on sets and the identity on atoms. For the first superstructure$,$ constructed using the procedure in section 2.1$,$ let ${\cal N} = {\cal X}.$ The second superstructure constructed using this procedure has as its ground set $Y_0 = \hyper {X_0}$ and$,$ as before$,$ $Y = \bigcup\{Y_n\mid n \in \omega\}.$ This leads to the G-structure ${\cal Y} = \langle Y, \in, = \rangle$ where since I apply this to logical operations this structure is call the {{\it Grundlegend Structure}}.\pars
\pars
Now to summarize. The consistency of {\bf ZF} implies the consistency of {\bf 
ZFH} and one can apparently use a model of {\bf ZFH} to obtain the nonstandard 
structure $\hyper {\cal M} = \langle \hyper {\cal N}, \in, = \rangle.$ The 
set $\hyper {\cal N}$ is dependent upon the of atoms $\nat,$ the atoms of 
{\bf ZFH} with the order induced by $\omega.$ Any sentence in a appropriate 
first-order language in which each quantified variable is restricted to an 
element of ${\cal N}$  (i.e. bounded variable) will$,$ when each constant is 
replaced by the * of the constant$,$ give a true statement about the structure 
$\hyper {\cal M}.$  Moreover$,$ $\hyper {\cal M},$ at the least$,$ has bounds for 
all standardly definable concurrent relations. For notation$,$ we denote 
for each $ n \in X_0,\  \Hyper n = n.$  In addition$,$ all properties of the * map as 
listed in [10]$,$ [17]$,$ [19]$,$ among other places$,$ hold true. Next some unusual 
names for G-structure objects will be adopted in order to reflect our 
application to languages and logics.\pars
Recall some of the basic terminology associated with $\cal Y.$ For each 
$A \in {\cal N},\ A$ is called a {{\it  standard entity.}} The set
$\hyper A$ is often called an {{\it (internal) standard entity}} or better 
still an {{\it extended standard entity}} in $\cal Y.$ If $b  \in \hyper A,\ 
A \in {\cal N},$ then $b$ is called an {{\it internal entity}}. Indeed$,$ $b 
\in {\cal Y}$ is internal iff there is some $X_n$ such that $b \in \hyper 
X_n.$ Any entity of $\cal Y$ which is not internal is called 
{{\it external}}. These terms are generally used throughout nonstandard 
analysis$,$ but for our present purposes they are modified as follows: 
Any entity of $\cal Y$ is a {{\it subtle object,}} some appropriate members 
of 
$\cal N$ are {{\it human objects}} and any entity in $\cal Y$ which is not 
the $\sigma$ of a member of $\cal N$ or the * of a member of $\cal N$ is a {{\it purely 
subtle object.}} Please refer to the basic references [11]$,$ [16]$,$ [19] for 
other terminology and the properties of *. So as to avoid symbolic confusion,
from this moment on$,$ the entire or the major part of any symbol used to represent objects within a language 
and within our intuitive model will be denoted by Roman type.
 \par \vfil\eject
{\leftline{\bf 2.4 General Interpretations}}
\medskip
Throughout  this work on ultralogics$,$  ${\rm L}_0$ will denote the usual set of 
propositional formulas (wff) constructed from the connectives $\neg,\ 
\lor,\ \land,\ \to,\ \iff,$ say as done by Kleene [7$,$ p. 108] and ${\rm L}_1$ is a set of predicate formulas with equality considered as an 
extension of ${\rm L}_0$ as say constructed on page 143 of [7]. We also use the 
usual assortment of set-theoretic
abbreviations when we consider the special predicates $\in$ and $=.$ Of 
course$,$ ${\rm L}_1$ is called a {{\it first-order language.}} Assume that 
${\rm L}_1 \subset {\cal W}$ and that the set of all predicate symbols is a subset of 
$\{{\rm P}_i\mid i \in \nat\}.$ It is important to realize that any intuitive 
set-theoretic deduction process$,$ and the like$,$ that is discussed relative to 
${\cal W}$ is to be embedded by the map $\theta$ to a corresponding process 
relative to $\cal E.$ This also applies to a member of $\cal W$ and the $i$ 
injection. The results of any *-transfer of statements which hold relative 
to $\cal E$ or $i[{\cal W}]$ are modeled 
in  $\hyper {{\cal M}}.$ Also$,$ any results relative to $\cal E$ or $i[{\cal 
W}]$  (i.e. with respect to standard objects) can be referred back to 
corresponding intuitive objects relative to ${\cal W}$ by means of either the 
maps $i^{-1}$ or $\theta^{-1}.$ Moreover$,$ in order to simplify notation 
somewhat any formal first-order statement that explicitly involves $i(w)$ 
individuals will be written with the $i$ deleted from the notation if no confusing 
results from such an omission. \pars 
For example$,$ the sentence\pars
\line{\hfil $ \forall x(x\in \nat \land x \geq 1 \to \exists y (y \in P 
\land \forall w(w \in \nat\land$\hfil} \smallskip 
\line{\hfil $0 <w\leq 
x\to y(w) = {\rm very,}\sp \land y(0) = {\rm just})))$\hfil (2.4.1)}\parm 
\noindent 
is a slight simplification of the following sentence\parm \line{\hfil $ 
\forall x(x\in \nat \land x \geq 1 \to \exists y (y \in P \land \forall 
w(w \in \nat\land$\hfil} \smallskip 
\line{\hfil $0 <w\leq x\to y(w) = 
i({\rm very,}\sp )\land y(0) = i({\rm just}))))$\hfil (2.4.2)}\parm 
Most of the 
following investigation is concerned with specific elements of $P$ and 
specific human reasoning processes with respect to $\cal E.$ Using the 
previous example$,$ the *-transfer yields\par 
\line{\hfil $ \forall x(x\in 
\hypernat \land x \geq 1 \to \exists y (y \in \hyper P \land \forall w(w 
\in \hypernat\land$\hfil} \smallskip 
\line{\hfil $0 <w\leq x\to y(w) = 
\Hyper i({\rm very,}\sp )\land y(0) = \Hyper i({\rm just}))))$\hfil (2.4.3)}\parm 
\noindent and which 
holds in $\hyper {{\cal M}}.$ Notice that we do not place * on the order 
relation $<$ since we assume that it is but an extension of the simple order 
$<$ on $\nat$ satisfying all of the same first-order properties. Also note that $\Hyper i({\rm very},\sp) = i({\rm very},\sp)$ etc. Let $\nu \in 
\hypernat - \nat = \nat_\infty$ (i.e. the infinite numbers). Then there exists 
in $\hyper P$ a *-partial sequence$,$ say $f,$ such that for each $w \in 
\hypernat,\ 0 < w\leq \nu,\ f(w) = {\rm very,}$ and $f(0) = {\rm just}.$ 
Hence even though members of $[f]$ are not  readable sentence in our sense$,$ we 
can read the elements in the range of $f$ as well as reading the intuitive 
ordering when $f$ is restricted to $[0,n],\ n \in \nat.$ This gives an 
intuitive interpretation for such an $f$ when it is so restricted to such 
standard segments as well as knowledge of the properties of the ordering when 
not so restricted. Observe also$,$ that if $f \in \hyper P - P,$ then there 
exists some $\nu \in \hypernat$ such that $f\colon [0,\nu ]\to \Hyper 
{i[\Hyper {\cal W}]},$ where $[0,\nu ] = \{x\mid (x \in \hypernat ) \land (0 \leq x \leq \nu 
)\}.$ Often, in our formal statements, parentheses are suppressed and the strength of connectives notion is used.\parm 
\noindent {\bf 2.5 Sets of Behavior Patterns} \par
\medskip
In certain applications of subtle consequence operators the following 
construction is useful. This is all relative to what is called {{\it 
adjective reasoning}} and any equivalent form. Let ${\rm B}^\prime$ denote a list of 
names or simple phrases that are used to identify specific behavior patterns. 
These terms are taken from a specific discipline language and are$,$ as usual$,$ to 
be considered as elements of $\cal W.$ For example$,$ the set ${\rm B}^\prime$ could 
be taken from the discipline called  {{\it psychology}} and each term 
could identify a 
specific human behavior pattern$,$ as general as such concepts as ``kind'' or 
``generous.'' You can also include any synonyms that might be equivalent 
to the members of ${\rm B}^\prime.$ Now consider B constructed as follows: an
element ${\rm b}\in {\rm B}$ if and only if b is a {{\it qualifiable}} form of a 
member of ${\rm B}^\prime.$ That is each ${\rm b} \in {\rm B}$ is a 
${\rm b}^\prime \in 
{\rm B}^\prime$ where ${\rm b}^\prime$ is written in a form so that it can be modified 
by the word {{\it very}}. (Or, such words as ``great,'' ``greater.'') Let ${\rm B}= {\rm C}_0,\ {\rm C}_1 = \{{\rm very,}\sp 
{\rm c}\mid {\rm c} \in 
{\rm C}_0\}.$ By induction let ${\rm C}_{n+1} =  \{{\rm very,}\sp {\rm c}\mid 
{\rm c} \in 
{\rm C}_n\}.$ Then an intuitive set of {{\it modified behavior patterns}} is the 
set ${\rm BP} = \bigcup \{{\rm C}_n \mid n \in \nat \}.$ The {{\it formal 
modified behavior patterns}} is the set ${\bf BP} = \bigcup \{
{\bf C_n}\mid n \in \nat\}.$ Notice that each ${\bf C_m}$ is a finite 
set.\pars
In certain cases$,$ the intuitive set BP is associated with a set of formal
propositional statements in $\cal W.$ Let ${\rm L}_0$ be our propositional language 
constructed from a denumerable set of atoms $\{{\rm P}_i\mid i \in \nat \}.$ 
Since B is finite$,$ then there exists an injection $j \colon {\rm B} \to 
\{{\rm P}_i\mid i \in \nat \}$ and $\{{\rm P}_i\mid i \in \nat \} - j[B]$ is 
denumerable. Let ${\rm V}\in \{{\rm P}_i\mid i \in \nat \} - j[{\rm B}].$ 
\noindent Let the 
symbol string ``very,$\sp$'' correspond to the partial formula 
``V $\land$''. Then proceed to construct ${\rm BP}_0$  
as follows:
${\rm E}_0 = j[{\rm B}],\ {\rm E}_{n+1} = \{({\rm V} \land {\rm x})\mid {\rm x} 
\in {\rm E}_n \}.$ Then$,$ finally$,$ ${\rm BP}_0 = \bigcup \{{\rm E}_n \mid
n \in \nat\}.$\pars
In what follows$,$ the modeling of human reasoning processes is often approached 
from two different points of view. First$,$ from the viewpoint of such sets is 
BP$,$ as well as many others$,$ we have the constructed set of  {{\it meaningful 
sentences}} in the sense of {Tarski.} Thus$,$ such strings or symbols become 
our formal language and a simple observer language (i.e. metalanguage) is used 
to investigate deductive processes on BP. These are mapped to the formal deductive 
processes on {\bf BP.} However$,$ many of these deductive processes on a given
BP can be associated 
with other formal processes in ${\rm L}_0,$ especially with respect to ${\rm BP}_0.$ 
Hence$,$ whenever possible it is acknowledged that there are at least two ``models''
for various BP type statements$,$ among others$,$ that are being in investigated.
The basic model (and probably the simplest) is that based on BP. Then a 
somewhat more complex model is based on ${\rm L}_0.$ The purest is probably more 
comfortable with the formal languages ${\rm L}_0$ and ${\rm L}_1.$ I feel$,$ however$,$ that 
BP is as meaningful a set of sentences formed by constructive methods as is 
the set ${\rm L}_0$ and the various forms in BP are easily recognized. \pars
\ddag This is an important fact. Let $X_{n+1}(X) = X_n(X) \cup \power {X_n(X)},\ n \geq 0,\ X_0(X) =X$ (Def. 2.1) and $X_{n+1} = \power {X_0 \cup\cdots X_n},\ n\geq 0,\ X_0 = X$ (Def. 2.3), where $X$ is the set of individuals. For Def. 2.3, we also have that $X_p \subset X_n,\ 1 \leq p\leq  n$ and $X_{n+1} = \power {X \cup X_n}, n\geq 0.$ We show by induction that $X_n(X) = X \cup X_n,\ n\geq 0.$ First, let $n = 0$ on the right. Then $X_0(X) = X,\ X_0 = X \Rightarrow X_0(X) = X \cup X_0.$ Now for the specific inductive form, let $n = 0.$ Then $X_1(X) = X \cup \power {X_0(X)} = X \cup \power {X} = X\cup \power {X_0} = X \cup X_1.$ Assume result holds for $n$. Then $X_{n+1}(X)= X_n(X) \cup \power {X_n(X)} = X\cup X_n \cup \power {X \cup X_n} = X\cup X_n \cup X_{n+1} = X \cup X_{n+1}$ and the result follows by induction.\parm

[1] (14 DEC 2012). The set $\cal W$ (and later ${\cal W}'$) was added to the ground set on this date. This has been done to provide an additional formal structure to enhance analysis.Using the members of a language itself as constituents of a ground set for a model is well established [13, p. 70]. However, it is the set $\cal E$ that is generally more significant for our purposes than members of the language itself since they represent the significant aspects of the formation of ``words'' whether they be formed by symbols, diagrams, images or coded sensory information. Hence, in this theory, members of $\cal E$ and $\Hyper {\cal E}$ still remain the basic form for a ``word'' or ``hyper-word.'' \pars

After developing the basic aspects of this approach, it was discovered that Robinson [15, section 3] also developed a nonstandard approach to sets of symbols.  I have noted this in more recent versions. (Also see  Geiser, J. T. (1968). "Nonstandard logic," J. Symbolic Logic 33(2):236-250.)  The idea of incorporating $\cal E$ as a way to include how languages are constructed is not part of the Robinson foundations. \pars

The set $\cal W$ can contain the language for various mathematics theories such as an appropriate portion  ${\cal T}(\omega)$ of the theory  of natural numbers. Each member of ${\cal T}(\omega)$ corresponds to objects in $\cal N$. As mentioned, one can consider members of $\cal W$ as written in a different color than any other symbols used for any other purposes. In some cases, the ``prime'' notion for the symbols expressing statements about members of $\cal N$ other than members of $\cal W$ is employed. For example, the expressions $2' <' 3'$ and ${\cal T}'(\omega')$ are in extended ${\cal W}'$. External to $\cal M$, one can state that the expression  $2' <' 3'$ is a member of ${\cal T}'(\omega')$.  Or we state that $2' <' 3'$ holds. (One can actually include an additional model for this purpose.) This corresponds directly to a statement $2 < 3$  that ``holds'' in $\cal M$ where $2, <, 3$ are names for the corresponding ``formal'' objects. In general, if used, ``primed'' statements of this type are expressed directly in terms of the corresponding ``not primed'' expressions. Robinson keeps the statements used to discuss behavior of the members of his set of symbols distinct from those in $\cal W$ by simply defining such a set and leaving the rest to ones intuition. \pars
The use of the embedding $i$ now seems of little significance. The embedding was used so that $\cal E$ could simply be considered as entities from the theory of natural numbers with its long history of empirical consistency. In the beginning of nonstandard analysis where simplified type theory was employed and formal set-theory was not considered, such a consistency notion might be useful. But since formal set-theory is now being considered, any consistency considerations depends upon the assume consistency of the set-theory axioms being employed. Hence, as demonstrated, the removing of the function $i$ from both the foundations and expressions should not effect any of the interpreted results. \pars

If $i$ is so removed, then members of $\cal E$ are still equivalence classes but they are now partial sequences of members of the language $\cal W$ rather than the codes produced by application of $i$. From the viewpoint of the nonstandard model, this would mean that rather than an ultraword being considered as a partial 
hyper-sequence of members of $\hypernat$ with some ``symbols'' being represented by members of $\hypernat - \nat,$ some *symbols are represented by members of $\Hyper {\cal W} - {\cal W}.$ \pars

There is, of course, a bijection $w\colon {\cal W} \to {\cal E},$ there $w(a) = [g]$ and there is a $f \in [g]$ such that $f \in T^0$ and $f(0) = i(a)$ (or simply $f(0) = a.$ This bijection may be useful for further developments of this Theory of Ultralogics. \par
\vfil\eject
\centerline{\bf 3. DEDUCTIVE PROCESSES}\par
\bigskip
\leftline{\bf 3.1 Introduction.}\par 
\medskip
We approach the investigation of various {special deductive processes} by 
defining them intuitively for some $\m A \subset {\cal W}$ as being
$\m k \subset F(\m A) \times \m A$ or $\m C\colon \power {\m A} \to \power {\m A}.$ These sets are 
all considered mapped to objects relative to $\power {{\cal E}} \times
\power {{\cal E}}$ for formal investigation. In at least one case$,$ a map 
$C\colon  \power {A} \to \power {A}$ is defined for each 
nonempty $A \in \hyper {{\cal N}}$ and it is shown (trivially) that such a map 
is a consequence operator. Any $C \in \hyper {{\cal N}}$ which satisfies (in 
$\hyper {{\cal N}}$) axioms (2)$,$ (3)$,$ (4) or their *-transform is a
{{\it subtle consequence operator}} or {{\it subtle reasoning process}}$,$ 
where for convenience $C$ is restricted to $A \subset \hyper {{\cal E}}.$\pars
\bigskip
\leftline{\bf 3.2 The Identity Process.}
\medskip
Let $A\subset {\cal E}$ be any nonempty set. For each $B \subset A,$ define
$I(B) = B.$ Obviously$,$ this is the {{\it identity operator}} from
$\power A$ onto $\power A.$ \parm
{\bf Theorem 3.2.1.} {\sl Let $A \subset {\cal E}$ be nonempty. Then the identity operator on 
$\power A$ is a consequence operator.}\pars
Proof. Let $B \subset A.$ Then $I(B) = B$ implies that $B \subset I(B) \subset 
A.$ Moreover$,$ $I(I(B))= I(B) = B$ for each $B  \subset A.$ Finally$,$ $B = I(B) 
= \bigcup \{F\mid F \in F(B)\} = \bigcup \{I(F)\mid F \in F(B)\}$ and the 
result follows. \qed
Let $H$ be a nonempty set of {Tarski type deductive processes.}  That is 
if $h \in H,$ then $h \subset F({\cal E}_1) \times {\cal E}_1$ for some ${\cal 
E}_1 \subset {\cal E}.$ Also let $H_0$ be a nonempty set of consequence 
operators on some $\power {{\cal E}_2}$ for ${\cal E}_2 \subset {\cal E}.$
Then $\hyper H \cup \hyper H_0 = D_0$ is considered a set of {subtle 
reasoning processes}. Notice that if $C \in H_0,$ then $\Hyper C \in D_0$ may 
not be a consequence operator under our definition. The first reason for this 
is that axiom (4) *-transforms to read that for every internal subset of $B 
\subset \hyper {{\cal E}_2},\ \Hyper C(B) =\bigcup \{ \Hyper C(F) \mid F\in 
\hyper F(B)\}.$ For the sentence\parm
\line{\hfil $\forall x(x \in \power {{\cal E}_2} \to \forall w(w \in {\cal 
E}_2 \to (w \in C(x) \iff$\hfil}\pars
\line{\hfil $\exists y(y \in F(x)\land w \in C(y))))).$\hfil (3.2.1)}\parm
\noindent holds in ${\cal M};$ hence in $\hyper {\cal M}.$ As is well known a 
*-finite set need not be finite. However$,$ there is at least one map from 
$\power A$ into  $\power A$  for any $A \subset \hyper {{\cal E}_2}$ which is a 
true consequence operator as shown by Theorem 3.2.1. Consider any infinite
$A \subset {\cal E}_2.$ Then no map $C \colon \power{\hyper A} \to 
\power{\hyper A}$ can be written as an extended standard map (i.e. the star of 
a standard map) from $\power A $ into $\power A.$  This follows from the next 
result.\par
{\bf Theorem 3.2.2.} {\sl Let infinite $A \subset {\cal E}_2$ and $G \colon
\power A \to \power A.$ Then there exists a subset of $\hyper A$ upon which
$\Hyper G$ is not defined.}\pars
Proof. Let infinite $A \subset {\cal E}_2,\  G \colon
\power A \to \power A$ and $D(G)$ be the domain of $G = P_1[G].$ For an 
appropriate $X_m$$,$ the 
sentence\par
\vfil
\centerline{}
\eject
$$\forall x(x \in X_m \to (x \in D(G) \iff x \in \power A)). \leqno (3.2.2)$$
holds in ${\cal M};$ hence in $\hyper {\cal M}.$ The *-transfer reads
$$\forall x(x \in \hyper {X_m} \to (x \in D({\Hyper G}) \iff x \in \Hyper 
{({\power A})})),\leqno (3.2.3)$$
when the elementary properties of the *-map are applied. The $\Hyper G$ is 
only defined on the internal subsets of $\hyper A.$ Since $\hyper A - A$ is 
external then this result follows. \qed
{\bf Corollary 3.2.2.1} {\sl There exist purely subtle reasoning 
processes.}\parm
As to the cardinality of ${\cal E},$ it follows immediately that since each 
$x  \in \cal E$ is finite$,$ then $\vert {\cal E} \vert= \vert \nat \vert.$ Note the 
following that will be used throughout this investigation. Recall that the 
identification $^\sigma X_0 = X_0$ is being used. Then if $f \in P_H$ it 
follows$,$ since $f$ is a finite sequence of members of $X_0$ that $\hyper f = f$ 
under this identification.  Also$,$ since for each $f \in P_H$ the set 
$[f]$ is finite$,$ then $\Hyper [f]= [f].$ Thus $^\sigma {\cal E} = {\cal E}.$ 
This reduction of finite sets of 
finite sets of partial sequences continues to other cases such as 
$^\sigma (F(P_H)) = \{ \hyper A\mid A \in F(P_H)\} = 
\{A\mid A \in F(P_H)\}= F(P_H).$ \pars
With the above results in mind$,$ it follows that each $x  \in \hyper {{\cal E}} 
- {\cal E}$ is a nonfinite *-finite subset of $\hyper P$ for the sentence
$\forall x( x \in {\cal E} \to x \in F(P))$ holds in $\hyper {{\cal M}}.$ 
Consequently$,$ $\hyper {{\cal E}} \subset \Hyper {(F(P))}$ and 
${\cal E} \subset F(P)$ imply that $\hyper {{\cal E}} - {\cal E} \subset 
\Hyper {(F(P))} - F(P).$  
Let infinite $A \subset {\cal E}.$ Then 
it is an important fact that
there exists a *-finite $B \in \Hyper (F(A))$ such that $^\sigma A \subset B 
\subset \hyper A \subset \hyper {{\cal E}}.$ For let $Q = \{(x,y)\mid y \in 
F(A)\land x \in A\land x \in y\}.$ Assume that 
$(a_1,b_1),\ldots,(a_n,b_n) \in Q.$ Then letting $b = b_1 \cup \cdots 
\cup b_n $ it follows that $(a_1,b),\ldots,(a_n,b) \in Q.$ Therefore$,$ $Q$ is a 
standard concurrent relation. Thus there exists some $B \in \Hyper (F(A))$ 
such that $\hyper x \in B$ for each $x  \in A.$ Now internal $B$ is not finite 
since $^\sigma A$ is not finite. Indeed$,$ as is well know $\vert B \vert = 
\vert [0,\nu ] \vert \geq 2^\omega,$ where $\nu \in \hypernat - \nat = 
\nat_\infty.$ Thus $\vert \hyper A \vert \geq 2^\omega.$ From the above 
remarks$,$ it also follows that for $E \in \hyper {{\cal E}} - {\cal E},\ 
\vert E \vert \geq 2^\omega.$ \pars
A few other useful results are easily obtained. For 
example,\parm
(i) if ${\cal P}\colon \power A \to Y$ and $B \in  \power A,$ then 
$\Hyper {(\power B}) = \hyper {{\cal P}}(\hyper B)$.\pars              
(ii) Consider the finite power set operator $F.$ If $F \colon \power A 
\to Y$ and $B \in  \power A,$ then 
$\Hyper {(F(B))} = \hyper {F}(\hyper B).$\pars 
(iii) If $C$ is a map from $\power A$ into $\power B,$ then 
for $D \in \power A$ it follows that $\Hyper {(C(D))} = \Hyper {C}(\hyper 
D).$\parm
The proofs of (i)$,$ (ii) and (iii) are easily obtained. Indeed$,$ all three follow from the 
formal definition of a map. 
 First$,$ assuming that $A,Y \in {\cal N}.$ All 
the objects with which we shall be concerned \break will also be members of 
${\cal N}.$  Indeed$,$ if necessary to obtain bounded sentences$,$ we know that there 
is some $n \in \nat$ 
such that everything needed to characterize (i)$,$ (ii)$,$ and (iii) are members 
of $X_0 \cup X_n.$ 
For example$,$ consider (i). Then the two sentences $\forall x(x \in \power A 
\to \exists y(y \in Y\land (x,y) \in {\cal P})),$ and $\forall x \forall y 
\forall z(x \in \power A \land y \in Y \land z \in Y \land (x,y) \in 
{\cal P} \land (x,z) \in {\cal P} \to y = z)$ imply$,$ by *-transfer$,$ that 
$\hyper {{\cal P}}$ is a map from $\hyper {{\power A}}$ into $\Hyper Y.$ 
Further$,$ $(B, {\power B}) \in  {\cal P}$ implies that $(\hyper B, \Hyper 
{(\power B)}) \in \hyper {{\cal P}}.$ Hence$,$ since $\hyper {\cal P}$ is a map,
mapping notation yields that $\Hyper 
{(\power B)}) = \hyper {{\cal P}}(\hyper B).$ Now (ii) follows in like manner.
Indeed$,$ the following set of sentences shows that $\hyper F$ generates the 
hyperfinite subsets of any internal subset of $\hyper A.$ \parm 
\line{\hfil $\forall x \forall w(x \in {\cal P}({\cal P}(A)) \land
w \in \power A \to ((A,x) \in F \iff$\hfil}\pars
\line{\hfil $x = \emptyset \lor\exists y \exists z 
(z \in X_n\land y \in \nat \land \forall v(v \in A\land v \in x 
\to$\hfil}\pars 
\line{\hfil$\exists i(i \in \nat\land 0 \leq i \leq y\land z(i) = v)))))$\hfil (3.2.4)}\parm
It is well known that$,$ in general$,$ if $A = \bigcup\{B_i\mid i \in \nat \},$  then 
$\hyper A \not= \bigcup\{\hyper B_i\mid i \in \nat \}.$  However$,$ if 
$A= \bigcup \{B\},$ then by *-transfer of the definition it follows that 
$\hyper A = \bigcup\{\hyper B\} = \Hyper (\bigcup \{B\}).$ This is incorporated in 
the proof of the next result. \parm      
{\bf Theorem 3.2.3 } {\sl Let $A \in {\cal N}, \ A \subset X_n.$ If $B$ is a partition of $A,$ 
then $\Hyper B$ is a partition of $\hyper A.$}\pars
Proof. The sentences \parm
\line{\hfil $\forall x(x \in X_n \to (x \in A \iff \exists y(y \in B \land  x \in 
y)));$\hfil}\pars
\line{\hfil $\forall x\forall y(x \in B\land y \in B \to  x = y\lor  x 
\cap y = \emptyset),$\hfil (3.2.5)} \parm
\noindent hold in ${\cal M}$; hence in ${\cal M}.$  Thus by *-transfer$,$ $\Hyper B$ is a 
partition of $\hyper A$ and $\hyper A= \bigcup\{\hyper B\} = \Hyper 
(\bigcup \{B\}).$ This completes the proof. \qed
\bigskip
\leftline{{\bf 3.3 Adjective Reasoning} (Also see page 35.)}
\medskip
Following the ideas of Tarski (20) it appears that the set {BP is a meaningful 
set of sentences}. Define for BP an intuitive deductive process as follows:
Let $\rm A \in F({\rm BP}).$ Then $\rm A \vdash_a {\rm b}, \ {\rm b} \in 
{\rm BP}$ if ${\rm b} \in \rm A$ or b is obtained from some ${\rm x} \in \rm A$ be 
removing a finite number $(\not= 0)$ of ``very,$\sp$'' strings from x. Due to its form, this 
process $\vdash_a$ is termed {{\it adjective reasoning}.}  Denote the 
relation in $F({\rm BP}) \times {\rm BP}$ obtained by $\vdash_a$ by the symbol
``$a$.'' Let nonempty $\rm B \subset {\rm BP}.$ Then for each ${\rm b} \in \rm B$ it 
follows that $(\{{\rm b}\},{\rm b}) \in a.$ Hence  ``$a$'' is singular and 
$C_a$ satisfies axiom (2) of the Tarski axioms. Let $(\rm A,{\rm b}) \in a$ and 
$\rm A \subset \rm B \in F({\rm BP}).$ Then $(\rm B,{\rm b}) \in a$ since $b$ is obtained 
entirely from an element in $\rm A.$ Assume that $(\rm A,{\rm b}_1),\ldots, (\rm A,{\rm 
b}_n) \in a$ and that $(\{{\rm b}_1,\cdots,{\rm b}_n\},{\rm c}) \in a.$ Then c 
is either some ${\rm b}_i,\ i = 1,\cdots,n;$ or c is obtained from some ${\rm 
b}_i$ by removing finitely many ``very,$\sp$'' symbol strings. But either this 
${\rm b}_i \in \rm A$ or this ${\rm b}_i$ is obtained from some ${\rm x} \in \rm A$ 
also by removing finitely many ``very,$\sp$'' symbol strings. Thus c is either 
an element of $\rm A$ or is obtained from $\rm A$ by removing finitely many 
``very,$\sp$'' symbol strings from a member of $\rm A.$  Thus $(\rm A,{\rm c}) \in 
a$ and $C_a$ is a consequence operator on $\power {{\rm BP}}$ by Theorem 
1.3.3. For the next results$,$ recall that when no confusion might occur
the set $^\sigma D$ is denoted by $D.$ \pars
A remark concerning notation is necessary. Two special abbreviations are 
used in certain explicit formal sentences. The first is the symbol $[x]$ for 
$x \in P.$ This denotes the unique object $z$ that satisfies the sentence
$$\exists !z(z \in {\cal E} \land x\in z \land x \in P),\eqno (3.3.1)$$
where $\exists !zA(z)$ means $\exists z(z \in {\cal E} \land \forall y
(y \in {\cal E}\to (A(y) 
\iff z = y))).$ Now in the first formula in the proof of Theorem 3.3.2 the 
formula $\exists !z_1 \exists !z_2 ((z_1 \in {\cal E}) \land (z_2 \in {\cal 
E}) \land (y \in z_1) \land (y_1 \in z_2) \land (z_2 \in {\bf 
C_a}(\{z_1\})))$ could be inserted. Also$,$ the notation $\{x\}$ denotes the unique singleton set 
that satisfies
the following for any $A \in {\cal N},$\parm
\line{\hfil $\forall z(z \in A \to \exists !x(x \in \power A \land 
$\hfil}\pars
\line{(3.3.2)\hfil$\forall w(w \in 
A\land w \in x \to w = z))).$\hfil }\parm
\noindent We could insert for $z_2 \in {\bf C_a}(\{z_1\})$ the formula\parm
\line{\hfil$\exists !z_3(z_3 \in \power {{\cal E}} \land \forall w_3( w_3 
\in {\cal E} \land w_3 \in z_3 \to$\hfil}\pars
\line{(3.3.3)\hfil $w_3 = z_1) \land z_2 \in {\bf C_a}(z_3)).$\hfil }\parm 
Of course$,$ these formulas are not inserted$,$ but the appropriate abbreviations 
are used when needed. Recall that only constants which represent elements of 
${\cal N}$ are  ``starred'' in either the *-transform or any explicit partial 
formula obtained from the more general statement. All the {internal objects} 
which are not standardly internal take on constant names from an extended 
language. Thus if $A \in \hyper {\cal E} - {\cal E},$ then $A = [a],$ where 
$a \in \hyper P.$ The same holds for any {singleton set.} If $a \in \hyper A - 
A,$ then we write $\{a\}.$ For each $A \in {\cal N}$  and any nonempty finite 
$\{a_1,\ldots,a_n\} \subset A$ it follows that $\Hyper \{a_1,\ldots,a_n\}=
\{\hyper a_1,\ldots,\hyper a_n\} \subset \{\hyper x \mid x\in A\} = {^\sigma 
A}.$ \pars
In what follows$,$ other simplifying processes are employed when writing 
formal sentences. In many cases$,$ these 
sentences  
do not appear to be written in the special bounded form. In all cases, the additional formal 
expressions can be easily added. In general$,$ this is done 
by the addition of another $A\ \land$ type expression and an 
equivalent formula obtained$,$ or $A \to $ when 
$\iff$ appears. Many of the missing 
expressions are of these types. Here is one example of this process. Let 
$T = i[{\cal W}].$\pars
Consider the set of all ``natural number'' intervals (i.e. segments)
$H_1 = \{[0,n]\mid n \in \nat \}.$ From the construction of the superstructure
it follows that there exists some $m \in \omega$ such that $\nat \times T 
\subset X_0 \cup X_m.$ Hence $\power {\nat \times T} \subset \power {X_0 \cup 
X_m}= X_{m+1},$ where $m \geq 1.$ Since no atoms are in $ \nat \times T$,
the set $\nat \times T \subset X_m.$ For each $x \in H_1,\ T^x \in \power {\nat 
\times T} \subset X_{m+1}.$ Obviously$,$ $H_1 \in X_2.$ Hence$,$ we also have that
$H_1 \in X_{m+1}.$ Observe that for each $x \in H_1,\ w \in T^x \subset x 
\times T \subset \nat \times T \subset X_0 \cup X_m$ implies since no atoms 
are involved that $w \in X_m \subset X_{m+1}.$ Thus$,$ all of the objects being 
considered in an expression of the type $w \in T^x$ are all members of the set 
$X_{m+1}.$ Notice that in the formal language $w \in T^x$ is a 
2-place predicate replaceable  by $\exists y 
\exists z( y \in \nat \land 0\leq y\leq x \land z \in T \land (y,z) \in w) 
\land \forall y_1 \forall z_1 \forall x_1(y_1 \in  \nat \land z_1 \in T \land
x_1 \in T \land (y_1,z_1) \in w \land (y_1,x_1) \in w \to x_1 = z_1)).$
Suppose that you have a formula with the expression  ``$\land (w \in 
T^x)$'' as a subformula. Then  replace it by  `` $\land
(w \in X_m) \land (w \in T^x).$'' Now the explicit specially constructed 
formula usually used in the literature for such bounded formula is obtained by 
expanding the finite sequences of ``$\land$'' into the equivalent   forms
``$(\ \to \ (\ \ldots$'' since recall that for the propositional calculus
an expressing such as ${\rm P} \land {\rm Q} \to {\rm S}$ is equivalent to
${\rm P} \to ({\rm Q} \to {\rm S}).$  With these processes$,$ all of the formula 
that seem to have quantified 
variables with missing bounding objects can be modified into an equivalent 
bounded form. Further, there are equivalent formula such as  $\forall x(x \in C \to \exists y (y \in x \ldots))$ where $C$ is standard that express the requirement that the quantified variables vary over members of our superstructure. Also, ${\cal N}$ and $\hyper {\cal N}$ are closed under the basic set-theoretic operations. 
\pars
For each $f \in P,$ let $[f]$ denote the equivalence class in ${\cal E}$ 
containing $f.$ For each $ g \in \hyper P,$ let $[g]$ denote the equivalence 
class containing $g$ and determined by the partition $\hyper {\cal E}$ of 
$\hyper P.$\parm
{\bf Theorem 3.3.1.} {\sl For each $f \in P,$ it follows that
$\Hyper [f] = [\hyper f] = [{^\sigma f}] = [f].$}\pars
Proof. Let $f \in P.$ Then $ f \in [f] \in {\cal E}$ implies that 
$\hyper f \in \hyper P$ and $ \hyper f \in \Hyper [f] \in \hyper {\cal E}$ by 
properties of the *-map. Now the sentence
$$f \in P \to \exists !z(z \in {\cal E} \land f \in z) \eqno (3.3.4)$$
holds in ${\cal M};$ hence in $\hyper {{\cal M}}.$ Thus there is a unique set 
$A \in \hyper {{\cal E}}$ such that $\hyper f \in A \in \hyper {{\cal E}}.$ 
This set is denoted by $[\hyper f]$ since it contains $\hyper f,$ and $\hyper 
{{\cal E}}$ is a partition. The uniqueness implies that $\Hyper [f] = [\hyper 
f]$ and the finite nature of $f$ yields that $[ {^\sigma f}] = [f]$ under our 
conventions. \qed
{\bf Theorem 3.3.2.} {\sl There exists a purely subtle $d \in \Hyper {{\bf 
BP}} - {\bf BP}$ such that $\Hyper {\bf C_a}(\{d\}) \cap {\bf BP} =$ an 
infinite set and $\Hyper {\bf C_a}(\{d\}) \cap (\Hyper {{\bf BP}} - {\bf 
BP})=$ an infinite set.}\pars
Proof. Let  ``just'' be a member of BP and consider the sentence \parm
\line{\hfil $\forall x(x \in \nat\land x > 0 \to \exists y(y \in T^x 
\land \forall w ( w\in \nat \land$\hfil}\pars
\line{\hfil $0< w\leq x \to y(w) = {\rm very,}\sp \land y(0) = {\rm 
just})\land $\hfil }\pars
\line{ \hfil $\forall z(z \in \nat\land z \leq x \to \exists y_1(y_1 \in 
T^z \land \forall w_1(w_1 \in \nat \land$\hfil}\pars
\line{\hfil $0 < w_1 \leq z \to y_1(w_1) = {\rm very,}\sp\land y_1(0) = {\rm 
just} 
\land$
\hfil}\pars
\line{\qquad\hfil $[y_1] \in {\bf C_a}(\{[y]\})))))).$\hfil (3.3.5)}\parm
\noindent which holds in ${\cal M};$ hence$,$ in $\hyper {{\cal M}}.$ So$,$ let 
$\nu \in \hypernat - \nat.$ Then there exists some *-partial sequence $ f\in 
(\Hyper T)^\nu - P$ such that $f(w) = i({\rm very,}\sp)$ for each 
$0 < w \leq \nu$ and $f(0) = i({\rm just}).$ Also for each $ n \in \nat,\ 
n>0,$ there exists a partial sequence $f_n \in T^n$ such that for each $ n \in 
\nat,$ where $0< w \leq n,\ f_n(w) = i({\rm very,}\sp)$ and $ f_n(0) = i({\rm 
just}).$ Notice that if $n,m \in \nat - \{0\}$ and $ n \not= m,$ then 
$\Hyper [f_n] \not= \Hyper [f_m].$ Now for each $ n\in \nat,\ n > 0,\ 
\Hyper [f_n] \in {^\sigma{\bf BP}}.$ Application of Theorem 3.3.1 implies that
$\Hyper [f_n] =[\hyper f_n] = [f_n]$ and the above sentence yields that for 
each such $ n \in \nat,\ \Hyper [f_n] \in \Hyper {{\bf C_a}}(\{[f]\}).$ 
Consequently$,$ $\Hyper {{\bf C_a}}(\{[f]\}) \cap {\bf BP} =$ an infinite 
set.\pars
Consider the infinite set $R= \{\nu_n\mid \nu_n=\nu - n\land n \in 
\nat\land n\geq 1\} \subset \hypernat - \nat.$ By *-transform of the 
above$,$ for each $ n \in \nat,\ n \geq 1,$ there exists  some $g_n \in (\hyper 
T)^{\nu_n} - P$ such that $g_n(w) = i({\rm very,}\sp)$ for each $w  \in 
\hypernat;\ 0< w \leq v_n < \nu$ and $g_n(0) = i({\rm just}).$ Observe that if 
$n,m \in \nat$ and $n \not= m,$ then $\nu_n \not= \nu_m.$ Moreover$,$ the 
following sentence\parm
\line{\hfil $ \forall x\forall y(x\in \nat \land y \in \nat\land x> 
0\land  y> 0 \land x \not= y \to$\hfil}\pars
\line{\hfil$\forall w \forall w_1 \forall z \forall z_1(w \in T^x\land w_1 
\in T^y \land z \in \nat \land z_1 \in \nat\land$\hfil}\pars
\line{\hfil $0<z \leq x\land 0<z_1\leq y\land w(z) = {\rm very,}\sp 
\land w_1(z_1) = {\rm very,}\sp  \land$\hfil}\pars
\line{\qquad\qquad\hfil $ w(0) = {\rm just}\land w_1(0) = {\rm just} \to [w] 
\not= [w_1]))$\hfil (3.3.6)}\parm
\noindent holds in $\cal M;$ hence in $\hyper {{\cal M}}.$ By *-transfer$,$ if
$n,m \in \hypernat -\nat,\ n,m > 0, n \not= m,$ then $[g_n] \not= [g_m].$ Since for each 
such $n \in \hypernat -\nat,\ [g_n]\in \Hyper {{\bf C_a}}(\{[f]\}),$ it follows that
$\Hyper {{\bf C_a}}(\{[f]\}) \cap (\Hyper {{\bf BP}} - {\bf BP}) =$ an infinite set. 
 \qed
{\bf Corollary 3.3.2.1} {\sl There exists a purely subtle $d \in \hyper {{\cal 
E}} - {\cal E}$ such that $\Hyper {{\bf C_a}}(\{d\}) \cap {\cal E} =$ an 
infinite set and $\Hyper {\bf C_a}(\{d\}) \cap (\hyper {{\cal E}} - {\cal 
E})=$ an infinite set.}\parm
Of course$,$ in order to apply these results to descriptions that involve 
members of BP an {interpretation procedure is required.}  We have previously 
discussed  the intuitive interpretation for any $[f] \in \hyper {{\cal E}} - 
{\cal E},$ where the range of $f \in \hyper P$ is a subset of $i[{\cal W}].$ 
Hence$,$ if $A \in {^\sigma {\cal E}},$ then $ A = \Hyper [g]=[g],$ where 
$[g] \in {\cal E}.$ Thus in the usual manner$,$ first interpret $\Hyper [g]$ to 
be $[g] \in {\cal E}$ and then proceed to the second step and interpret
$[g]$ by selecting any $f \in [g]$ and applying our previously discussed 
inverse procedure. Clearly$,$ this interpretation method is a one-to-one 
correspondence from a subset of $\hyper {{\cal E}}$ into ${\cal W}.$\pars
The concept of adjective deduction$,$ which is obviously isomorphic to a subsystem 
of ordinary propositional deduction$,$ was originally  introduced to give a 
{measure of the strength of various behavioral properties}. These intuitive 
strengths may not be codifiable by a numerical measure. Thus$,$ intuitively$,$ 
 ``very,$\sp$very,$\sp$bold'' is  a stronger concept than  ``very,$\sp$bold''.
The exact same process can be applied to physical concepts as well. Even 
though it may not be possible to measure the combined strengths of all of the 
intuitive forces that my be altering the appearance of a physical 
entity such as a thunderhead$,$ the term ``very,$\sp$'' could be replaced by  
other terms such as ``greater,$\sp$'' or ``weaker,$\sp$'' coupled with terms 
such as  ``force'' and the like. The same type of $\vdash_a$ analysis  
would follow.\pars
With respect to the above remarks$,$ later in this book$,$ we consider the 
reasoning 
process called simply  ${\rm S},$ which is an {axiomatically presented subsystem of 
propositional deduction}. The process ${\rm S}$ is closely associated with 
adjective reasoning$,$ if the set BP is constructed in a different manner and 
from different objects. One of the minor problems with these constructions is 
their relation to formal languages and {the use of parentheses} within such 
formal languages. Another illustration of the use of these nonstandard 
methods that does parallel Robinson's original work along this line requires 
BP to be formally embedded into a propositional language with the {insertion 
and removal of such parentheses}.\par
\bigskip
\leftline{\bf 3.4 Propositional Reasoning}
\medskip
Let ${\rm S_0}$ denote the consequence operator determined by the usual propositional 
deduction $\vdash$ as defined on say pages 108-109 of [7] (i.e. Group A1 
deduction).  
Technically$,$ ${\rm BP}_0 \not\subset {\rm L}_0$ since  ${\rm BP}_0$ is 
constructed without use of parentheses. Let ${\rm L}^\prime = {\rm L}_0 \cup {\rm BP}_0.$ 
Extend $\vdash$ and $\rm S_0$ in the obvious manner. Let ${\rm A} \in {\rm BP}_0 -
j[{\rm B}].$ Then ${\rm A} = {\rm V} \land \ldots \land {\rm V} \land {\rm 
b},$ where ${\rm b} \in j[{\rm B}]$ and there are $n \geq 1$ connectives 
$\land.$ We now consider inserting parentheses in the following manner called 
the {{\it insertion procedure.}} (1) Moving from left to right put a  ``('' 
before each V$,$ keeping count of the number of   ``('' so placed. (2) 
Place the same number of  ``)'' after the  ``b'' as your count in step (1). 
Denote this new symbol by ${\rm A}_(.$ Note that ${\rm A}_( \in {\rm L}_0.$\pars
{\it Example.} Suppose that you are given ${\rm A} = {\rm V} \land {\rm V} 
\land {\rm V} \land {\rm V} \land {\rm b}.$ Then ${\rm A}_( = 
({\rm V} \land ({\rm V} 
\land ({\rm V} \land ({\rm V} \land {\rm b})))).$  \pars
This process of considering a method of inserting parentheses and doing it in 
an ordered effective manner is no more complex and no less effective than 
Kleene's concept of  ``closure with respect to (just) 
${\rm x}_1,\ldots,{\rm x}_q$'' on page 105 of [8]. Now to define in the 
obvious manner $\vdash^\prime \subset F({\rm L}^\prime) \times {\rm L}^\prime.$ 
First$,$ consider ${\rm BP}_{0(} = \{ \rm x_( \mid \rm x \in {\rm BP}_0 - 
j[{\rm B}]\} \cup j[{\rm B}].$ Then ${\rm BP}_{0(} \subset {\rm L}_0.$ 
For $\rm F \in F({\rm L}^\prime),$ consider $\rm F_( = \{\rm x_( \mid
\rm x \in \rm F \cap ({\rm BP}_0)\} \cup (\rm F - ({\rm BP}_0 - j[{\rm B}])) \subset {\rm 
L}_0.$ Then (i) if $\rm F_( \vdash {\rm B} \in {\rm L}_0,$ define
$\rm F \vdash^\prime {\rm B}.$ (ii) If $\rm F_( \vdash {\rm B}_1 \in {\rm L}_0$ 
and D is ${\rm B}_1$ with {\it all of the parentheses removed} and ${\rm D} 
\in {\rm BP}_0,$ then let $\rm F \vdash^\prime {\rm D}.$ (iii) Finally$,$ remove 
superfluous parenthesis if you wish [7$,$ p. 74]. Only the procedure in this 
paragraph is to be used to obtain a ${\rm B} \in {\rm L}^\prime$ such that 
$\rm F \vdash {\rm B}.$\pars
Obviously$,$ $\vdash^\prime$ is closely related to $\vdash,$ since it is well 
known that for ${\rm A,\ B,\ C} \in {\rm L}_0, \vdash ({\rm A} \land ({\rm B} 
\land {\rm C})) \iff (({\rm A }\land {\rm B}) \land {\rm C}).$ By an abuse of 
notation we often write $\vdash$ for $\vdash^\prime,$ $\rm S_0$ for ${\rm 
S_0}^\prime$ and ${\rm L}_0$ for ${\rm L}^\prime$ as well as 
suppressing parentheses insertion and removal 
for the elements of ${\rm BP}_0.$ The next result follows from the fact that 
for ${\rm A,\ B} \in {\rm L}_0,\ {\rm A} \land {\rm B} \vdash {\rm B}.$\parm
{\bf Theorem 3.4.1} {\sl There exists a purely subtle  $d \in \Hyper {{\bf 
BP}_0} - {\bf BP}_0$ such that $\Hyper {{\bf S_0}}(\{d\}) \cap {\bf BP_0} = $ an 
infinite set and also $\Hyper {{\bf S_0}}(\{d\}) \cap (\Hyper {{\bf 
BP_0}} - {\bf BP_0}) = $ an infinite set.}\pars
Proof. Make the following changes in the formal first-order sentences 
explicitly given in the proof of Theorem 3.3.2. First$,$ let $ {\rm b} \in j[{\rm 
B}].$ Now for the every  ``just'' substitute the symbol  ``b''. Then for every
``very,$\sp$'' string substitute the symbols ``V $\land$\ ''.  With these 
substitutions made$,$ the proof is exactly as for Theorem 
3.3.2. \qed
{\bf Corollary 3.4.1.1} {\sl There exists a purely subtle  $d \in \Hyper {{\bf 
L_0}} - {\bf L_0}$ such that $\Hyper {{\bf S_0}}(\{d\}) \cap {\bf L_0} = $ an 
infinite set and also $\Hyper {{\bf S_0}}(\{d\}) \cap (\Hyper {{\bf 
L_0}} - {\bf L_0}) = $ an infinite set.}\parm
[{\it Remark:} The above theorems for the propositional consequence operator 
also hold for 
the  consequence operator ${\rm S}$ and other such variations discussed 
later in this book.]\pars 
It is easily shown that for a propositional formula$,$ say A$,$ that 
${\rm A} \vdash^\prime {\rm A} \land \cdots \land {\rm A},$ with any 
$n \in \nat$ number of connectives $ \land.$ \parm
{\bf Theorem 3.4.2} {\sl For any $q \in {\bf L_0}$ it follows that 
$\Hyper {{\bf S_0}}(\{q\}) \cap {\bf L_0}= $ an infinite set and 
$\Hyper {\bf S_0}(\{q\}) \cap (\Hyper {{\bf L_0}} - {\bf L_0}) =$ and infinite 
set.}\pars
Proof. Let $ q = [f],\  f \in T^0,\  f(0) = i({\rm A}),\ {\rm A} \in {\rm 
L}_0.$ The sentence \parm
\line{\hfil$\forall z( z \in \nat \land z > 0 \to \exists y(y \in T^z \land 
\forall w(w \in \nat \land 0 < w \leq z \to$\hfil}\pars
\line{\hfil$y(w)= i({\rm A}\ \land\ ) \land y(0) = i({\rm A}) \land 
[y] \in {\bf S_0}(\{q\}))))$\hfil (3.4.1)}\parm
\noindent holds in ${\cal M}$; hence in $\hyper {\cal M}.$ Now proceed in the same 
manner as in the proof of Theorem 3.3.2$,$ making the obvious changes$,$ starting 
with the statement$,$ ``Also for each $ n \in \nat,\ n > 0,$ . . . . '' This 
completes the proof. \qed\medskip

\leftline{\bf 3.5 Modus Ponens Reasoning}
\medskip
The reasoning termed {{\it Modus Ponens}}  (MP) is$,$ of course$,$ the major step 
in propositional deduction. One can$,$ however$,$ get more basic than $\rm S_0$ and define 
MP reasoning to produce a subsystem of $\rm S_0$ in the following (intuitive) manner.
Simply let MP be the same deduction process as determines $\rm S_0$ but with no axiom 
schemata. Use the symbol MP to represent the consequence operator obtained 
from this process. Then$,$ for each $\rm A \subset {\rm L}_0,$ it follows that $\rm A \subset
{\rm MP}(\rm A) \subset {\rm S_0}(\rm A).$ Thus for each internal $ B \subset \Hyper 
{{\bf L_0}}, \  B \subset \Hyper {{\bf MP}}( B) \subset \Hyper {{\bf S_0}}( B)
\subset \Hyper {{\bf L_0}}.$ \pars
Besides applying MP to ${\rm L}_0,$ it is straightforward to apply it to 
certain meaningfully constructed collections of intuitive readable sentences.
For example$,$ consider the set of symbols $\rm B_1 = \{ {\rm If}\sp{\rm 
perfect,}\sp{\rm then}\sp\rm x. \mid \rm x \in {\rm BP}\}.$
Now apply MP to any
finite subset of ${\rm BP} \cup {\rm B_1} \cup \{{\rm perfect}\}.$ Clearly$,$ we 
can associate MP deduction  formally to ${\rm BP}_0$ in a meaningful way. 
Simply let $\rm c \in \{{\rm P}_i \mid i \in \nat \} - (\{{\rm V}\} \cup j[{\rm 
B}])$ and ${\rm B}^{\prime\prime} = \{\rm c \to \rm x \mid \rm x \in {\rm BP}_0 \},$ etc. 
We leave to the reader the simple consequences of MP deduction in this 
case.\par
\bigskip
\leftline{\bf 3.6 Predicate Deduction}
\medskip
In this section$,$ predicate deduction in ${\rm L}_1,$ say as defined by Kleene 
on page 82 page [7]$,$ is briefly discussed relative to lengths of formal 
proofs. Robinson mentions [15$,$ p. 25]$,$ what is well known from G\"odel's work,
that using formal predicate deduction there is for each $ n \in \nat$ a 
readable sentence in ${\rm L}_1$ that is provable as a theorem  from the empty 
set of hypotheses$,$ but {requires n or more steps}.\pars
Let ${\rm S}_1$ denote the operator determined by predicate deduction with 
respect to ${\rm L}_1.$ Then  ${\rm S}_1(\emptyset)$ is the set of all 
provable formula (i.e. theorems). Of course$,$ all the properties of  ${\rm 
S}_1$ are now referred to ${\cal E}.$ Hence there exists a relation $R_{L_1}
\subset \nat \times {\bf S_1}(\emptyset)$ with the property that
$(x,y) \in R_{L_1}$ iff $ y \in {\bf S_1}(\emptyset),\ y \in {\bf L_1},\ x \in 
\nat$ and $x =$ the length of a formal proof that yields $y \in {\bf 
S_1}(\emptyset).$\parm
{\bf Theorem 3.6.1} {\sl For each $\nu \in \hypernat - \nat$ there is a subtle
$d \in \Hyper {{\bf L_1}},\ d \in {\bf S_1}(\emptyset)$ and for each $\lambda
\in \hypernat,\ \lambda <\nu,\ (\lambda ,d) \notin \hyper {R_{L_1}}.$ 
Moreover$,$ there exists some $\nu_0 \in \hypernat - \nat$ such that 
$\nu_0 \geq \nu$ and $(\nu_0,d) \in \hyper {R_{L_1}}.$}\pars
Proof. As stated above the sentence\parm
\line{\hfil $\forall x(x \in \nat \to \exists y(y \in {\bf S_1}(\emptyset) 
\land \forall w(w \in \nat\land 0\leq w\leq x \to $\hfil}\pars
\line{(3.6.1)\hfil $(w,y) \notin R_{L_1}) \land \exists z(z \in \nat \land z \geq 
x \land (z,y) \in R_{L_1})))$\hfil }\parm
\noindent holds in $\cal M$; hence in $\hyper {{\cal M}}.$ The result follows 
by *-transfer.\parm
We now investigate a little  more fully what is meant by the ``{length of a 
formal proof}.'' There exists a partial sequence $f^\prime$ of elements of 
${\rm L}_1$ such that the domain $D(f^\prime) = [0,n],\ n \in \nat,$ [rather 
than $ n \in \nat$]$,$ the range of 
$f^\prime = 
{\rm Rn}(f^\prime) \subset {\rm L}_1$ and $f^\prime (n) = \rm A \in {\rm 
S}_1(\emptyset)$ and the length of the formal proof that yields $\rm A \in {\rm 
S}_1(\emptyset)$ is $ n + 1.$ Of course$,$ $f^\prime$ actually gives the 
elements of ${\rm L}_1$ that appear in such a specific formal proof. Now 
relate this intuitive partial sequence $f^\prime$  to a corresponding partial 
sequence in $P$ in the following manner. Let $(x,y) \in f_{L_1}$ iff $
y = i(\rm w)$ and $(x,\rm w) \in f^\prime.$ Denote by $P_{L_1} \subset P$ the set of 
all such length of proof sequences. Then $(x,y) \in \hyper {R_{L_1}}$ iff
$x \in \hypernat,\ y \in \Hyper {{\bf S_1}}(\emptyset)$ and there is some 
$f \in \hyper {P_{L_1}}$ such that $y = [f(x)].$ By *-transfer the 
{hyperlength of the proof} would be $x + 1.$ The set $P_{L_1}$ may be used to 
characterize the concept intuitively associated with the proof length for 
objects in ${\rm S}_1(\emptyset).$ These sequences have other properties as 
well but these will not be considered in this investigation. 
With this in mind$,$ then  $\hyper {P_{L_1}}$ represents the {subtle concept of 
proof length} for elements in $\Hyper {{\bf S_1}}(\emptyset).$ It's the proof 
length concept we employ in one application of the results from this chapter. 
Theorem 3.6.1 can now be stated in an alternate form.\parm
{\bf Theorem 3.6.2} {\sl For each $ \nu \in \hypernat - \nat$ there is a 
subtle $ d  \in \Hyper {{\bf S_1}}(\emptyset)$ such that for each $ \lambda 
\in \hypernat -\nat,\  \lambda < \nu,$ there does not exist some $ g \in 
\hyper {P_{L_1}}$ such that $d = [g(\lambda)].$ Moreover$,$ there exists some 
$\nu_0 \in \hypernat -\nat$ such that $ \nu_0 \geq \nu$ and some $ f \in 
\hyper {P_{L_1}}$ such that $ d = [f(\nu_0)].$}\parm
{\bf Corollary 3.6.2.1} {\sl There exists $ d \in \Hyper {{\bf 
S_1}}(\emptyset),\ f \in \Hyper {P_{L_1}}$ and $ \nu \in \hypernat -\nat$ 
such that $ d = [f(\nu)]$ and for each $ g \in \Hyper {P_{L_1}}$ and 
each $\nu^\prime \in \hypernat$ such that $d = [g(\nu^\prime)],\ \vert 
D(g)\vert
\geq 2^\omega.$}\parm
[{\it Remarks}: It should be apparent to the reader that statements that hold in 
${\cal M}$ relative to consequence operators or deductive processes are 
obtained from the corresponding intuitive reasoning processes by application 
of $\theta.$ The proofs that these statements hold in the intuitive case are 
straightforward or obvious$,$ and are omitted in all cases. Also$,$ you might 
wonder about the term ``ultralogics'' since it has not been specifically 
defined as yet. The term is reserved for various special subtle consequence 
operators to be used in various cosmogony investigations that will be 
discussed later in this book.] 
\parm

(Adjective reasoning can also be determined by a general logic system.  (Herrmann, R. A. General Logic-Systems and Finite Consequence Operators, Logica Universalis, 1(2006):201-208 (Partial paper at  http://arxiv.org/abs/\hfil\break math/05012559).) Let $\rm x \in BP$ have $\r n > 0$  $\rm very,\sp$ strings to the left of a $\rm c \in C_0$ (page 26). Then a rule of inference $\rm R_x$ for $x$ is constructed by reduction as follows: remove one $\rm very,\sp$ from x. Write the result as $\rm x_1$. Then let $\rm (x,x_1)\in R_x$. Continue this finite reduction until $\rm c in C_0$ is obtained. Hence, the last member of $\rm R_x$ so constructed is $\rm (x,c).$ By definition, the set of all such finite binary relations $\rm R_x$ obtained for each such x yields a general logic-system. (This is not a unique construction.) From this system, the corresponding consequence operator $\rm C_a$ is obtained.)
\vfil\eject
\centerline{NOTES}

\vfil\eject.
\centerline{\bf 4. SPECIAL DEDUCTIVE PROCESSES}\par
\bigskip
\leftline{\bf 4.1 Introduction.}\par 
\medskip
There are certain words that intuitively denote an {upper [resp. lower] bound 
to such concepts  as ``stronger'' [resp. ``weaker'']}. With respect to certain 
philosophic studies$,$ one such concept is the notion of ``{perfect}'' when 
associated with a language like BP. In what follows$,$ this ``perfect'' 
associated 
with BP is used as a prototype for these other cases. Two types of deductive 
processes associated with this prototype will be introduced$,$ a very trivial 
one followed by a much more interesting and significant procedure. 
\bigskip
\leftline{\bf 4.2 Reasoning From the Perfect Type W}
\medskip
First$,$ an  
intuitive extension of BP is defined. Let ${\rm BPC|} = {\rm BP} \cup \{{\rm 
perfect}\}$ and for convenience denote the readable 
string ``perfect'' by the single c. Now we {define type W reasoning from the 
perfect} by considering an intuitively defined operator$,$ $\Pi_W,$ 
from 
$\power {{\rm BPC|}}$ into  $\power {{\rm BPC|}}.$\pars
For any finite ${\rm F} \subset {\rm BPC|}$:\parm
 (i) if ${\rm c} \in {\rm F},$ then $\Pi_W({\rm F}) = {\rm BPC|};$\pars
 (ii) if ${\rm c} \notin {\rm F},$ then $\Pi_W({\rm F}) = {\rm F};$\pars                                                          
(iii) and for arbitrary ${\rm A} \subset {\rm BPC|},$ let\parm 
\line{\hfil $\Pi_W ({\rm A}) = \bigcup \{ \Pi_W({\rm F}) \mid {\rm F} \in F({\rm A}) \}$\hfil}\parm
{\bf Theorem 4.2.1} { \sl The map $\Pi_W \colon \power {{\rm BPC|}} \to \power 
{{\rm BPC|}}$ is a consequence operator.}\pars
Proof. Let ${\rm A} \subset {\rm BPC|}.$ Clearly$,$ axiom (4) holds by the definition. 
Let ${\rm a} \in {\rm A}.$ Then $\{{\rm a}\} \in F({\rm A}).$ Now if ${\rm a} \not= {\rm 
c},$ then $\Pi_W(\{{\rm a}\}) = \{{\rm a}\}.$ If ${\rm a} = {\rm c},$ then 
$\Pi_W (\{{ a}\}) = {\rm BPC|}.$ In these two cases$,$ (iii) of the definition 
yields that ${\rm a} \in \Pi_W({\rm A}).$ Thus$,$ even when $ {\rm A} = 
\emptyset,$ it follows that ${\rm A} \subset \Pi_W({\rm A}) \subset {\rm BPC|}$ and axiom 
(2) holds. \pars
Since axiom (4) holds and ${\rm A} \subset \Pi_W({\rm A}),$ it follows that $\Pi_W({\rm A}) 
\subset \Pi_W(\Pi_W({\rm A})).$ Now either $\Pi_W({\rm A}) = {\rm A};$ in which case
$ \Pi_W(\Pi_W({\rm A})) = \Pi_W({\rm A}) ={\rm A}$ or $\Pi_W({\rm A}) = {\rm BPC|}$; in which case
$ \Pi_W(\Pi_W({\rm A}))= \Pi_W({\rm BPC|}) = {\rm BPC|}= \Pi_W({\rm A}).$ Thus axiom (3) 
holds and this completes the proof. \qed
Recall that $T =i[{\cal W}]$ and if ${\rm w} \in {\cal W},$ then $f_w \in P$ 
denotes the partial sequence 
which is an element of $T^0$ and $f_w(0) = i(w),\ {\bf w} = 
[f_w].$
Also$,$ due to their finitary character$,$ each $ x \in {\cal E}$ is often 
identified
\noindent 
 with $\hyper x \in {^\sigma {\cal E}}.$ \parm
{\bf Theorem 4.2.2} {\sl For each internal $B \subset {\bf BPC|}$ if
${\bf c} = [\hyper f_{c}] = [f_c],$ then $\Hyper \Pi_W(B) = \Hyper {{\bf BPC|}} 
= \Hyper {{\bf BP}} \cup \{{\bf c}\} = \Hyper {{\bf BP}} \cup \{[f_c]\},$ where 
$[f_c] \in {^\sigma {\cal E}}$ and ${^\sigma {\cal E}}= {\cal E}$ under the 
basic identification of $^\sigma \nat$ with $\nat.$}\pars
Proof. Simply consider the sentence\parm
\centerline{}
\vfil
\eject
\line{(4.2.1)\hfil $\forall x(x \in \power {{\bf BPC|}} \land {\bf c} \in x \to
{\bf \Pi_W}(\{[f_c]\}) = {\bf BPC|})$\hfil }\parm
\noindent that holds in ${\cal M};$ hence in $\hyper {{\cal M}}.$ The result 
follows by *-transfer. \qed
{\bf Corollary 4.2.2.1} {\sl The set $\Hyper {{\bf \Pi_W}}(\{[f_c]\}) = 
\Hyper {{\bf \Pi_W}}(\{[\hyper f_c]\})= \Hyper {{\bf BPC|}}.$}\parm
{\bf Corollary 4.2.2.2} {\sl For each  $b \in i[{\rm B}]$ and each 
$\nu \in \hypernat - \nat$ there exists a subtle $f^b \in 
(\Hyper T)^{\nu} - P$ such that for each $x \in \hypernat,$ where 
$0< x \leq \nu,\ f^b(x) = i({\rm very,}\sp)$ and $f^b(0) = b.$ Moreover,
$[f^b] \in \Hyper {{\bf \Pi_W}}(\{[f_c]\}).$}\parm
For each $b \in i[{\rm B}]$ and a fixed $\nu \in \hypernat - \nat,$ apply the 
axiom of choice and let $f^b$ denote one of the subtle objects that exists by 
Corollary 4.2.2.2 and satisfies the stated properties. Since B is finite$,$  
the set $F_\nu = \{[f^b]\mid b \in i[{\rm B}]\}$ is internal. The next result is 
obvious.\parm
{\bf Theorem 4.2.3} {\sl For each $\nu \in \hypernat - \nat,$ internal
$F_\nu \subset \Hyper {\bf \Pi_W}(\{[f_c])\}).$}\parm
Observe that there exist$,$ at least$,$  $2^\omega$ distinct $F_\nu$ sets.\par
\bigskip
\leftline{\bf 4.3 Strong Reasoning From the Prefect}
\medskip

For the second type of reasoning from the perfect$,$ our attention will be 
restricted to 
${\rm L}^\prime= {\rm L_0 \cup BP_0}$ and the set ${\rm BP}_0$ that 
bijectively corresponds to BP. 
Let a specific ${\rm c} \in \{{\rm P}_i 
\mid i \in \nat\} - (\{{\rm V}\} \cup j[{\rm B}]).$ Correspond c to the 
readable sentence  ``prefect.'' Let \parm
\line{\hfil $ {\rm C}_( = \{({\rm c} \to {\rm x}_() \mid {\rm x} \in {\rm 
BP_0- j[B]}\} 
\bigcup \{({\rm c} \to {\rm x})\mid \rm x \in j[{\rm B}]\} \bigcup \{{\rm c}\}.$\hfil}\smallskip
 \line{\hfil $ {\rm C} = \{{\rm c} \to {\rm x} \mid {\rm x} \in {\rm 
BP_0-j[B]}\} 
\bigcup \{{\rm c} \to {\rm x}\mid \rm x \in j[{\rm B}]\} \bigcup \{{\rm c}\}.$\hfil}
 \parm
 Why do we go through the following exercise of inserting and removing 
parentheses so as to conform more closely to the formula of a formal language?
The basic reason is related to some of the results later in this book that 
refer to counting of symbols by means of the partial sequences. Clearly$,$ 
parenthesis insertion does correspond to the increase strength idea of 
adjective reasoning$,$ as does the ordering of the ${\rm very,}\sp$ symbols by the 
partial sequences. However$,$ in certain deductive processes$,$ all of the axioms 
for the propositional logic are not used. Hence even though it is certainly 
of no importance$,$ due to equivalence$,$ when all of the usual axioms are used
to write a formal $({\rm V} \land ({\rm V} \to {\rm b}))$ as ${\rm V} 
\land {\rm V} \to {\rm b},$ it may not be possible to establish this 
equivalence for these restrictive deductive processes. The process we now 
outline simply removes this formal difficulty at the cost of a more involved 
finitary process.\pars
Let ${\rm BPC}_0 = {\rm BP}_0 \cup {\rm C}$ and ${\rm BPC}_{0(} = {\rm BP}_{0(}
\cup {\rm C}_(.$ The axioms are elements of the set ${\rm Ax} = \{({\rm V} \land 
{\rm x}) \to {\rm x} \mid {\rm x} \in {\rm BP}_{0(}\}$ with the suppression of the outer most 
parentheses for simplicity in application of MP. Let $\vdash_\pi$ denote
 ordinary propositional deduction but {\it only 
using the axiom set {\rm Ax} and only formula from the set} ${\rm Ax} \cup 
{\rm BPC}_{0(}$ {\it in the steps of any proof.} For a specific ${\rm A}_( 
\in {\rm BP}_{0(},$ 
let A be the element of ${\rm BP}_0$ formed by removing all parentheses from 
${\rm A}_(.$ Define $\Pi \colon \power{{\rm BPC}_0} \to 
\power {{\rm BPC}_0}$ as follows: Let finite ${\rm F} \subset {\rm BPC}_0 
\subset {\rm L}^\prime$ and 
${\rm F}_( = ({\rm BP_0}\cap {\rm F})_( \cup ({\rm C} \cap {\rm F})_(.$ Then for any 
${\rm D} \subset 
{\rm BPC}_0, \  {\rm X} \in \Pi ({\rm D}) $ iff there exists a finite ${\rm F} 
\subset {\rm D}$ 
and ${\rm A} \in {\rm BPC}_{0(}$ such that ${\rm F}_( \vdash_\pi {\rm A}_($ and if 
${\rm A}_( \in {\rm C}_(,$ then ${\rm X} = {\rm A}$ or if ${\rm A}_( \in {\rm 
BP}_{0(},$ then ${\rm A} = {\rm X}.$ It is not difficult to show that $\Pi$ 
is a consequence operator since it is the restriction of formal $\vdash_\pi$ to ${\rm BPC}_{0(}$ and $\Pi$ is called {{\it strong reasoning from the 
perfect.}}\pars

By way of a reminder$,$ *-transfer and the fact that ${\rm E}_0 = j[{\rm B}]$ and 
$i[{\rm E}_0]$ are finite imply that $ f \in \Hyper {{\bf BP_0}} - \Hyper 
{{\bf E_0}} $ iff there exists some $ n  \in \hypernat, n > 0,$ and $b \in 
i[{\rm E_0}]$ such that $f(w) = i({\rm V} \land)$ and $f(0) = b$ for each $ w 
\in \hypernat$ such that $0 < w \leq n.$ Further$,$ for each $ b  \in i[{\rm 
E}_0]$ and each $ b \in  i[{\rm E}_0]$ and each $ n \in \hypernat, n > 0,$ 
there exists a $f \in \Hyper {{\bf BP_0}} - \Hyper {{\bf E_0}}$ such that 
for each $ w \in \hypernat,$ where $0 < w\leq n,\ f(w) = i({\rm V} \land )$ 
and $f(0) = b.$ Notice also that $[g] \in \Hyper {{\bf E_0}}$ iff 
$ g \in (\Hyper T)^0,\ g(0) = b \in i[{\rm E}_0]$ and 
$\vert [g] \vert = 1.$ The set $\Hyper {{\bf E_0}}$ being nonempty and finite implies that
$\Hyper {{\bf E_0}} = \{ \Hyper [g_1], \ldots, \Hyper [g_n]\} = 
{^\sigma{\bf E_0}} = {\bf E_0}.$ Finally$,$ due to the identification of each 
specific $\Hyper i({\rm e})$ with $i({\rm e}),$ it follows that $\hyper g = g,$ where$,$ as 
usual$,$ this follows from the finitary character of each equivalence class. \pars
We know that for fixed $ n  \in \hypernat$  and any $ {b} \in i[{\rm E}_0]$ 
there 
exists a unique $g^{b} \in (\Hyper T)^n$ such 
that $g^{b}(0) = {b}$ and if 
$ w \in \hypernat$ and $0 < w\leq n,$ then $g^{b}(w) = i({\rm V}\land).$
Let $G_n = \{[g^{b}]\mid {b} \in i[{\rm E}_0]\}$ for $n \in \hypernat - \{0\}.$
Now $G_n$ has the same general properties as the 
previously defined set $F_n.$ In particular$,$ each $G_n$ is internal and if $ 
\nu \in \nat_\infty,$ then  $G_\nu$ is purely subtle. \parm
{\bf Theorem 4.3.1.} { \sl If $n \in \hypernat - \{0\},\ m \in \hypernat - 
\{0\}$ are such that $0 < m\leq n$ and $ f \in (\Hyper T)^m$ has the property 
that for each $w \in \hypernat,$ where $0 < w \leq m;\ f(w) = i({\rm V}\land)$ 
and $f(0) = {b}  \in i[{\rm E}_0],$ then $[f] \in \Hyper {\bf \Pi}(G_n).$ 
Moreover$,$ if $[g] \in i[{\rm E}_0],$ then $[g] \in \Hyper {\bf 
\Pi}(G_n).$}\pars
Proof. First$,$ since $G_n$ is internal$,$ $ n \in \hypernat - \{0\},$ and 
$\Hyper {\bf BP_0} \subset \Hyper {\bf BP_0} \cup \hyper {\bf C} = \Hyper {\bf 
BPC_0}$ and $G_n \subset \Hyper {\bf BP_0},$ it follows that $G_n$ is in the 
domain of $\Hyper {\bf \Pi}.$ \pars
Let ${\rm A} \in {\rm BP}_0,\ {\rm A} = {\rm V}\land \cdots \land {\rm V}\land {\rm 
b}, {\rm b} \in {\rm E}_0$ and that are $ n \geq 1,\ ( n \in \nat)$ 
connectives $\land.$ We prove by induction that for any $ n \in \nat$ such 
that $ 0 < w \leq n,$ the symbol string ${\rm B} = {\rm V}\land \cdots \land {\rm V}\land {\rm 
b}$ with $w \geq 1$ connectives $\land$ or ${\rm B}={\rm b}$ has the property 
that ${\rm A}_( \vdash_\pi {\rm B}_(,$ where if ${\rm B}= {\rm b},$ then ${\rm B}_( = 
{\rm b}.$\pars
Case 1. Let $n = 1.$ Then ${\rm A} = {\rm V}\land {\rm b}.$ Consider
${\rm A}_( = ({\rm A} \land {\rm b}).$ The following is a proof that 
${\rm A}_( \vdash_\pi {\rm b}.$ (i) $({\rm V} \land {\rm b}),$ (ii)
$({\rm V} \land {\rm b}) \to {\rm b},$ 
(iii) b. A proof composed of step (i) only yields the trivial result that 
${\rm A}_( \vdash_\pi {\rm A}_).$\pars
Case $(n+1).$ Suppose that the result holds for $n,$ and ${\rm A} = {\rm V}\land
\cdots \land {\rm V} \land {\rm b}$ has $n+1$ connectives $\land.$ The formula
$\rm A_( = (\rm V\land (\rm A\cdots (\rm V \land \rm b)\cdots )).$  Let $\rm B_( =
(\rm V\land (\rm A\cdots (\rm V \land \rm b)\cdots ))$ have $n$ connectives $\land.$ 
The following is a proof that $\rm A_( \vdash_\pi \rm B_(.$ (i) $\rm A_( =
(\rm V \land \rm B_(),$ (ii) $(\rm V \land \rm B_() \to \rm B_(,$ (iii) $\rm B_(.$ 
Thus $\rm A_( \vdash_\pi \rm B_(.$ From the induction hypothesis$,$ $\rm B_( 
\vdash_\pi \rm E_(,$ where $\rm E_($ has $w$ connectives $\land$ such that 
$0 < w\leq n,$ or $\rm E_( = \rm b.$ Since $\vdash_\pi$ is transitive$,$ it 
follows that $\rm A_( \vdash \rm E_(.$ The trivial proof using step (i) yields 
that $\rm A_( \vdash_\pi \rm A_(,$ and the basic result above follows by 
induction. \pars
We have shown that for each $b \in i[\rm E_0],$ the following sentences 
\parm
\line{\hfil $\forall y \forall x\forall w\forall z ((x \in \nat)\land(x 
>0)\land (y \in \nat) \land (0<y\leq x) \land $\hfil}\pars
\line{\hfil $( w \in T^x) \land \forall w_1((w_1 \in \nat)\land (0<w_1\leq x) 
\to (w(w_1) = i({\rm V} \land))\land$ \hfil}\pars
\line{\hfil $(w(0)= b))\land ( z \in T^y) \land \forall z_1 ((z_1 \in 
\nat)\land (0 < z_1 \leq y) \to$ \hfil}\pars
\line{(4.3.1)\hfil $z(z_1) = i({\rm V} \land)) \land (z(0) = b)\to ([z] \in {\bf
\Pi}(\{[x]\}))).$\hfil} \parm
\line{ \hfil $\forall x \forall y\forall w\forall z((x\in \nat)\land(x>0)\land
(w \in T^x) \land \forall w_1((w_1 \in \nat) \land$ \hfil}\pars
\line{\hfil $(0 < w_1 \leq x) \to w(w_1) = i({\rm V} \land)) \land (w(0) = 
b)\land $\hfil}\pars
\line{(4.3.2)\hfil $(z \in T^0) \land (z(0) = b) \to ([z] \in {\bf \Pi}
(\{[x]\})))$\hfil}\parm
\noindent hold in $\cal M,$ hence in $\hyper {\cal M}.$ Since the singleton 
subsets  of $G_n$ are *-finite$,$ it follows that $\bigcup \{\Hyper {\bf \Pi}
(\{[g^{b}]\})\mid b \in i[{\rm E_0}]\} \subset \Hyper {\bf \Pi}(G_n)$ by *-
transfer of Axiom (4). Let $n \in \hypernat-\{0\},\ f \in (\Hyper T)^m,\ 
0< m \leq n,\ f(w) = i(\rm V \land)$ for each $w \in \nat,\ 0 < w \leq m$ and 
$f(0) = b.$ Then by *-transfer of sentence (4.3.1)$,$ we have that
$$ [f] \in \Hyper {\bf \Pi}(\{[g^{b}]\}) \subset \Hyper
{\bf \Pi}(G_n). \leqno (4.3.3)$$ \par
For $g\in T_0$ such that $g(0) = b,$ *-transfer of sentence (4.3.2) yields
$\Hyper [g] = [\hyper g] = [g] \in \Hyper {\bf \Pi}(\{[g^{b}]\}) \subset 
\Hyper {\bf \Pi}(G_n)$ and this completes the proof. \qed
{\bf Corollary 4.3.1.1} {\sl If $\nu \in \nat_\infty,$ then ${\bf BP_0} \cap
G_\nu = \emptyset,\ {\bf BP_0} \subset \Hyper {\bf \Pi}(G_\nu)$ and
${\bf BP_0} \cup G_\nu \subset \Hyper {\bf \Pi}(G_\nu).$} \parm
For $\rm A \in \rm L_0,$ let $\#(\rm A)$ denote the length of the formula A.
\parm
{\bf Theorem 4.3.2} {\sl Let $\rm B \subset {\rm BP}_0,\ \rm A = \rm V \land 
\cdots \land \rm V \land \rm b,\ \rm b \in \rm E_0,$ where there are $n \geq 1$ 
connections $\land$ and for each $ \rm Z \in \rm B, \#(\rm A) > \#(\rm Z).$ Then 
$\rm A \notin \Pi (\rm B).$}\pars
Proof. (Note: In this proof certain of the indicated parentheses may be 
superfluous.) Assume the hypothesis of the theorem. We show that
there does not exist a finite $\rm F \subset \rm B$ such that $\rm F_( \vdash_\pi 
\rm A_(.$ Assume that there exists a finite $\rm F \subset \rm B$ such that $\rm F_( 
\vdash_\pi \rm A_(.$ A relation $R_m \subset {\rm BPC}_{0(} \times {\rm 
BPC}_{0(}$ is called an {\it m-chained sequence of {\rm MP} processes} if 
there exists an $ m\in \nat,\  m\geq 1$ such that $(X,Y) \in R_m$ has the 
form $(\rm A_i,(\rm A_i) \to (\rm A_{i-1}))$ for $ i= 1,\ldots,m,$ where $\rm A_0 
= \rm A_(.$ Also each $\rm A_i = \rm V \land (\rm A_{i-1})$ and $\rm A_i$ is a step 
in the proof of $\rm F_( \vdash_\pi \rm A_(,$ where $ i= 1,\ldots, m.$ We now show 
by induction that$,$ for each $m \in \nat,\ m \geq 1,$ there exists an 
m-chained sequence of MP processes in the proof that $\rm F_( \vdash_\pi \rm 
A_(.$\pars
{\it Case} m =1. We know that $\rm A_( \in {\rm BP}_{0(}$ implies that $\rm 
A_($ is not an instance of the use of an axiom since no axioms appear in ${\rm 
BP}_{0(}.$ Since $\rm A \notin \rm B$ then $\rm A_( \notin \rm F_(.$  Thus $\rm A_($ 
being the last step in the proof implies  that $\rm A_($ is the conclusion of 
an MP process with premises $(\rm D) \to (\rm A_()$ and D. Assume that D = c. 
Then the single step which contains the primitive c could not be an axiom nor 
an assumption since $\rm c \notin \rm B.$ Hence c would be the conclusion of a 
prior  MP process. Therefore a prior step would be of the form $(\rm E) \to \rm 
c.$ This is impossible; all steps must be formula in ${\rm  Ax} \cup {\rm 
BPC}_{0(}.$ Thus $\rm D\not= \rm c$ implies that $(\rm D)\to (\rm A_()$ is an 
instance of an axiom or the conclusion of a prior MP process. However$,$ since 
no step can be of the form $(\rm E) \to ((\rm D)\to (\rm A_()),$ it follows that
$(\rm D) \to (\rm A_() $ must be an instance of an axiom. Consequently$,$ $\rm D = \rm 
V \land (\rm A_() = \rm A_1 \in {\rm BPC}_{0(}$ and D is a step in the proof. 
Therefore$,$ $\rm R_1 = (\rm A_1,(\rm A_1) \to (\rm A_0)),\ \rm A_0 = \rm A_($ is a 1-
chained sequence of MP processes. \pars
Assume the result holds for m.\pars
{\it Case} m + 1. Let $\rm R_m = \{(\rm A_m,(\rm A_m) \to (\rm A_{m-1})), 
\ldots, 
(\rm A_1,((\rm A_1)\to (\rm A_0))\}$ be an m-chained sequence of MP 
processes. Now $\rm A_m = \rm V\land (\rm A_{m-1})$ implies by a simple induction
proof that $\#(\rm A_m) > \#(\rm A_0)= \#(\rm A_().$ Hence$,$ $\rm A_m \notin \rm B_(.$ 
Thus
$\rm A_m \notin \rm F_(.$ Moreover$,$ $\rm A_m \notin \rm C,\ \rm A_m \notin {\rm 
Ax}$ for the primary connective is $\land.$ This implies that $\rm A_m \in {\rm 
BP}_{0(}$ and $A_m$ must be the conclusion of some MP process with premises
$(\rm D)\to (\rm A_m)$ and D. As in case m =1$,$  it follows that $\rm D \not= c$ 
and that $(\rm D) \to (\rm A_m)$ must be an instance  of an axiom. Consequently$,$ 
$\rm D = \rm V \land (\rm A_m).$ Let $\rm D = \rm A_{m+1}.$ Since D is a step in 
the proof$,$ $\rm D \in {\rm BPC}_{0(}.$ Thus$,$ $\rm R_{m+1} = \{(\rm A_{m+1},
(\rm A_{m+1}) \to (\rm A_m))\} \cup \rm R_m$ is an $m+1$-chained sequence of MP 
processes.\pars
The length of the proof that $\rm F_( \vdash_\pi \rm A_($ is some  finite 
number$,$ say $n \in \nat,\ n\geq 1.$ The above shows that there exists an 
n+1-chained sequence of MP processes for this proof. Since 
$\vert P_1(\rm R_{n+1}) \vert = n+1$ (note that for each $i,\ j \in \nat$ such 
that $0\leq i < j \leq n+1,\ \#(\rm A_i) < \#(\rm A_j))$ and each element of $P_1(\rm 
R_{n+1})$ is a distinct step in the proof$,$ this contradicts the fact that 
the proof length is n. Consequently$,$ there does not exist a finite $\rm F \subset 
\rm B$ such that $\rm F_( \vdash_\pi \rm A_(.$ Therefore $\rm A \notin \Pi(\rm B)$ 
and this completes the proof of this theorem. \qed
Under our embedding$,$ Theorem 4.3.2 is interpreted by $\theta$ as embedded into 
$\cal M.$ When this is done the 
length of a formula A$,$ $\#(\rm A),$ is the length of the preimage A of the map i 
associated with the 
special partial sequence $f_{\rm A}(0) = i(\rm A).$  Let 
${\cal G} = \{G_n\mid (n > 0) \land (n \in \nat )\} \cup \{G_0\},\ G_0 = \bf 
E_0,$ and $\vert {i[\rm E_0]} \vert = m+1.$ Indeed$,$ $\rm E_0 = j[B]= \{\rm 
b_0,\ldots, \rm b_m \}.$ Let ${\cal G}\subset X_n.$ The following sentences hold in ${\cal M};$ hence in $
\hyper {\cal M}.$\parm

\line{\hfil $\forall x(x \in X_n \to ((x \in {\cal G}) \iff \exists y\exists y_0\cdots 
\exists y_m((y \in \nat)\land (y \geq 0)\land (y_0 \in T^y) $\hfil }\pars
\line{\hfil $\land \cdots \land (y_m \in T^y) \land \forall w((w \in \nat) 
\land (0<w\leq y) \to (y_0(w) =$\hfil }\pars
\line{\hfil $i({\rm V} \land)) \land \cdots \land (y_m(w) = i({\rm V} \land)))\land 
(y_0(0) = {\rm b_0})$\hfil }\pars 
\line{\hfil \qquad \qquad$ \land \cdots \land (y_m(0) = {\rm b_m} \land ([y_0] \in x) \land 
\cdots \land ([y_m]\in x)))).$ \hfil (4.3.4)} \parm 
\noindent Sentence (4.3.4) can also be written as \parm
\line{ \hfil $ \forall x(x \in X_n \to (x \in {\cal G} \iff \exists y((y \in \nat) \land 
(y\geq 0) \land A(y))));$\hfil }\pars
\line{\qquad\hfil $ \forall y((y \in \nat)\land (y\geq 0) \to A(y)),$  \hfil (4.3.5)
}\parm
\noindent where $A(y)$ is the obvious expression taken from (4.3.4). The objects that 
exist for each   ``y'' in the $A(y)$ expression (i.e. the $y_j \in T^y,\ j = 
0, \ldots, m$) are unique with respect to the property expressed in $A(y).$ 
Obviously$,$ for each $ n\in \nat,\ G_n \in \Hyper \cal G.$ Moreover$,$ there 
exists a bijection $F\colon \nat \to \cal G$ such that $F(n) = G_n.$ Now let $ n,\ 
m \in \nat,\ n , m$ and $\rm A =\rm V \land \ldots \land \rm b,\ \rm b \in \rm E_0$ 
has $m$ connectives $\land$. Then $\#( \rm A ) > \#( \rm Z )$ 
for each $[f_Z] \in G_n.$ It follows from Theorem 4.3.2 that $[f_A] \notin {\bf 
\Pi} 
(G_n).$ \parm
{\bf Theorem 4.3.3} {\sl If $ n,\ m \in \hypernat,\ 0\leq n< m, \ b \in  
i[E_0]$ and $ f \in (\Hyper T)^m$ such that $ f(w) =i(\rm V \land)$ for each 
$ w\in \hypernat, \ 0 < w \leq m,$ and $f(0) = b,$ then it follows that 
$[f]\notin \Hyper {\bf \Pi}(G_n).$}\pars
Proof. Let $ b \in i[\rm E_0].$ From the above discussion$,$ the following 
sentence\pars
\line{ \hfil $\forall x \forall y \forall z ((x\in \nat) \land (y\in \nat) 
\land (0\leq y < x) \land (z \in T^x)\land $\hfil}\pars
\line{ \hfil $\forall w ((w \in \nat) \land (0 < w \leq x) \to (z(w) = i(\rm V 
\land)) \land$ \hfil}\pars
\line{(4.3.6)\hfil $ (z(0) = b)) \to ([z] \notin {\bf \Pi}(F(y))))$\qquad\hfil 
}\parm
\noindent holds in $\cal M$; hence in $\hyper {\cal M}.$ The result follows by 
*-transfer. \qed
{\bf Corollary 4.3.3.1} {\sl For each $n \in \hypernat,\ \Hyper {\bf 
\Pi}(G_n) \not= \Hyper {\bf \Pi}(\Hyper {\bf BP_0})$ and 
$\Hyper {\bf \Pi}(G_n) = \bigcup \{G_x \mid( x\in \hypernat) \land (0\leq x\leq 
n)\}.$}\pars
Proof. Since $G_n \subset \Hyper {\bf BP_0},$ $\Hyper {\bf \Pi}(G_n) \subset
\Hyper {\bf \Pi}(\Hyper {\bf BP_0}).$ Theorems 4.3.1 and 4.3.3 along with the 
above discussion completely characterizes the elements of $\Hyper {\bf 
\Pi}(G_n).$ This completes the proof. \qed
For $\nu \in \nat_\infty,$ let ${\cal G}_\nu= \bigcup \{G_x \mid (x \in  
\nat_\infty) \land ( x < \nu)\}.$ Then ${\cal G}_\nu$
is a {\it purely subtle external object.} This follows from the fact that
$\Hyper {\bf \Pi}(G_\nu)$ is internal$,$ ${\bf BP}_0$ is external and 
$\Hyper {\bf \Pi}(G_\nu) ={\cal G}_\nu \cup {\bf BP}_0\cup G_\nu.$ Moreover$,$ 
observe that if $\lambda,\nu \in \nat_\infty, \nu > \lambda,$ then $\Hyper 
{\bf \Pi}(G_\lambda) \cap G_\nu = \emptyset$ and that ${\cal G}_\nu$ and ${\bf 
BP}_0$ are not in the domain of $\Hyper {\bf \Pi}$ unless we extend $\Hyper 
{\bf \Pi},$ say by the identity operator.\pars
{\bf Theorem 4.3.4} {\sl The set $\Hyper {\bf \Pi}(\Hyper {\bf BP}_0) = 
\Hyper {\bf BP}_0.$}\pars
Proof. (Note once again that some superfluous parentheses may have been added 
to some formula in this proof.) It is know that  $\Hyper {\bf BP}_0 \subset 
\Hyper {\bf \Pi}(\Hyper {\bf BP}_0).$ Let finite $\rm F \subset {\rm BP}_0,\ 
\rm A \in \rm C$ and assume that $\rm F_( \vdash_\pi \rm A.$ Then $ \rm A = \rm c$ 
or
$\rm A = \rm c\to (\rm x), \ \rm x \in {\rm BP}_{0(}.$ \pars
{\it Case} 1. Assume that $ \rm A= \rm c.$ Since $ \rm c\notin \rm F_( \subset
{\rm BP}_{0(}$ and $ \rm A = \rm c$ is not an instance of an axiom$,$ $\rm A = \rm c$ 
must be the conclusion of an MP process. Thus a prior step is of the form 
$(\rm D) \to \rm c.$ This is impossible for $(\rm D) \to \rm c \notin {\rm Ax} \cup 
{\rm BP}_{0(}.$ \pars
{\it Case} 2. Assume that $\rm A = \rm c \to (\rm x),\ \rm x \in {\bf BP}_{0(}.$ 
Again $ \rm c \to (\rm x)$ is the 
conclusion of an MP process. This is impossible since no formula of the type
$(\rm D) \to (\rm c\to (\rm x))$ is an element of ${\rm Ax} \cup {\rm BP}_{0(}.$ 
Hence by *-transfer of the appropriate first-order sentence$,$ after the 
$\theta$ embedding$,$ it follows that $\Hyper {\bf \Pi}(\Hyper {\bf BP}_0) 
\subset \Hyper {\bf BP}_0.$ \qed
{\bf Corollary 4.3.4.1} {\sl The set $\Hyper {\bf \Pi}(\Hyper {\bf BP}_0)
{\subset\atop\not=} \Hyper {\bf \Pi}(\Hyper {\bf BPC}_0).$}\parm
It is easy to see that $\Pi (\rm C) = {\rm BPC}_0.$ For let $\rm A \in  {\rm 
BP}_0$ and consider the proof (1) c$,$ (2) $\rm c \to (\rm A_(),$ (3) $\rm A_(.$
Thus $\{\rm c, \rm c\to (\rm A_() \} \vdash_\pi \rm A_($ yields that $\rm A \in 
\Pi (\rm C).$ Hence $\Pi (\rm C) = {\rm BPC}_0.$ Also$,$ $\Hyper {\bf \Pi} 
(\hyper {\bf C}) = \Hyper {\bf BPC}_0.$\parm
{\bf Theorem 4.3.5} {\sl Let internal $A \subset \Hyper {\bf BP}_0$ and 
internal $B\subset \Hyper {\bf BPC}_0.$ Then $\Hyper {\bf \Pi}(A \cup B) 
= \Hyper {\bf BPC}_0$  
iff $\hyper {\bf C} \subset B.$}\pars
Proof. For the sufficiency$,$ let internal $A \subset \Hyper {\bf BP}_0$$,$ 
internal $\hyper {\bf C} \subset B.$ Then $A \cup B$ is internal and 
$\Hyper {\bf \Pi}(\hyper {\bf C}) = \Hyper {\bf BPC}_0 \subset \Hyper {\bf 
\Pi}(B) \subset \Hyper {\bf \Pi}(A \cup B) \subset \Hyper {\bf BPC}_0.$ Thus
$\Hyper {\bf BPC}_0 = \Hyper {\bf \Pi}(A \cup B).$\pars
For the necessity$,$ assume that internal $A \subset \Hyper {\bf BP}_0,$ 
internal
$B \subset \Hyper {\bf BPC}_0$ and that $\Hyper {\bf \Pi}(A \cup B) = 
\Hyper {\bf BPC}_0.$  Let $\rm A_1 \subset {\rm BP}_0,\ \rm B_1 \subset {\rm 
BPC}_0$ and $\Pi(\rm A_1 \cup \rm B_1) = {\rm BPC}_0.$ It follows from Theorem 4.3.4
that $\rm B_1 \not\subset {\rm BP}_0.$ Indeed$,$ given any finite $\rm F \subset \rm 
A_1 \cup(\rm B_1 \cap {\rm BP}_0).$ If $ \rm D \in \Pi(\rm F),$ then $\rm D \in 
{\rm BP}_0.$ Thus only for a finite $\rm F_1 \subset \rm B_1 \cap \rm C$ can 
there be an $\rm E \in \rm C$ such that $ \rm E \in \Pi(\rm F_1).$ Hence all that 
needs to be shown is that $\rm C\subset \Pi (\rm B_1 \cap \rm C)$ implies     
that $\rm B_1 \cap \rm C = \rm C.$ So$,$ assume that $\rm B_1 \cap \rm C \not= \rm C.$ 
Hence either $\rm c\notin \rm B_1 \cap \rm C$ or there exists some $\rm A_(  \in 
{\rm BP}_{0(}$ such that $\rm c \to (\rm A_() \notin (\rm B_1 \cap \rm C).$ \pars
{\it Case} 1. Assume that $\rm c\notin \rm B_1 \cap \rm C$ and $\rm F$ is any 
finite
subset of $\rm B_1 \cap \rm C$ such that $\rm F \vdash_\pi \rm c.$ Of course$,$ c is 
the last step in a formal proof. c is the conclusion of some MP process since 
c is not an assumption nor an axiom. Thus some formula of the form $(\rm D) \to
\rm c$ must be in a prior step in the formal proof. This is impossible since no formula of this form 
is an element of ${\rm BPC}_{0(}.$\pars
{\it Case} 2. Assume that there exists some $\rm A_( \in {\rm BP}_{0(}$ such 
that $\rm c \to (\rm A_() \notin \rm B_1 \cap \rm C$ and there exists finite $\rm F  
\subset \rm B_1 \cap \rm C$ such that $\rm F \vdash_\pi \rm c\to (\rm A_().$ Again 
$\rm c \to (\rm A_() $ is not 
an assumption nor an axiom. Consequently$,$ $\rm c \to (\rm A_()$ is the conclusion 
of an MP process. Thus there exists some formula of the form $(\rm D) \to (\rm c 
\to (\rm A_())$ in a prior step. Again this is impossible.\pars
These two cases imply that $\rm B_1 \cap \rm C =\rm C.$ Therefore$,$ $\rm C 
\subset \rm B_1$ implies the sentence\parm
\line{\hfil $ \forall x\forall y((x\in \power {{\bf BP_0}}) \land (y \in \power 
{{\bf BPC_0}}) \land$\hfil}\pars
\line{\qquad\qquad\hfil$ ({\bf \Pi}(x \cup y) = {\bf BPC_0}) \to ({\bf C} \subset 
y))$ \hfil (4.3.7)}\parm
\noindent holds in $\cal M;$ hence in $\Hyper {\cal M}.$ The result follows 
from *-transfer. \qed
Note that all of the results in this section hold for ${\rm BP}$ and 
${\rm BPC},$ where $\rm C$ is constructed without parentheses. \par
\bigskip
\leftline{\bf 4.4 Order}
\medskip
We briefly look at two special types of order relations$,$ the  ``number of 
symbols'' order and the  ``better than'' order. Previously the concept of the 
length of a formula or word A (i.e. \#(A)) was introduced. This type of order 
has  few properties unless it is restricted to certain interesting types of 
subsets.\pars
Let nonempty $\rm B, \ \rm D \subset {\rm BPC}$ (or ${\rm BPC}_0$)$,$ then define 
$\rm B\leq_\# \rm D$ 
if for each $\rm b  \in \rm B$ and for each $\rm d \in  \rm D$$,$ it follows that
$\#(\rm b) \leq \#(\rm d).$ This order is obviously a pre-order in the sense 
that it is reflexive and transitive. However$,$ in general$,$ it should probably 
not be considered a  partial 
order since antisymmetry does not imply set equality although it does imply 
that all the symbol strings have equal length in both B and D. Also other 
pre-orders of this type appear not to be partial orders for the same reason. 
If $\leq_\#$ is restricted to certain collections of sets$,$ then it does 
become a useful partial order under set equality.\pars 
Consider the collection $\{G_n \mid n \in \hypernat \}.$ Then  the 
pre-order $\leq_\#$ restricted to this set is isomorphic to the simple order of 
$\hypernat.$ Indeed$,$ $G_n \leq_\# G_m$  iff $ n\leq m,$ where $n ,\ m\in 
\hypernat $ and $\leq $ is the usual extension of the simple order induced on 
$\nat$ by $\omega.$ Moreover$,$ notice that $G_n = G_m$ iff $n = m,$ and
$G_n \not=_\# G_m$ iff $G_n \cap G_m = \emptyset.$ \pars
For the collection $\{{\cal G}_\nu \mid \nu \in \nat_\infty \},$  it follows 
that ${\cal G}_\zeta \subset {\cal G}_\lambda$ iff $\zeta \leq \lambda.$ 
Thus $\{{\cal G}_\nu \mid \nu \in \nat_\infty \}$ is ordered by inclusion when 
the simple order of the subscripts is considered. Notice also that $\bigcap \{ 
{\cal G}_\nu \mid \nu \in \nat_\infty \} = \emptyset.$ \pars
Let fixed $\nu \in \nat_\infty.$ Then there exist infinitely many ${\cal 
G}_\lambda$ which differ only be a finite set of subtle objects. Simply 
consider the set $\{ {\cal G}_{\nu + n}\mid n \in \nat \}.$ If $ n
,\ m \in \nat $ and $ m > n,$ then $\vert {\cal G}_{\nu + m} - {\cal G}_{\nu + 
n} \vert = (m - n) \vert G_0 \vert\in \nat$. Also there 
exist 
infinitely many sets  ``longer than'' any $G_{\nu +n},$ where $ n \in \nat$ or 
strictly containing any ${\cal G}_{\nu +n}.$ To see this consider
$\nu^2 < \nu^3 < \cdots < \nu^n < \cdots,\ n \in \nat,$ and observe that 
$\nu^2 - \nu = \nu(\nu - 1) \in \nat_\infty.$ Thus the length of an interval 
$[\nu^m,\nu^n],\ m < n,$ is an infinite natural number and for any $n \in 
\nat,$ this implies that $\nu + n < \nu^2.$ Hence for any $ n \in \nat$ such 
that $ n > 1,$ it follows that ${\cal G}_{\nu +n} \subset {\cal G}_{\nu^n}$ 
and $G_{\nu + n} \leq_\# G_{\nu^n}.$ It is also interesting to note that 
for each pair $ \nu,\ \lambda$ of infinite natural numbers$,$ such that $\nu 
\leq \lambda,$ $\Hyper {\bf \Pi}(G_\nu) \subset \Hyper {\bf \Pi}(G_\lambda)$ 
and conversely.  
\pars

The ``better than'' order is only defined for comparable readable
sentences. For this research$,$ the domain of definition is restricted to the 
set $\bf BP_0$ [resp. {\bf BP}]. Two elements $[f],\ [g] \in {\bf BP_0}$ 
[resp. {\bf BP}] are {{\it comparable}} if there exists $ b \in i[\rm E_0]$ 
[resp. $i[\rm B]$] such that $f^b \in [f]$ and $g^b \in [g].$ Recall that
$f^b$ and $g^b$ are unique element of $T^n$ and $T^m, $ respectively$,$ where   
n and m count the number  ``$\rm V \land$'' [resp.  ``${\rm very,}\sp$''] 
symbol strings. The $f^b,\ g^b$ are restricted to $T^n,\ T^m,$ where $\nu =n,m >0.$ For example, $0< x \leq n,\ f^b(x) = i({\rm very,}\sp)$ and $f^b(0) = b \in i[{\rm B}].$ For two comparable objects $[f],\ [g]$ define $[f] \leq_B [g]$
if $ n\leq m.$ Two nonempty sets $A,\ D \subset {\bf BP_0}$ [resp. {\bf BP}] 
have the property that $A \leq_B D$ if for each $[f] \in A$ there exists some
$[g] \in D$ such that $[f] \leq_B [g].$ 
This is the {\it better than} pre-order and usually
$[f] \leq_B [g]$ is stated as follows: ``$[g]$ is better than $[f]$'' or some 
similar  expression. Actually, for the ${\bf BP_0}$ [resp. {\bf BP}], the ``better than 
order'' is a partial order and$,$ in some cases$,$ it is equivalent to the 
$\leq_\#$ order. Of course$,$ $\leq_B$ and $\leq_\#$ are *-transferred to 
$\Hyper {\cal M}.$\pars
For each $\bf b \in B$, let $C_b=\{x\mid (x \in {\bf BP})\land({\bf b} \leq_B x\}.$\parm
{\bf Theorem 4.4.1} {\sl There exists a purely subtle $c \in \Hyper C_b$ such that $C_b \Hyper {\leq_B}\{c\}$}.\pars
Proof. The sentence\hfil\break
\line{\hfil$\forall x(x \in \nat \to \exists y(y\in T^x\land[y] \in C_b\land\forall z\forall w(z\in \nat\land z\leq x\land w \in T^z \land$\hfil}\parm
\line{(4.4.1)\hfil $[w] \in C_b \to  [w] \leq_B [y]))$\hfil}\medskip
\noindent holds in $\cal M$; hence, in $\hyper {\cal M}.$ \pars
Let $\nu \in \hypernat - \nat,$ then there is 
a $f \in (\Hyper T)^\nu$ and a purely subtle $c =[f] \in \Hyper C_b$, where $[f]$ satisfies the remainder of the *-transformed (4.4.1) statement. Let $a \in C_b \subset \Hyper {C_b}.$ Then there is some $m \in \nat$ and some $g \in T^m$ such that $a = [g]$. Thus, since $m \leq \nu,$ then $[g] \Hyper {\leq_B [f]}.$ Consequenely, $C_b \Hyper {\leq_B} \{c\}.$\qed

The following is somewhat trivial and is not formalized as a theorem. Consider the usual representation for $c = [f],\ f \in (\Hyper T)^\nu,\ \nu \in \hypernat - \nat.$  Intuitively, members of $\Hyper C_b(\{c\}),$ are obtained by removing *-finitely many (including 0) $i({\rm very,|||})$ from $c.$ Let $n \in \nat.$ Then $\{x\mid(x \in \hypernat)\land (n\leq x \leq \nu)\}$ is *-finite. By *-transfer of the appropriate sentence, you have the following for each $m \in \nat$. If $m = 0$, then $[g] \in \Hyper C_a(\{c\}),$ where $g(0) = \b b.$ If $m \geq 1,$ then $[g] \in \Hyper C_a(\{c\}),$ where $g(0) = \b b$ and, for each $j \in \nat,$ such that $1 \leq j \leq m,$ $g(j) = i({\rm very,|||}).$ Thus, for $\Hyper C_a(\{c\}),$ which is simply a restriction of *-propositional deduction, one has that attribute b as well as all of the $\rm very,|||\cdots very,|||b$ attributes are rationally related to $c.$ When Theorem 4.4.1 is interpreted, then $c$ is stronger than, better than, greater than, b or any of these standard strengthens of the basic b. \pars
  
\centerline{\bf 5. CONSEQUENCE OPERATORS }\par
\bigskip
\leftline{\bf 5.1 Basic Definitions}\par 
\medskip
Recall once again the 
Tarski [21] cardinality independent axioms for a {{\it finite consequence 
operator}} ${\rm C}\colon \power {\rm A} \to \power {\rm A}$   
on a nonempty set of meaningful sentences $\rm A.$   \pars %===SIMPLE TABLE OF 
%ONE BOX WITH MULTIPLE LINES %IT IS IN DISPLAYED FORM WHEN PRINTED 
$$\vbox{\offinterlineskip \hrule \halign{&\vrule#& \strut\quad#\hfil\quad\cr 
height2pt&\omit&\cr &${\rm (2)\ If}\ \rm B \subset \rm A,\ {\rm then}\ 
\rm B \subset {\rm C}(\rm B) \subset \rm A,$&\cr &\ &\cr &${\rm (3)\ If}\ \rm B 
\subset \rm A,\ {\rm then}\ {\rm C(C(B))=C(B)}, $&\cr &\ &\cr &${\rm (4)\ If  
}\ \rm B \subset \rm A,\ {\rm then}\ {\rm C(B) = \bigcup \{C(F)\bigm\vert 
F} \in {\it F}(\rm B)\}$&\cr       
height2pt&\omit&\cr} \hrule}$$                            
\pars
The modern theory of {{\it consequence operators (the term finite dropped)}}
alters axiom (4) and replaces it by axiom \pars
%===SIMPLE TABLE OF 
%ONE BOX WITH MULTIPLE LINES %IT IS IN DISPLAYED FORM WHEN PRINTED 
$$\vbox{\offinterlineskip \hrule \halign{&\vrule#& \strut\quad#\hfil\quad\cr 
height2pt&\omit&\cr &${\rm (5).\ If}\ \rm B,\ \rm D  \subset \rm A,\  
{\rm if}\ \rm B \subset \rm D \ {\rm then}\  {\rm C(B) \subset C(D)},$&\cr       
height2pt&\omit&\cr} \hrule}$$\pars  
   
It is very important to state that {\bf axioms (2)$,$ (3) and (4) imply axiom 
(5).} Thus a finite consequence operator is a consequence operator$,$ but not 
conversely. All things that can be established for consequence operators 
without any further axioms hold for finite consequence operators. For this 
reasons$,$ some of the following results will be established for consequence 
operators in general. Of course$,$ consequence operators need not be restricted 
only to objects that are considered to be language. In the theory here being developed, $\r A$ can be of two types. Either $\r A \subset {\cal W}$ or $\b A \subset {\cal E}.$  I shall$,$ however$,$ continue
to use Roman notation for all of the objects related to ${\cal W}$ so as to differentiate them from 
the other mathematical entities. In all that follows$,$ the symbol 
${\cal C^\prime}$ will denote the set of all consequence operators defined on some 
specified $\power {\rm A}$
and the symbol ${\cal C}_f^\prime$ the set of a finite consequence operators. 
Obviously$,$ ${\cal C}_f^\prime \subset 
{\cal C}^\prime,$ where$,$ if no mention is made of any other possible$,$ it will also be 
assumed that each member of ${\cal C}^\prime$ is defined on the same $\power {\rm 
A}.$\pars
On the sets ${\cal C}^\prime$ and ${\cal C}_f^\prime$$,$ we can define a significant partial 
order. For $\rm C_1,\ C_2 \in {\cal C}^\prime$ let $\rm C_1 \leq C_2$ if $\rm C_1(X) \subset C_2 
(\rm X)$ for each $X \in  \power {\rm A}.$ This partial order$,$ I term 
{{\it the stronger than order.}} The partial order defined on ${\cal 
C}_f^\prime$ is the restricted stronger than order.\par
A great deal has been discovered about  algebras $\langle  
{ \cal C}^\prime, \leq \rangle$ and $\langle { \cal C}_f^\prime, \leq \rangle.$ 
For example$,$ one can define a compatible meet operation as follows: 
For each $\rm C_1,\ C_2 \in {\cal 
C}^\prime$ [resp. ${\cal C}_f^\prime$]$,$ let the map $\rm C_1 \wedge C_2 \colon \power 
{\rm A} \to \power {\rm A}$ be  defined by $\rm (C_1 \wedge C_2)(X) = C_1(X) \cap 
C_2(X),$  where $\rm X \subset \rm A.$  
Each of these algebras has the same upper unit 
and the same lower unit. The lower unit is but the identity map on $\power
{\rm A}.$ The upper unit $\rm U$ is the map defined by $\rm U(X) = A$ for each 
$\rm X\subset A.$ These algebras are both meet semi-lattices.\pars
Our interest in the above two algebras is not in any deep investigation into\break
\vfil
\centerline{}
\eject
\noindent
there different properties but$,$ rather$,$ will be restricted two chains. In 
general$,$ $\langle {\cal C}^\prime, \leq \rangle$ is not closed under 
composition. [24] However$,$ for chains there is a very simple 
relation between the stronger than order and composition.\parm

{\bf Theorem 5.1.1} {\sl Let ${\cal D} \subset {\cal C}^\prime.$ Then 
${\cal D}$ is a chain in $\langle {\cal C}^\prime, \leq \rangle$ iff for each 
$\rm C_1, C_2 \in {\cal D}$ either the composition $\rm C_1C_2 = C _1$ or $\rm C_2C_1 = 
C_2.$}\pars
Proof. For the necessity$,$ assume that hypothesis. Suppose that 
$\rm C_1 \leq C_2.$ Then for each $\rm X \subset A,\ X \subset C_1(
X)\subset C_2(X)$ implies that $\rm C_2(\rm X) \subset C_2(C_1(\rm X))\subset C_2(C_2(\rm 
X)) =C_2(\rm X).$ Hence $\rm C_2C_1 = C_2.$ In like manner$,$ if $\rm C_2 \leq C_1.$ \pars
For the sufficiency$,$ let $\rm C_2C_1 = C_2.$ Then for each $\rm X \subset \rm A,\ 
C_1(\rm X) \subset C_2(C_1(\rm X)) = (C_2C_1)(\rm X) = C_2(\rm X)$ implies that
$\rm C_1 \leq C_2.$ In like manner for $\rm C_1C_2=C_1$ and this  completes the proof. 
\qed
\bigskip 
\leftline{\bf 5.2 Basic $\sigma$ Properties} 
\medskip Since consequence operators are relations between sets$,$ it becomes 
more essential to incorporate$,$ to a certain degree$,$ {the $\sigma$} operator 
into much of our discussion. Since we wish to maintain symbolic consistency and 
avoid trivialities$,$   
 assume that nonfinite $A\subset \cal W.$ One important
 result that will 
be used many 
times without further elaboration$,$ uses the finitary construction of the 
equivalence class $[f] \in \cal E.$ From our previous discussion$,$ a readable 
sentence $[f]$ behaves as follows: $\Hyper {[f]} = [\hyper f] = [f]$ under the 
identification of the natural numbers. Now if $ \emptyset \not= A \subset \cal 
E,$ then $^\sigma A =\{\Hyper {[f]}\mid [f] \in A \} = \{[f] \mid [f] \in A \} 
= A.$ The following result brings 
together various facts relative to the $\sigma$ operator all of which follow 
easily from the definitions and characterizing properties. The proofs will be 
omitted.\parm
{\bf Theorem 5.2.1} \pars

(i) {\sl Let $A \in \cal N.$ Then 
$^\sigma (F(A)) = F(^\sigma A)$. If also $A \subset ({\cal W} \cup {\cal E}),$ then $^\sigma 
(F(A)) = F(A).$}\pars

(ii) {\sl Let ${\rm C} \in {\cal C}^\prime,\ \r B \subset \r X 
\subset \cal W.$} \pars 

(a) {\sl $^\sigma (\b C(\b B)) = \b C(\b B).$}\pars 

(b) {\sl $\Hyper {\b C}\bigm| \{\hyper {\b A}\mid {\b A} \in \power X \}=
 \{ (\hyper {\b A},\hyper {\b B}) 
\mid (\b A,\b B) \in \b C \} = {^\sigma \b C}.$} \pars 

(c) {\sl If $\r F \in \r F(\r B),$ 
then $^\sigma (\b C(\b F)) \subset (^\sigma \b C )(^\sigma {\b F}) = 
(^\sigma \b C)(\b F).$ Also 
$^\sigma (\b C(\b B)) \subset {(^\sigma \b C) (\hyper {\b B})}$ and$,$ in general$,$ 
$^\sigma (\b C(\b B)) \not= (^\sigma \b C)(\hyper {\b B}),\ ^\sigma (\b C(\b F)) \not= 
(^\sigma \b C)(F).$} \pars 

(d) {\sl If $C \in {\cal C}_f^\prime,$ then 
$^\sigma (\b C(\b B))= \bigcup \{^\sigma (\b C(\b F)) \mid \b F \in {^\sigma (F(\b B))}\} 
= \bigcup \{ ^\sigma (\b C(\b F)) \mid \b F \in F(^\sigma {\b B})\} = \bigcup \{\b C(\b F) 
\mid \b F \in F(\b B) \}.$}\parm

(A duplicate theorem holds, where $A \in \cal N$ and $C \in \cal C$, where $\cal C$ set of consequence operators defined on subsets of W if W is included as a subset of the ground set. The difference is that the ``bold'' notion does not appear.)\parm

Throughout the reminder of this section$,$ in order to escape trivialities$,$ we 
{remove the upper unit $\rm U$ from the collections of consequence operators.}
Let ${\cal C} = {\cal C}^\prime - \{\rm U\}$ and ${\cal C}_f =
 {\cal C}_f^\prime -  \{\rm U\}.$ One of the consequences 
of this last 
requirement shows that if $D\in \cal N$ satisfies the axioms for a 
consequence operator and $G \subset \cal C,$ then there 
does not exist a consequence operator $C \in \Hyper {\b G}$ such that 
${^\sigma \b D} = C.$ To see this simply note that from Theorem 5.2.1 
${^\sigma \b D}$ is defined on extended standard sets while each member of 
$\Hyper {\b G}$ is defined on the internal subsets of $A.$ Since $A$ is not 
finite$,$ there exists internal subsets of $A$ that are not equal to any 
extended standard set. {\bf There is also one useful general fact.}
 Consider {\it 
any} sets $A, B \in \cal N$ such that $A \subset B.$ Then $\hyper A - {^\sigma 
A} \subset \hyper B - {^\sigma B}.$ For suppose that there exists some
$X \in  (\hyper A - {^\sigma A})$ such that $X \in {^\sigma B}.$ Then 
$X = \hyper D$ for some $D \in B.$ Thus $\hyper D \in \hyper A$ implies that
$D \in A.$ From this$,$ we have the contradiction that $X = \hyper D \in 
{^\sigma A}.$ Also note that there does not exist $D \in \cal N$ such that
$\hyper D\in (\hyper B - {^\sigma B})$ (i.e. each member of $(\hyper B - 
{^\sigma B})$ is an {internal pure nonstandard object} or a 
{pure subtle object}.)
\par
\bigskip
\leftline{\bf 5.3 Major Results}
\medskip For the algebras $\langle {\cal C}, \leq \rangle$ and
$\langle {\cal C}_f, \leq \rangle$ {two types of chains} will be studied.
Denote by  $\rm K$ any nonempty chain contained in either of 
these algebras and by  $\rm K_\infty$ a chain 
with the following property. For each $\rm C \in K_\infty$ 
there exists $\rm C^\prime\in K_\infty $ such that $\rm C < C^\prime.$ \parm
{\bf Theorem 5.3.1} {\sl There exists $C_0 \in \Hyper {\b K}$ such that for each 
${\rm C \in K},\ \Hyper {\b C} \leq C_0.$ There exists some $C_\infty \in \Hyper {\b 
K_\infty}$ such that $C_\infty$ is a purely subtle consequence operator and for 
each $C \in K_\infty,\ \Hyper {\b C} < C_\infty.$ Each member of $\Hyper {\b 
K}$ and $\Hyper {\b K_\infty}$ are subtle consequence operators.}\pars 
Proof. Let $R = \{(x,y)\mid (x \in \b K) \land (y \in \b K)\land ( x\leq y) \} $ 
and  
$R_\infty =\{ (x,y) \mid (x \in \b K_\infty) \land (y \in \b K_\infty) \land (x 
< y) \}.$ In the usual manner$,$ it follows that $R$ is concurrent on $\b K$ and 
$R_\infty$ is concurrent on $\b K_\infty.$ Consequently$,$ there is some $C_0 
\in \Hyper {\b K}$ and some $C_\infty \in \Hyper {\b K_\infty}$ such that for 
each $\rm C \in  K$ and each $\rm C^\prime \in  K,\ \Hyper { \b C} \leq C_0,\ 
\Hyper {\b C^\prime} < C_\infty$ since $\Hyper {\cal M}$ is an enlargement. 
Further$,$ it follows that $C_\infty \in \Hyper {\b K_\infty} - {^\sigma \b 
K_\infty}$ implies that $C_\infty$ is a purely subtle consequence operator. 
Note that each member of $\Hyper {\b K} \cup \Hyper {\b K_\infty}$ is defined 
on the set of all internal subsets of $\Hyper {\b A}.$ This completes our proof. \qed 
 Notice that  ${\rm C_\infty}$ is {stronger than} or  
``more powerful than'' any $\rm C \in K_\infty$ in the following sense. If $\rm B 
\subset \rm A,$ then for each $\rm C \in K_\infty,$ it follows that 
$\b C(\b B) \subset \Hyper {(\b C(  \b B))} = \Hyper {\b C}( \Hyper {\b B}) 
\subset C_\infty (\Hyper {\b B}).$ Also for each $\rm C \in K_\infty$ 
there exists some internal $E_C  \subset \hyper {\b A}$ and 
$\Hyper {\b C}(E_C){\subset \atop \not=} C_\infty (E_C).$ 
Recall that for $\rm C \in \cal C,$ a 
set $\rm B \subset \rm A$ is a {{\it C-deductive system}} if $\rm C(B) = B.$
Also$,$ when we write the *-operator on any map $f$ in the form $\hyper {f(x)}$ 
this always means $(\hyper f)(x)$ rather than $\Hyper {(f(x))}.$ \parm
{\bf Theorem 5.3.2 } {\sl Let ${\rm C} \in {\cal C}_f$ and 
$\rm B \subset \rm A\subset \cal W.$ Then there exists a *-finite $F \in 
\Hyper (F(\b B))=\hyper F(\Hyper {\b B})$ 
such that $\b C(\b B) \subset \Hyper {\b C} ( F) \subset \Hyper {\b C}(\Hyper 
{\b B}) = \Hyper {(\b C(\b B))}$ and $\Hyper {\b C}(F) \cap \b A = \b C(\b B) = 
\Hyper {\b C}(F) \cap \b C(\b B)$}.\pars
Proof. Consider the binary relation $Q = \{(x,y) \mid (x \in \b C(\b B)) 
\land (y \in F(\b B)) \land (x \in \b C(y))\}.$ By axiom (4)$,$ the domain of $Q$ 
is $\b C(\b B).$ Let $(x_1,y_1), \ldots, (x_n,y_n) \in Q.$ By Theorem 1 in 
[5$,$ p. 64] (i.e. axiom (5)) we have that $\b C(y_1) \cup \cdots \cup
\b C(y_n) \subset \b C(y_1 \cup \cdots \cup y_n).$ Since $F = y_1 \cup \cdots 
\cup y_n \in F(\b B),$ then $(x_1,F),\ldots, (x_n,F) \in Q.$ Thus $Q$ is 
concurrent on $\b C(\b B).$ Hence there is some $F \in \Hyper {(F(\b B))}$ such 
that $^\sigma(\b C(\b B)) = \b C( \b B) \subset \Hyper {\b C}(F)\subset
\Hyper {\b C}(\Hyper {\b B}) = \Hyper {(\b C(\b B))}.$ Since $^\sigma A = A, 
\Hyper {\b C}(F) \cap \b A = \b C( \b B)= \Hyper {\b C}(F) \cap \b C( \b B).$ 
\qed
{\bf Corollary  5.3.2.1 } { \sl If ${\rm C} \in {\cal C}_f$ and $\rm B \subset \rm A \subset 
\cal W$ is a C-deductive system$,$ then there exists a *-finite $F \subset 
\Hyper {\b B}$ such that $\Hyper {\b C}(F) \cap \b A = \b B.$}\parm
{\bf Corollary 5.3.2.2} {\sl Let $ {\rm C} \in {\cal C}_f.$ Then there exists a *-
finite $F \subset \Hyper {\b A}$ such that for each 
$\rm B \subset \rm A,\  \Hyper {\b C}(F) \cap \b B = \b B.$}\pars
Proof. In Theorem 5.3.2$,$ let the  ``B'' be equal to A. Then there exists some 
*-finite $F \subset \Hyper {\b A}$ such that $\Hyper {\b C}(F) \cap \b A = 
\b C(\b A) = \b A.$ Thus $\Hyper {\b C}(F) \cap \b A \cap \b B = \Hyper {\b 
C}(F) \cap \b B =\b A \cap \b B = \b B.$ \qed
{\bf Theorem  5.3.3} { \sl Let $\rm B \subset \rm A \subset \cal W.$\pars
{\rm (i)} There exists a *-finite $F_B \in \Hyper {(F(\b B))}$ and a subtle 
consequence operator $C_B \in \Hyper {\b K}$ such that for all $C \in K,\ 
^\sigma(\b C(\b B)) = \b C(\b B) \subset C_B(F_B).$\pars
{\rm (ii)} There exists a purely subtle consequence operator $C_B^\infty \in 
\Hyper {\b K_\infty}$ such that for all $\rm C \in K_\infty,\ 
^\sigma (\b C(\b B)) = \b C(\b B) \subset C_B^\infty (F_B).$}\pars
Proof. (i) Consider the binary relation $Q = \{((x,z),(y,w))\mid (x \in \b K) 
\land ( y \in \b K) \land (w \in F(\b B)) \land (z \in x(w))\land (x(w) \subset
y(w))\}.$ Let nonempty $\{ ((x_1,z_1),(y_1,w_1)), \ldots, ((x_n,z_n),(y_n,w_n))\} 
\subset Q.$ Notice that $F = w_1 \cup \cdots \cup w_n \in F(\b B)$ and 
the set $R = \{x_1,\ldots, x_n\}$ has a largest member $D$ with respect to the
$\cal E$ embedded $\leq$ ordering for the consequence operators. It follows 
that $z_i \in x_i(w_i) \subset x_i(F) \subset D(F)$ for each $i = 1,\ldots, 
n.$ Hence $\{((x_1,z_1),(D,F)),\ldots,((x_n,z_n),(D,F))\} \subset Q$ implies 
that $Q$ is concurrent on its domain. Consequently$,$ there exists some 
$(C_B,F_B) \in \Hyper {\b K} \times \Hyper {(F(\b B))}$ such that for each 
$(x,z) \in $ domain of $Q,\ (\Hyper {(x,z)},(C_B,F_B)) \in \Hyper Q.$ 
Therefore$,$ each $(u,v) \in $ $^\sigma$(domain of $Q),\ ((u,v),(C_B,F_B))\in 
\Hyper Q.$ Let arbitrary $C \in K$ and $\rm b \in \rm C(\rm B).$ Then there exists 
some $F^\prime \in F(\rm B)$ such that $\rm b \in C(F^\prime).$ Thus $(\Hyper {\b 
C}, \Hyper {\b b}) \in $ $^\sigma$(domain of $Q).$ Consequently$,$ for each $C \in K$ 
and $\rm b \in C( \rm B),\ \b b = \Hyper {\b b} \in \Hyper {\b C}(F_B) \subset 
C_B(F_B).$  This all implies that for each $C \in K,\ ^\sigma (\b C(\b B)) = 
\b C(\b B) \subset C_B(F_B).$ \pars
(ii) Change the relation $Q$ to $Q^\prime$ be adding the additional 
requirement to $Q$ that $x \not= y.$ Replace the $D$ in (i) with any 
$D^\prime$ that is greater than and {\bf not equal} to the largest member of $R.$
Such a $D^\prime$ exists in $\b K_\infty$ from the definition of $\rm K_\infty.$
Continue the proof in the same manner as in (i) to obtain $C_B^\infty$ and 
$F_B.$ The fact that $C_B^\infty$ is a purely subtle consequence operator 
follows as in the proof of Theorem 5.3. \qed
{\bf Corollary 5.3.3.1} {\sl There exists a {\rm [}resp. purely{\rm ]} subtle 
consequence operator $C_A \in \Hyper {\b K}$ {\rm [}resp. $ \Hyper{\b 
K_\infty}$ {\rm ]} and a *-finite $F_A \in \Hyper {(F(\b A))}$ such that for all 
$\rm C \in K$ {\rm [}resp. $\rm K_\infty$ {\rm ]} and each $\rm B \subset \rm A,\ \b B 
\subset \b C(\b B) \subset C_A(F_A).$}\pars
Proof. Simply let the ``B'' in Theorem 5.3.3 be equal to A. Then there exists 
a [resp. purely] subtle $C_A \in \Hyper {\b K}$ [resp. $\Hyper {\b K_\infty}$] 
such that for all $\rm C \in K$ [resp. $\rm K_\infty$]$,$ $\b C(\b A) \subset 
C_A(F_A).$ If $\rm B \subset \rm A$ and $C \in K$ [resp. $K_\infty$]$,$ then 
$\rm B \subset C(\rm B) \subset C(\rm A).$ Thus for each $\rm B \subset \rm A$ and 
$\rm C \in K$ [resp. $\rm K_\infty$]$,$ $\b B \subset \b C(\b B) \subset C_A(F_A)$ and 
this completes the proof. \qed
Relative to the above results$,$ it is well known that for $\b B \subset \b A$ 
that there exists a *-finite $F_1\subset \Hyper {\b B}$ such that
$\b B \subset F_1 \subset \Hyper {\b B}.$ Thus for any $\rm C\in \cal C$ it 
follows that $\b B \subset \b C(\b B) \subset \Hyper {\b C(F_1)} \subset 
\Hyper {\b C}(\Hyper {\b B}).$ One significance of the above results is that 
the $C_B^\infty$ is purely subtle and$,$ thus$,$ not the same as any extended 
standard consequence operator.\par
\bigskip
\leftline{\bf 5.4 Applications}
\medskip
In what follows$,$ let denumerable $\rm L$ be a language constructed from a 
denumerable set of primitive symbols $\{\rm P_i \mid i \in \nat \}.$
 As to the construction of $\rm L$  
it is$,$ 
at least$,$ constructed from the binary operation $\to.$ Deduction over $\rm 
L$ is defined in the usual  sense.  
Only finitely many steps are 
allowed$,$ and if any axiom schema 
are used$,$ then they do not yield statements  of the form $\rm P_i$ or  
$\rm P_i \to \rm P_j,\ i \not=j.$  Further$,$ deduction from premises is also allowed. 
There are many 
examples of such languages. Propositional languages with denumerably many 
atoms. Indeed$,$ in a predicate language 
with$,$ at least$,$ one predicate the list of all predicates can be 
considered the set of primitives from which $\rm L$ is constructed. Of 
course$,$ simple natural languages are isomorphic to $\rm L$ in the usual 
sense. There will be one modification$,$ however. The modification is in a 
rule of inference. Define the 
$\rm MP_n,\ n \in \nat$ rule of inference on 
$\rm L$ as follows:\parm
{\leftskip=0.5in \rightskip=0.5in \noindent If two previous steps of a 
demonstration (or proof) are of the form $\rm A,\ \rm A \to \rm B$ where  
for each $\rm P_i$ in the primitive expansions of $\rm A,\ \rm A \to \rm B,\  i\leq 
n,$ then the formula $\rm B$ may be written down as the next step. No other 
type of MP rule is used.\par}\parm 
Given a set of hypotheses ${\cal H} \subset \rm L$ and $\rm X  \in \rm L,$ denote by 
the symbol ${\cal H} \rm {\vdash_n X}$ 
 this deductive process. It is immediate that $\rm \vdash_n$ determines a 
finitary consequence operator $\rm C_n$ on $\power {\rm L}.$ Suppose that $\vdash$ 
deduction on $\power {\rm L}$ has all of the above properties with the 
exception that the MP rule of inference is the ordinary modus ponens in 
unrestricted form. Let ${\rm S_0}$ denote the consequence operator determined by 
$\vdash.$\pars
It is a simple matter to show that for any $\rm B \subset L,\ 
{S_0}(B) = \bigcup \{C_n(B) \mid n \in \nat \}.$ Let $\rm B \subset \rm L.$ Suppose 
that $\rm X \in {\rm S_0}(\rm B).$ Then $\rm B \vdash X.$ Now in the formal proof of this 
fact when all of the formula are written in primitive form there is a maximum 
$\rm P_i$ subscript$,$ say $ \rm m \in \nat.$ It follows immediately that 
the same steps yield a formal proof that $\rm B \vdash_m X$ using the 
$\rm MP_m$ in place of any MP step that appears in the $\vdash$ formal 
proof. Thus $\rm X \in \bigcup \{C_n(B) \mid n \in \nat \}.$ Clearly$,$ if 
$\rm B \vdash_n X,$ then $\rm B \vdash X.$ This implies 
$\rm {S_0}(B) = \bigcup \{C_n(B) \mid n \in \nat \}.$ Another interesting 
result is that if $\rm B$ is an ${\rm S_0}$-deductive system$,$ then $\rm {S_0}(B) = 
\bigcup \{C_n(B) \mid n \in \nat \}= B$ implies that for each 
$\rm n \in \nat,\ C_n(B) = B.$ Thus $\rm B$ is also a $\rm C_n$-deductive
system.\pars
Let $\rm n,\ m \in \nat,\ n \leq m,\ B \subset L.$ Suppose that 
$\rm X \in C_n(B).$ If there are any $\rm MP_n$ steps in the formal proof$,$ 
then these steps can also be obtained by application of $\rm MP_m.$ On the 
other hand$,$ if no steps were obtained by the $\rm MP_n$ rule$,$ then the exact 
same steps yield a formal proof that $\rm B \vdash_m X.$ From this we have that for 
each $\rm B \subset L,\ C_n(B) \subset C_m(B).$ Hence$,$ $\rm C_n \leq 
C_m.$ Therefore$,$ $\rm \{C_n \mid n \in \nat\}$ is a chain of consequence 
operators. \pars
Now to show that this chain is of type $\rm K_\infty.$ Let $\rm n < m$ and let
$\rm B = \{P_n, P_n \to P_m \}.$ First$,$ no member of $\rm B$ can be 
obtained as an instance of an axiom. 
Further$,$ it cannot be the case that $\rm B \vdash_n P_m$ for $\rm MP_n$ 
does not apply to $\rm P_n \to P_m$ or$,$ indeed$,$ any formula 
containing $\rm P_m.$ Therefore$,$ $\rm P_m \notin C_n(B).$ Obviously,
$\rm P_m \in C_m(B).$ Hence$,$ $\rm C_n(B)  \not= C_m(B)$ implies that
$\rm C_n < C_m.$ Thus  this chain is of type $\rm K_\infty.$ Further$,$ note that 
$\rm P_m \in {S_0}(B).$ Thus for each $\rm C \in \{C_n \mid n \in \nat\}$ 
there exists some $\rm B 
\subset L$ such that $\rm C(B) \not= {S_0}(B)$ and$,$ clearly $\rm C \leq {S_0}.$ 
Hence$,$ in general$,$ $\rm C < {S_0}$ for all $\rm C \in \{C_n \mid n \in \nat\}.$\parm
{\bf Theorem 5.4.1} {\sl Let $\rm L$ and $\rm K_\infty =  \{C_n \mid n \in 
\nat \}$ be defined as above. Then there exists a purely subtle consequence 
operator $C_\infty \in \Hyper {\bf K_\infty}$ and a *-finite
$F\in \Hyper {(F(\bf L))}$ such that for each $\rm B \subset \rm L$ and each 
$\rm C\in K_\infty$ \pars
{\rm (i)} $\b C(\b B) \subset C_\infty(F),$\pars
{\rm (ii)} ${\bf S_0}(\b B) \subset C_\infty (F) \subset \Hyper {{\bf S_0}(F)} \subset
\Hyper {{\bf S_0}}(\Hyper {\b L}),$\pars
{\rm (iii)} $ \Hyper {{\bf S_0}}(F) \cap \b L = C_\infty(F) \cap \b L.$} \pars 
Proof. (i) is but Corollary 5.3.3.1. From (i)$,$ it follows that 
$\bigcup \{^\sigma(\b C(\b B)) \mid C \in K_\infty \} =  
\bigcup \{(\b C(\b B)) \mid C \in K_\infty \} = {\bf S_0}(\b B) = {^\sigma ({\bf 
S_0}(\b B))} \subset C_\infty(F)$ and the first part of (ii) holds. By 
*-transfer $C_\infty < \Hyper {{\bf S_0}}$ and $C_\infty$ and $\Hyper {{\bf 
S_0}}$ are 
defined on all internal subsets of $\Hyper {\b L}.$ Hence,
$C_\infty (F) \subset \Hyper {{\bf S_0}}(F) \subset \Hyper {{\bf S_0}}(\Hyper {\b L})$ 
and this completes (ii). (iii) follows immediately from (ii) and this 
completes the proof. \qed  
For this application$,$ let $\rm L$ be a predicate type language and $\rm M$ any 
set-theoretic structure in which the predicates and constants are interpreted 
in the usual manner. A finite consequence operator defined on 
$\power {\rm L}$ is {{\it sound}} for $\rm M$ if whenever $\rm B \in \power {\rm L}$ 
has the property that $\rm M \models B,$ then $\rm M \models C(B).$ As usual$,$ 
$\rm T(M) = \{x \mid (x \in L) \land (M \models x) \}.$ Obviously,
if $\rm C$ is sound for $\rm M,$ then $\rm T(M)$ is a C-deductive system. \pars
Corollary  5.3.2.2 implies that there exists *-finite $F \in \Hyper {(\bf 
T(M)})$ such that $\Hyper {\b C(F)} \cap \b L = {\bf T(M)}.$ Notice that $F$ 
being *-finite implies that $F$ is *-recursive. Moreover$,$ $F$ is a *-axiom 
system for 
$\Hyper {\b C(F)}$ and we do not lack knowledge about the behavior of $F$ since 
any formal property about $C$ or recursive sets$,$ among others$,$ must hold true 
for $\Hyper {\b C}$ or $F$ when properly interpreted. If $\rm L$ is a 
first-order language with at least one predicate$,$ then its associated consequence 
operator $\rm S_1$ is sound for first-order structures. Theorem 5.4.1 not only 
yields a *-finite $F_1$ but a purely subtle consequence operator $C_1$ such 
that $F_1$ is a *-axiom system for $C_1(F_1)$ and $\Hyper {{\bf S_1}}(F_1).$ In 
this case$,$ we have that $\Hyper {{\bf S_1}}(F_1) \cap \b L = {\bf T(M)}=C_1(F_1) 
\cap \b L.$ As strange as it may appear$,$ by use of internal and external 
objects$,$ the nonstandard logics $\{\Hyper {\b C}, \Hyper {\b L}\},\ \{C_1,
\Hyper {\b L}\}, \ \{\Hyper {{\bf S_1}}, \Hyper {\b L}\}$ technically by-pass
a portion of G\"odel's first incompleteness theorem. Of course$,$ this 
incompleteness theorem still holds under an internal interpretation. \pars
By definition $\rm b \in {\rm S_0}(\rm B), \rm B \subset \rm L$ iff there is a finite 
length proof of $\rm b$ from the premises $\rm B.$ Thus for each $b \in 
\Hyper {\bf (T(M))}$ there exists a *-finite length proof of b from the 
*-finite 
$F_1.$ If we let $\Hyper {\cal M}$ be an enlargement with the 
$\aleph_1$-isomorphism property$,$ among others$,$ then each *-finite length proof 
is either externally finite or externally infinite. Further$,$ all externally 
infinite proof lengths would be of the same cardinality. \pars
$\{$Remark: Using the customary notation in this chapter$,$ the relation 
$\leq$ has not been starred in $\hyper {\cal N}.$ If this omission is 
confusing$,$ the * can be easily inserted. When these two different order 
relations are compared$,$ the * notation becomes necessary. For example,
the relation $\Hyper \leq$ in $\hyper {\cal N}$ is NOT an extension$,$ in the 
usual sense$,$ of the 
relation $\leq$ as defined in $\cal N,$ although it is an extension of 
$^{\sigma}\leq.$  Also notice that if we had restricted our attention 
to ${\cal C}_f,$  then the partial order $\leq$ is characterized totally 
by the finite subsets of A. This is useful  since that $\Hyper \leq$ is 
characterized by the *-finite subsets of $\hyper {\rm A}.$ It's clear that 
our concept of a consequence-type operator must be generalized slightly. 
Let $\cal B$ and 
${\cal B}_0$ be two families of sets. Then if $f\colon {\cal B} \to {\cal 
B}_0$ satisfies axioms (2)(3)(4) or (2)(3)(5) or the *-transform of these 
axiom systems$,$ then $f$ is a subtle consequence operator.
I also point out 
that$,$ unfortunately$,$ there are many typographical errors in 
reference [24].$\}$\par
\vfil
\eject
\centerline{NOTES}

\vfil\eject
\centerline{\bf 6. ASSOCIATED MATERIAL}\par
\bigskip
\leftline{\bf 6.1 Perception}\par 
\medskip
In this section$,$ the theory of ultralogics is applied to one aspect of 
 {{\it subliminal perception.}} What is needed is an interpretation scheme. 
When subsets of $\hyper {\cal E}$ are concerned the {{\it conscious objects}} 
are 
subsets [resp. elements of] $^\sigma {\cal E} = \cal E.$ The {{\it subconscious 
objects}} are nonstandard internal subsets [resp. elements of]
$\hyper {\cal E}.$ Moreover$,$ subconscious objects can contain conscious 
objects and the union of a subconscious set and a finite  conscious set is a 
subconscious set. The {{\it unconscious objects}} are external nonstandard 
subsets of $\hyper {\cal E}.$ Like definitions apply to members of $\hyper 
{\cal E} \times \hyper {\cal E}$ and so forth. In what follows$,$ only strong 
reasoning from the perfect is considered. You may assume that it is defined 
on a natural language  ``very,$\sp$'' or a formal language
``$\rm V \land$'' and the like. \pars
As to some sort of interpretation procedure the following seems adequate. Let 
$\lceil \ \rceil $ denote an {interpretation symbol.} First$,$ we have 
{{\it subperception and the better than ordering.}} Let $A \subset 
\hyper {\cal E}$ be one of the above defined objects in the domain of $
\Hyper {\b \Pi}$. Let internal $D \subset \Hyper 
{\b \Pi}(A)$ and standard $^\sigma E \subset
\Hyper {\b \Pi}(A)$. Assume that $^\sigma E 
\leq_B D$ and that each member of $i^{-1}[E]$ is a sentence which is 
distinctly comparable by the  ``very,$\sp$'' symbol string. [Note I am not 
differentiating between the object and a constant representing that object.] 
One might interpret the following: $\lceil \forall x((x \in {^\sigma E}) \to 
\exists y(( y \in D)\land (x \leq_B y)))\rceil:$ ``You are (I am$,$ we are$,$ etc)
subperceptibly aware that for each conscious (known) object (element$,$ member)
of (in) $\lceil ^\sigma E \rceil$ there exists an object (element$,$ member) of 
(in)
$\lceil D \rceil$ which is better than that conscious object of (in) $\lceil 
^\sigma E \rceil.$ ''\pars
Let $a \in {^\sigma E}.$ Then: $\lceil \exists y ((y \in D) \land (a \leq_B 
y))\rceil:$  ``You are (I am$,$ we are) subperceptibly aware that there exists a 
conscious (known) object (element$,$ member) of (in)
$\lceil D \rceil$ which is better than $\lceil a \rceil$.'' Note that
$a = [f_0] \in {\bf BP}$ and $\lceil a \rceil = i^{-1}(f_0(0)).$ Another 
example is $\lceil ^\sigma E \leq_B D \rceil:$ ``You are (I am$,$ we 
are) subperceptibly aware that 
$\lceil D \rceil$ is better than $\lceil ^\sigma E\rceil.$ ''\pars
For another example$,$ let $^\sigma F \subset \Hyper {\b\Pi}(A).$ Then 
$\lceil (\vert {^\sigma E} \vert < \vert {^\sigma F} \vert ) \land
(^\sigma E \leq_B D)\land (^\sigma F \leq_B D)\rceil:$  ``The result that 
$\lceil D \rceil$ is better than $\lceil ^\sigma F \rceil  $ is a stronger 
subperceptible property than $\lceil D \rceil$ is better than $\lceil ^\sigma 
E \rceil.$''\pars
We also have the idea of {{\it general subperception}}. In this case$,$ we use 
some of the meaningful set-theoretic terminology. Let internal nonstandard 
$A,\ B \subset
\Hyper {\bf BP}.$ Now any elementary set-theoretic relation 
\noindent existing between 
$A$ and 
$B$ can be subperceptibly interpreted as $\lceil A \subset B \rceil:$ 
``You are (I am$,$  we are) subperceptible aware of the following: $\lceil A 
\rceil$ is contained in  $\lceil B \rceil.$ '' Also you might interpret 
relations between standard objects as a complete awareness.\par
\bigskip
\leftline{\bf 6.2 Existence}
\medskip 
Some philosophers of science differentiate between theoretical entities and 
those\break
\vfil
\centerline{}
\eject
\noindent that are assumed to exist in objective reality. In the original work in 
ultralogics$,$ these two concepts were disjointly modeled. This was 
done as follows: Consider $^\sigma [f_0] = [f_0] \in \hyper  {\cal E}$ to be the 
unique partial sequence with the property that
$i^{-1} (f_0(0)) = {\rm externally\sp exists\sp in\sp reality} = 
\lceil [f_0] \rceil.$  For $A \subset \Hyper {\bf BP},$ let $(A)_R = \{(x, [f_0]) 
\mid x \in A \}.$ Then for any $E \subset \Hyper {\bf BP},$ define $RR(E,A) = \{(x,y)\mid 
(x \in E)\land (y \in (A)_R\}$ to be the {{\it realism relation}}.
These definitions are then 
extended to $\Hyper {\bf BP} \times \Hyper {\bf BP}$ in such a manner
that $(A)_R \times (B)_R$ is considered to be isomorphic to $(A \times B)_R.$ 
I now believe that this is a waste of effort. The  difference lies in the 
interpretation and not in the mathematical structure. Thus$,$ under the 
interpretation$,$ if one wishes to differentiate between these two concepts$,$ one 
simply includes ``existence in objective reality'' as a 
part of the interpretation for some entities and the statement  ``theoretical 
entities'' for other distinct entities. \par
\bigskip
\leftline{\bf 6.3 An Alternate Approach}
\medskip
What is presented in this section is mainly of historical interest although 
this author's first research into nonstandard analysis used this alternate 
approach. This approach utilizes a pseudo-set theory and has essentially been 
replaced by the superstructure approach. Some years ago$,$ certain applications 
employ this 
alternate approach due to its use of a basic language that is somewhat more expressive 
than the $\in, =$ language. However$,$ what might be gained in an additional 
freedom of expression will lead to a more complex array of extensions$,$ 
definitions and the requirement that extreme care  be exercised.\pars

All of our constructions are within $\bf ZFC.$ We utilize the {{\it 
transitive closure operator,}} denoted by $TC.$ 
Let $V$ be a set. (Note: This definition also applies to atoms and sets containing atoms.) The 
transitive closure of $V$ is obtained by an inductive construction using 
the union operation.  Let $V_0 = V$ and for each $i \in \nat,$
 let $V_{i+1} = \cup V_i.$ Then the transitive set $TC(V)  = \bigcup \{ V_i \mid i \in \nat \}.$ 
The set $TC(V)$ for a set $V$ has the property that if $W$ is another 
transitive set such that $V\subset W,$ then $V \subset TC(V) \subset W.$ 
Define the {{\it superstructure 
operator}}$,$ denoted by $SS,$ on $V$ as $SS(V) = \bigcup \{ U_i 
\mid i \in \nat \},$ where $U_0 = V, \ U_{i+1} = U_i \cup \power {U_i}, i 
\in \nat.$ Recall that this is the first type of superstructure defined in 
Chapter 2. To correspond to our previous investigation$,$ let $V ={\cal W} \cup \nat$ and 
${\cal N}_1 = SS(TC(V)).$ (If $V$ is a set of atoms, then $TC(V) = V.$) 
Let the structure ${\cal M}_1 = \langle {\cal N}_1, \in, =, {\rm ap}, {\rm pr} 
\rangle,$ where $\in,\  =$ are 
the usual set-theoretic membership and set equality relations restricted to 
${\cal N}_1$ and ap$,$ pr are two ternary relations$,$ the  ``applying a function 
to its argument'' and  ``ordered pair creation''$,$ respectively. (Of course$,$ =$,$ 
ap$,$ pr can all be defined in terms of $\in.$ ) Notice that ${\cal M}_1$ is a 
fragment of our {\bf ZFC} model. \pars
Consider a $\kappa$-adequate ultrafilter$,$ where $\kappa > \vert {\cal N}_1 
\vert.$ By Theorem 7.5.2 in [19] or 1.5.1 in [9]$,$ such an ultrafilter $\cal U$ 
exists in our {\bf ZFC} model and is determined by the indexing set $J = 
F(\power \kappa).$ By the ultrapower or ultralimit construction$,$  a first-order
structure ${\cal M}_2 = ({\cal N}^J, \in_{\cal U},=_{\cal U}, {\rm ap}_{\cal 
U},{\rm pr}_{\cal U})$ is obtained of the same type as is ${\cal M}_1$ but 
${\cal M}_2$ is a nonstandard model for the set of all sentences$,$ $K_0,$ in 
our first-order language L$,$ with predicates $\in,=,$ ap$,$ pr which hold in
${\cal M}_1.$ (Theorem 
3.8.3 in [19]) Note that the cardinality of 
the set of constants in L $\geq \vert {\cal M}_1 \vert.$  Further$,$ the members 
in ${\cal N}^J$ are interpreted by constants in an extended language ${\rm 
L}^\prime$. By 
the axioms of our {\bf ZFC} set-theory$,$ the relation  ``='' is an equivalence 
relation with substitution for $\in,$ ap$,$ pr and$,$ hence$,$ $=_{\cal U}$  
has these properties for $\in_{\cal U}, {\rm ap}_{\cal U}, {\rm pr}_{\cal U}.$ 
Consequently$,$ we shift the structure ${\cal M}_2$ (i. e. contract it) [13$,$ p. 
83] and obtain a structure $^\prime {\cal M} = \langle ^\prime{\cal N}_1, 
\eps, 
=, {^\prime{\rm ap}},^\prime{\rm pr} \rangle,$ where = is the original equality 
in our {\bf ZFC} model. Note that members of $^\prime {\cal N}$ are still 
interpreted by constants in ${\rm L}^\prime$ as before.
[In [4] and [11]$,$ a structure isomorphic to $^\prime {\cal M}$ is obtained by 
application of the compactness theorem for a first-order language.] \pars
We next isomorphicly embed ${\cal M}_1$ into $^\prime{\cal M}$ in the same 
manner as outlined in [4$,$ p. 22]. However$,$ please notice that the following 
notation differs from that used in this reference. First$,$ let $I$ be the 
the original interpretation map from L onto ${\cal M}_1$ and let $^\prime I$ 
by the composition of the extended ultrapower interpretation map and the 
contraction interpretation map restricted to the constants in L.\pars
 Now for each 
``a'' that is a constant in L$,$ define $^\sigma (I(a)) = {^\prime I(a)} \in 
{^\prime {\cal N}}.$ Assume that ${\rm a,\ b}$ are constants in L that 
represent the same element  and consider the well-formed formula $(\rm a = \rm 
b).$ Then ${\cal M}_1 \models (\rm a = \rm b)$ iff $^\prime {\cal M} \models (\rm 
a= \rm b)$ implies that $^\sigma (I(a)) = {^\sigma (I(b))}.$ Thus the map 
$\sigma$ is well-defined. Again let a$,$ b be constants in L. 
Then  ${\cal M}_1 \models \neg(\rm a = \rm b)$ iff $^\prime {\cal M} \models 
\neg(\rm a = \rm b)$ implies that $^\sigma(I(\rm a)) \not= {^\sigma (I(\rm b))}$ 
iff 
$I(\rm a) \not= I(\rm b).$ Thus $\sigma$ is injective. It is immediately clear 
that ${\cal M}_1 \models (\rm a \in \rm b)$ iff $^\prime {\cal M} \models
(\rm a \in \rm b)$ implies that $I(\rm a) \in I(\rm b)$ iff $^\sigma(I(\rm a)) 
\ {\eps}\  ^\sigma(I(\rm b))$ and$,$ in like manner$,$ for the relations   ``ap'' and   
``pr''. Consequently$,$ $\sigma$ is an isomorphic embedding of ${\cal M}_1$ into 
$^\prime{\cal M}.$ For convenience in all that follows$,$ we suppress the 
interpretation map notation and simply use the constants of the language L  
(and the extended language $\rm L^\prime$) to represent members in ${\cal N} =
\{x \mid \exists \rm a(\rm a\in \rm L) \land (^\sigma(I(\rm a)) = {^\prime I(\rm 
a))} \land ( x = {^\prime I(\rm a))} \}.$ Of course$,$ ${\cal N} \subset
{^\prime {\cal N}}$ and ${\cal M}_1$ is isomorphic to ${\cal M} = \langle 
{\cal N},{(\eps})/\sigma, =, {( ^\prime {\rm ap})}/\sigma, ({^\prime {\rm 
pr }})/\sigma \rangle.$ [Note: $(\eps)/\sigma$ is the relation 
restricted to members of ${\cal N}$$,$ etc.]\pars
Now to continue this construction. First$,$ by Theorem 1.5.2 in [9]$,$ 
$^\prime{\cal M}$ is an enlargement. For each $p \in {^\prime {\cal N}},$ 
let $\hyper p = \{x \mid (x \in {^\prime {\cal N}})\land (x\ \eps\ p) 
\}$ and $\hyper {\cal N} = \{\hyper p 
\mid p \in {^\prime {\cal N}} \}.$ Notationally$,$ let $\hyper{\hyper {\cal N}} 
= {^\prime {\cal N}}\cup {\cal N}_1 \cup \hyper {\cal N}$ 
and define $^0{\cal N} = SS(TC(\hyper{\hyper{\cal N}})).$ Obviously,
$ {\cal N} \cup {\cal N}_1 \cup\  {^\prime {\cal N}}\cup \hyper {\cal N} \cup
\hyper{\hyper {\cal N}} \subset TC(\hyper{\hyper{\cal N}}).$ Since 
$TC(\hyper{\hyper{\cal N}})$ is a member of $SS(TC(\hyper{\hyper{\cal N}}))$$,$ 
then this implies$,$ by Theorem 2.6 (iii) [4]$,$ that ${\cal N}, {\cal N}_1, 
{^\prime {\cal N}}, \hyper {\cal N}, \hyper {\hyper {\cal N}} \in {^0{\cal 
N}}.$ Finally$,$ applying Theorem 2.6 (ii) and 2.10 in [4] to the transitive set 
$^0{\cal N},$ one obtains that $^0{\cal N }$ is closed under finite power set 
iteration and finite Cartesian products. The structure $^0{\cal M} = \langle 
{^0{\cal N}}, \in, =, {\rm ap},  {\rm pr} \rangle$ that is a fragment of our 
{\bf ZFC} model is the G-structure in this alternative approach. \pars
The model $^0{\cal M}$ contains all of the set-theoretic objects needed for 
this investigation. The fact that we are only interested in semantic 
consistency allows us to consider all of the structure $\cal M$ as the 
standard model in which sentences from L are interpreted and $^\prime{\cal M}$ 
the nonstandard model for sentences from L. Of course$,$ we can always return to 
${\cal M}_1$ by application of $\sigma^{-1}.$ \pars
There is one important notational convention that is continually employed.
The $\sigma$ map is suppressed when considering members of ${\cal N}.$ That is 
to say that for each constant  $\rm a  \in L,$ $^\sigma (I(\rm a)) = {^\prime I(\rm 
a)} = \rm a.$ The use of $\cal M$ can be made more efficient since there should be 
no great difficulty if you consider $(\eps)/\sigma,(^\prime{\rm pr})/\sigma, 
(^\prime{\rm ap})/\sigma$ to be the same as the {\bf ZFC} model relations 
$\in, {\rm ap}, {\rm pr}$ for there is no first-order differences between 
these structures. \pars
The fact that we are actually working with the restricted $\eps,\ 
^\prime{\rm ap},\ ^\prime{\rm pr}$ can be determined by the additional result 
that the only objects to which the restriction of these relations apply 
are members of $\cal N.$ The only other objects to which nonrestricted  
$\eps,\ ^\prime{\rm ap},\ ^\prime{\rm pr}$ apply are elements 
of $^\prime{\cal N}$  
that are not members of $\cal N.$ The actual $\in, {\rm ap}, {\rm pr}$ 
are used in all other contexts such as the following important definition as 
previously stated. For each $A \in {^\prime {\cal N}},\ \hyper A = 
\{ (x \in {^\prime {\cal N}})\land (x\ {\eps\ A)}\}.$ This is the 
beginning of certain technical features for this model. What is significant
as we define some of these technical terms is that all of the objects within 
this and other nonstandard investigations are set-theoretic members of 
$^0{\cal M}$ and$,$ of course$,$ $^0{\cal M}$ is in our {\bf ZFC} model.\pars
Rather than force the reader to seek out references [4] or [11]$,$ I reproduce 
here the more significant definitions required to relate many of our 
results to the $\in, {\rm ap}, {\rm pr}$ operations within the structure 
$^0{\cal M}.$ Let $p \in \ {^\prime {\cal N}}.$ If $ p\notin {\cal N},$ then 
$p$ is called a {{\it nonstandard object or entity}}. If $p \in \cal N,$
then $p$ is called a {{\it standard object}}. If $S \subset \hyper P$ and 
there 
exists some $Q \in \ {^\prime{\cal N}}$ such that $S = \Hyper Q,$ then $S$ is 
called an {{\it internal object or set}}. Observe that $P\in {^\prime{\cal 
N}}$ that is not a $^\prime$atom is an internal subset of itself. Also it is 
often the case that each element of $\hyper P$ is called {\it internal} for if 
$P \in {^\prime {\cal N}},$ then there exists some $X_n$ such that $P\    
{\eps}\ X_n $ and $X_n$ is $^\prime$-transitive. Thus if $p\  
{\eps}\ P,$ then $p\ {\eps}\ X_n$ implies that $p \in\  {^\prime{\cal 
N}}.$ Intuitively internal means that there exists a symbolic name in 
L$^\prime$ for the object that generates$,$ under the given definitions$,$ the 
second corresponding object.\pars
This generation of the second corresponding object is of a special nature. Let 
$f \in {^\prime{\cal N}}$ be an $^\prime$n-ary relation where $ n > 1.$ Thus 
$f$ satisfies in L$^\prime$ the appropriate sentence that defines such a 
object. Extend $f$ in the following manner. Let $f^\star = 
\{(a_1,\ldots,a_n)\mid(a_1 \in {^\prime {\cal N }}) \land \cdots \land
(a_n \in {^\prime{\cal N}}) \land ((a_1 {^\prime{\rm pr}}\cdots {^\prime{\rm 
pr}} \ a_n)\ {\eps}\ f)\}.$ In general$,$ $\hyper f \not= f^\star.$  For the 
many properties associated with this definition$,$ refer to references [4] [11].
I note that
in [4] one of the important properties for such an extension of $f$
relative to the i'th projection [Theorem 4.5 (vii)] is stated on one side of 
the equation incorrectly. However$,$ the proof goes through correctly and one 
should correct the statement of that small portion of the theorem to show that 
the i'th projection of the n-ary relation $f^\star = \Hyper ($of the i'th 
projection of $f$). [Note: there are two theorems in [4] that are proved 
incorrectly$,$ even though the theorem statement is correct. The proofs were 
corrected when these results were published.]  Any n-ary relation that is 
produced by an extension that has the $^\star$ on the right is called 
{{\it an internal n-ary relation.}} Notice that what this actually means is 
that there is a name for the $^\prime$n-ary relation in the extended 
language $\rm L^\prime.$ \pars
Our major interest and application for this model will be confined to objects 
in
$\cal E$ as well as in $\hyper {\cal E} -{\cal E},$ and a fixed power set 
iteration  
or 
Cartesian products of these objects. The use of the *-ing process is 
different in this model than it is in the model utilized in the previous 
sections of this chapter and previous chapters  
of this book. For example$,$ it is important to realize that the set $f \in \cal 
E$ is a finite set of functions. Thus $\forall x(x \in {[f]} \iff (x = a_1) \lor
\ldots \lor (x = a_n))$ holds in $\cal M.$ Therefore$,$ $\forall x(x \in \Hyper 
[f] \iff (x = a_1) \lor
\ldots \lor (x = a_n))$ implies that $\Hyper {[f]} =[f]$ and 
$[f]$ is internal. Observe that the symbols $a_1,\ldots,a_n$ do not carry the 
* notation  as would be necessary$,$ prior to our identification process$,$ 
in the previous sections of this chapter and previous chapters. Further$,$ each 
$ g \in [f]$ is a finite set of ordered pairs$,$ as previously. Thus $\hyper g = 
g$ and $g$ is internal. \pars 
Even though the above property seems to be a nice property$,$ 
the nonstarring of standard objects$,$ it turns out that the partition concept 
must be handled differently. Indeed$,$ Theorem 3.2.3 is not true in this model. 
If $A,\ B \in \cal N$ and $B$ is a partition of $A,$ then $\Hyper B$ is not a 
partition of $\hyper A.$ What is needed is to consider the set $D = \{\hyper x 
\mid x \in \Hyper B \}.$ The $D$ is a partition for $\hyper A,$ but $D$ is not 
in general an internal set. On the other hand$,$ each element of $\Hyper B$ is 
an internal set as is each finite subset. It is interesting to note that we 
require a different extension definition for the consequence operators when 
$^0{\cal M}.$\pars
Let $C\colon \power A \to \power B$ be a standard set-valued map. One must be 
more careful with the extensions of such set-valued maps than the other maps
since the types of objects contained in the ordered pairs are of 
significance. Observe that $C^\star$ is composed of ordinary ordered pairs of 
$^\prime$sets and as such these sets are in $^0{\cal M}$ and contain 
$^\prime$elements. However$,$ these sets may also contain ordinary elements as 
well. Assume that $(P,Q) \in C^\star,$ and that $P,\ Q$ are not finite 
standard sets$,$ then 
$(P,Q) \not= (\hyper P, \hyper Q).$ Consequently$,$ in general$,$ a map such as 
$C$ can be extended to a map $\hyper {\underline{C}}=\{(\hyper P, \hyper Q)
\mid (P,Q) \in C^\star \}.$ In this case$,$ $\hyper P,\ \hyper Q$ contain only 
$^\prime$elements from $P$ and $Q.$ This different interpretation occurs 
because the  ``starring'' process in the first model used in this analysis is 
distinct from the  ``starring'' process as employed with respect to $^0{\cal 
M}.$ In fact$,$ the * process in the first model is a renaming of the standard 
language objects as  they are  interpreted within that model and is a member of the 
extended language $\rm L^\prime.$ The other members of $\rm L^\prime$ are 
restricted to internal members of our model. \pars
More importantly$,$ with the first model all of the relations $\eps, \ 
^\prime{\rm pr},\ ^\prime{\rm ap}$ have been replaced by the ordinary $\in,\ 
{\rm pr}\ {\rm ap}$ within the {\bf ZFH} model and the standard model has been 
embedded into the structure that what would have been the $\eps,\ ^\prime{\rm pr},\ 
^\prime{\rm ap}$ defined objects denoted by members of L are so altered that 
they become the original relations restricted to entities that are 
isomorphicly related to the original standard objects. One can say that the 
first model alters the objects with a ``minimal'' language change. This 
alternate approach requires a much larger language change but is more expressive 
in character. \pars
The above extension processes lead to three distinct objects 
$\hyper C,\ C^\star$ and $ \hyper {\underline{C}}$ within $^0{\cal M}.$ If 
$C \colon \power A  \to \power B$ is a consequence operator$,$ then 
$C^\star$ and $\hyper {\underline{C}}$ are interesting but distinct member of 
$^0{\cal M}.$ They both satisfy (extended) Tarski type axioms. If we were to 
continue this development$,$ then the map $\hyper {\underline{C}}$ appears to be 
the most appropriate for such an investigation. However$,$ $\hyper 
{\underline{C}}$ is$,$ in general$,$ an external object. It does not satisfy the 
internal defining method for $^\prime$n-ary relations that requires the 
variables to vary over $^\prime{\cal N}.$ Observe that $\hyper 
{\underline{C}}$ is defined for all of the internal subsets (within $^0{\cal 
M}$) of $\hyper A.$ With respect to the first model$,$ $\hyper C$ is restricted 
to its internal subsets of $\hyper A$ as well. Now $C^\star$ is an internal 
map in $^0{\cal M}$ that is defined on  ``internal entities'' that are 
$^\prime$subsets of $A$ and it yields $^\prime$subsets of $B$ which$,$ when 
viewed 
from the structure $^0{\cal M},$ could contain many non-$^\prime$elements. 
Notice that we cannot obtain any information about these other objects by 
simply transferring$,$ by the *-transfer method$,$ information from the standard 
model. These objects could be investigated by a more careful analysis of the 
exact construction of $^0{\cal M}.$ \pars
Finally$,$ this alternate approach also requires a more specific definition for 
the  ``standard restriction process.'' The definition of the  ``standard restriction'' for 
$^0{\cal M}$ would depend upon which type of extension $\hyper C,\ C^\star,\ 
{\rm or}\ \hyper {\underline{C}}$ is used. It is clear that all of this as 
well as the appropriate definitions for human$,$ subtle and purely subtle 
entities can be successfully accomplished.\parm
[NOTE: In Chapters 7 - 11, in most cases, the symbol $A_1 = i[{\cal W}].$ Also in Chapter 9, the structure ${\cal M}_1$ being used is incorrectly denoted. The superstructure $\cal N$ is to have either ${\cal W} \cup \real$ or ${\cal W} \cup {\cal Q}$ as a ground set. It is useful to extend the language to ${\cal W}'$ that includes symbols for $\real$ or $\cal Q.$ Further in Section 9.1, $r \in \Hyper A_1 \simeq \hyperreal$ should read $r \in \Hyper A_1.$ Notice that if one chooses to use ${\cal W}'$, then $r$ corresponds to an $r' \in \Hyper {\cal W}'.$ \pars
\vfil
\eject
\centerline{\bf CHAPTERS 1---6 REFERENCES}
\medskip
\noindent {\bf  1} Barwise J. (ed.) {\it Handbook of Mathematical Logic,} North-Holland,
Amsterdam$,$ 1977.\pars
\noindent {\bf  2} Birkhoff$,$ G. D.$,$ A set of postulates for plane geometry based on 
scale and protractor$,$ Annals of Math. 33(1932)$,$ 329---345.\pars
\noindent {\bf  3} Hamilton$,$ S. G.$,$ {\it Logic for Mathematicians,} Cambridge University
Press$,$ New York$,$ 1978.\pars
\noindent {\bf  4} Herrmann$,$ R. A.$,$ Nonstandard Topology. Ph. D. Dissertation$,$ American 
University$,$ 1973. (University Micro Film \# 73-28,762)\pars
\noindent {\bf  5} Jech$,$ T. J.$,$ {\it Lectures in Set Theory,} Lecture Notes in 
Mathematics Vol.  217$,$ Springer-Verlag$,$ Berlin$,$ 1971.\pars
\noindent {\bf  6} Jech$,$ T. J.$,$ {\it The Axiom of Choice,} North-Holland$,$ Amsterdam$,$ 
1973. \pars
\noindent {\bf  7} Kleene$,$ S. C.$,$ {\it Introduction to Metamathematics,} D. Van Nostrand 
Co.$,$ Princeton$,$ 1950.\pars
\noindent {\bf  8} Kleene$,$ S. C.$,$ {\it Mathematical Logic,} John  Wiley and 
Sons$,$ Inc.,
New York$,$ 1967.\pars
\noindent {\bf  9} Luxemburg$,$ W. A. J.$,$ A general theory of monads. in {\it Applications of 
Model Theory to Algebra$,$ Analysis and Probability,} (ed. Luxemburg)$,$ Holt$,$ 
Rinehart and Winston$,$ New York$,$ (1969)$,$ 18---86.\pars                        
\noindent {\bf  10} Luxemburg$,$ W. A. J.$,$ What is nonstandard analysis$,$ in {\it Papers in 
the Foundations of Mathematics} No. 13 Slaught Memorial Papers$,$ Amer. Math. 
Monthly$,$ (June-July 1973)$,$ 38---67.\pars
\noindent {\bf  11} Machover$,$ M. and J. Hirschfeld$,$ {\it Lectures on Non-standard 
Analysis,} Lecture Notes in Mathematics$,$ Vol. 94$,$ Springer-Verlag$,$ Berlin$,$ 
1969.\pars
\noindent {\bf  12} Markov$,$ A. A.$,$ Theory of algorithms$,$ Amer. Math. Soc. Transl.$,$ 
Ser. 2$,$ 15(1960)$,$ 1---14.\pars
\noindent {\bf  13} Mendelson$,$ E.$,$ {\it Introduction to Mathematical Logic}$,$ 2'nd ed.,
D. Van Nostrand Co.$,$ New York$,$ 1979.\pars
\noindent {\bf  14} Rasiowa$,$ H. and R. Sikorski$,$ {\it The Mathematics of 
Metamathematics,} Polska Akademia Nauk$,$ Monografie Methematyczne$,$ Tom 41$,$ 
Warsaw$,$ 1963. \pars
\noindent {\bf  15} Robinson$,$ A.$,$ On languages which are based on non-standard arithmetic,
Nagoya Math.$,$ 22(1962)$,$ 83---118.\pars
\noindent {\bf  16} Robinson A.$,$ {\it Non-Standard Analysis}$,$ (2'nd ed.) North-Holland$,$ 
1974.\pars
\noindent {\bf  17} Robinson$,$ A.$,$ and E. Zakon$,$ A set-theoretic characterization of 
enlargements$,$ in  {\it Applications of 
Model Theory to Algebra$,$ Analysis and Probability,} (ed. Luxemburg)$,$ Holt$,$ 
Rinehart and Winston$,$ New York$,$ (1969),109---122.\pars    
\noindent {\bf  18} Stoll$,$ R.$,$ {\it Set Theory and Logic,} W. H. Freeman and Co.$,$ San 
Francisco$,$ 1963. \pars
\noindent {\bf  19} Stroyan$,$ K. D. and W. A. J. Luxemburg$,$ {\it Introduction to the 
Theory of Infinitesimals,} Academic Press$,$ New York$,$ 1976.\pars
\noindent {\bf  20} Suppes$,$ P. {\it Axiomatic Set Theory,} Van Nostrand$,$ 1960$,$ (Reprint
Dover$,$ 1972.)\pars
\noindent {\bf  21} Tarski$,$ A.$,$ {\it Logic$,$ Semantics$,$ Metamathematics,} (Papers from 
1923---1938)$,$ Oxford University Press$,$ New York$,$ 1956.\pars
\noindent {\bf  22} Thue$,$ A. Probleme \"uber Ver\"anderunger von Zeichenreihen nach
gegebenen Regeln$,$ Skrifter utgit av Videnskapsselskapet i Kristiania$,$ I. 
Metematisk---naturvidenskabelig Klasse 1914.10. \parm
\centerline{Additional References}
\medskip
\noindent
\noindent {\bf  23} Hurd$,$ A. E. and P. A. Loeb$,$ {\it An Introduction to Nonstandard 
Analysis,} Academic Press$,$ Orlando$,$ 1985.\pars
\noindent {\bf  24} Herrmann$,$ R. A.$,$ Nonstandard consequence operators$,$ Kobe J. Math. 
4(1987)$,$ 1---14.\pars
 \end